\documentclass[twoside,british,3p]{elsarticle}
\usepackage{amsmath}
\usepackage{amsthm}
\usepackage{amssymb}
\usepackage{fontspec}
\usepackage{array}
\usepackage{booktabs}
\usepackage{mathrsfs}
\usepackage{stmaryrd}
\usepackage{graphicx}
\usepackage{geometry}
\usepackage[bookmarks=true,bookmarksnumbered=false,bookmarksopen=false,
 breaklinks=false,pdfborder={0 0 0},pdfborderstyle={},backref=false,colorlinks=false]
 {hyperref}
\hypersetup{pdftitle={An entropy-stable splitting strategy for a hyperbolic dispersive shallow-water model on unstructured meshes}}

\makeatletter

\providecommand{\tabularnewline}{\\}

\makeatother

\providecommand{\widecheck}[1]{\check{#1}}

\journal{Journal of Computational Physics}

\makeatother

\providecommand\lemmaname{Lemma}
\theoremstyle{plain}
\newtheorem{lem}{\protect\lemmaname}
\providecommand\propositionname{Proposition}
\newtheorem{prop}{\protect\propositionname}
\providecommand\remarkname{Remark}
\theoremstyle{remark}
\newtheorem{rem}{\protect\remarkname}
\usepackage{polyglossia}
\setdefaultlanguage[variant=british]{english}
\begin{document}
\begin{frontmatter}
\title{An entropy-stable splitting strategy for a hyperbolic dispersive shallow-water
model on unstructured meshes}
\author[imt-insa-cnrs]{Frédéric Couderc}
\ead{couderc@math.univ-toulouse.fr}
\author[icj,iuf]{Arnaud Duran\corref{cor1}}
\ead{arnaud.duran@univ-lyon1.fr}
\author[lama]{Maria Kazakova}
\ead{maria.kazakova@univ-smb.fr}
\author[shom]{Bastien Lombard}
\ead{bastien.lombard@shom.fr}
\author[imt-insa-cnrs]{Pascal Noble}
\ead{noble@insa-toulouse.fr}
\author[imt-insa-cnrs,shom]{Rémy Baraille}
\ead{remy.baraille@shom.fr}
\author[imt-insa-cnrs]{Jean-Paul Vila}
\ead{vila@insa-toulouse.fr}
\address[imt-insa-cnrs]{Univ Toulouse, INSA Toulouse, CNRS, IMT, Toulouse, France.}
\address[icj]{Université Lyon 1, Centrale Lyon, INSA Lyon, Université Jean Monnet,
CNRS, ICJ UMR5208, 69622 Villeurbanne, France.}
\address[iuf]{Institut Universitaire de France (IUF).}
\address[lama]{LAMA, UMR 5127, Université Savoie Mont Blanc, France.}
\address[shom]{SHOM - Service Hydrographique et Océanographique de la Marine, France.}
\begin{abstract}
This paper presents an entropy-stable splitting strategy for the numerical
solution of a hyperbolic approximation of the Serre--Green--Naghdi
(SGN) model, leading to an efficient solver for dispersive water waves
on unstructured meshes. The model was introduced by Escalante et al.
(2019) and justified shortly thereafter by Richard (2021) as an asymptotic
model of the free-surface Euler equations under a weak compressibility
assumption. While the parent model includes full topography terms
and enhanced dispersive properties, the present study focuses on its
mild-slope form without enhanced dispersion, which retains an exact
energy equation, the cornerstone of the stability analysis. The splitting
yields a shallow-water subsystem with two passive scalars, governed
by the \textit{slow} gravity waves, and an acoustic subsystem carrying
the dispersive effects, governed by the \textit{fast} acoustic waves.
The gravity step adapts an existing entropy-stable scheme, while an
original scheme is derived for the acoustic step; both dissipate the
energy at the discrete level, under local fully discrete stability
criteria on general unstructured meshes. Complementary linear analyses
select the low-dissipation calibration of the gravity step used in
the computations, outside the scope of the nonlinear estimate. The
acoustic step, which requires no global linear solve, is sub-iterated
after each gravity time step; the resulting solver is three to six
times more expensive than a classical shallow-water solver, depending
on the chosen acoustic speed. Numerical validations, ranging from
academic benchmarks to a realistic coastal configuration, highlight
the properties of the approach and pave the way for operational use.
\end{abstract}
\begin{keyword}
Nonlinear dispersive waves \sep Serre--Green--Naghdi equations
\sep Hyperbolic relaxation \sep Splitting schemes \sep Entropy-stable
schemes \sep Coastal flows\cortext[cor1]{Corresponding author}
\end{keyword}
\end{frontmatter}

\section{Introduction}\label{sec:introduction}

Fast and accurate numerical simulation of dispersive water waves remains
a very challenging problem, motivated simultaneously by academic,
safety, and economic issues. Climate change is leading to rising sea
levels and increasing frequency and intensity of storm surges, the
combination of which considerably increases coastal risks \citep{oppenheimer2019,vousdoukas2018}.
Both science and engineering must rely on advanced numerical tools
to enable real-time forecasting whenever possible, and this applies
both to large-scale simulations and to local high-resolution zooms.
The stakes are high and involve major theoretical and technical challenges.
This notably implies the ability to accurately and efficiently capture
nonlinear and dispersive effects, and, in some cases, turbulent ones,
as well as to account for phenomena that remain difficult to reproduce
precisely, such as coastal flooding or infra-gravity waves. For completeness,
a comprehensive description of these mechanisms also requires accounting
for wave-breaking phenomena and their effects on coastal regions,
considering large domains with short-wavelength waves over long times.
Another important aspect concerns the design optimization of protective
coastal structures, with high-fidelity simulations of very local areas
over long times to cover the full range of scenarios.

Water waves involve three-dimensional and turbulent free-surface flows
with complex geometries, making a direct numerical simulation based
on the fundamental equations of fluid mechanics prohibitively expensive,
even on modern CPU or GPU architectures. For this reason, the nonlinear
shallow-water equations, also known as the Saint-Venant equations,
are commonly used for hydraulic or tidal wave simulations. They are
obtained by depth-averaging the fundamental equations of fluid mechanics
over the water column under the assumption of \textit{long} wavelengths
compared to the water depth, and significantly reduce the computational
cost, the system having one dimension less. Nevertheless, waves induced
by wind forcing, together with bathymetry effects encountered in realistic
situations, generate wavelengths that are much shorter than tidal
ones. This means that the classical shallow-water equations are no
longer valid in this range of \textit{moderate} wavelengths compared
to the water depth, notably since they miss their dispersive character.
In such regimes the characteristic speed is indeed a function of the
wavelength. The dispersion relation of water waves was derived from
potential flow theory by Airy and stands as a reference. The Serre--Green--Naghdi
(SGN) equations are a good candidate to capture dispersive effects.
These equations were first obtained by Serre \citep{serre_contribution_1953}
and extended in 2D by Green \& Naghdi \citep{green1976}. They can
be derived from the free-surface Euler equations as an asymptotic
model in the shallow-water regime (see \citep{lannes2013}), and,
in contrast to the original Boussinesq model \citep{boussinesq_theorie_1872},
without any assumption on the wave amplitude.

Over the past two decades, considerable work has been devoted to efficiently
discretizing the SGN equations, based on classical methods such as
finite volumes (FV), finite differences (FD), discontinuous Galerkin
(dG), finite elements (FE), or combinations of them, often referred
to as \textit{hybrid} methods. Some examples in 1D include \citep{dong2016}
(dG), \citep{docarmo2018,lemetayer2010} (hybrid FV/FD), \citep{mitsotakis2014}
(FE), \citep{li2014} (hybrid dG/FE). In 2D on Cartesian grids, one
may refer to the hybrid splitting methods \citep{bonneton2011,popinet2015},
or the dG approaches \citep{li2019,samii2018} for instance. Other
works also exploit the projection structure of the SGN equations \citep{fernandez2018,noelle2022}.
The list is of course not exhaustive. In 2D on unstructured meshes,
the list of existing works is much more limited. To the authors' knowledge,
the dG schemes \citep{duran_discontinuous_2017,marche_combined_2020}
and the hybrid FV/FE approach \citep{kazolea2023} remain among the
few available. Generally speaking, since the SGN equations are not
hyperbolic, these methods rely on an elliptic operator whose numerical
inversion is a well-known source of technical difficulty. Although
work has been carried out to mitigate these effects \citep{gavrilyuk2024,lannes_new_2015},
the inversion remains delicate, tends to dominate the computational
cost and limits scalability on parallel architectures. This is especially
the case in two dimensions: in widely used operational solvers such
as FUNWAVE-TVD \citep{shi2012} or the unstructured SGN code Uhaina
\citep{filippini2018}, the elliptic solve typically accounts for
a large share of the overall computing time.

Some hyperbolic Boussinesq-type models have been derived in recent
years, making it possible to partially circumvent these difficulties.
These strategies can be seen as approximations of the SGN equations
by hyperbolic relaxation models. One of the first works in that direction
was proposed by Favrie \& Gavrilyuk \citep{favrie_rapid_2017}. The
idea is to return to the Euler--Lagrange formulation of the system
on the basis of an \textit{augmented} Lagrangian, involving a parameter
ensuring the formal convergence to the original SGN equations when
it tends to infinity. This work has been followed by successive improvements
and extensions (see \citep{Tkachenko2023} for a 2D formulation, and
\citep{busto2021,guermond2022} for an extension to mild slopes and
full topography terms). At the same time, similar strategies have
been proposed, with the possibility of covering several families of
dispersive models. We refer for instance to \citep{escalante_efficient_2019,escalante2020},
where the relaxation method is introduced through a notion of artificial
compressibility. This pseudo-compressible strategy is also at the
core of operational non-hydrostatic ocean models such as SYMPHONIE-NH
\citep{auclair2011} and CROCO \citep{auclair2018}, albeit in a different,
three-dimensional modelling context. The model proposed in \citep{escalante_efficient_2019}
was derived shortly after in \citep{richard_extension_2021} from
the free-surface Euler equations in the shallow-water regime and under
a weak compressibility assumption. This formulation accounts for all
the bathymetry terms and improves the dispersive properties, following
a method similar to the one initially proposed in \citep{bonneton2011},
in the spirit of the enhanced Boussinesq-type models \citep{madsen1992,nwogu1993}.
Generally speaking, this new class of models appears to be better
suited to the construction of efficient schemes for dispersive flows:
no elliptic step, existence of methods adapted to hyperbolic problems,
better suitability for parallel computing. The counterpart is the
presence of a stiff source term that may significantly restrict the
time step. The splitting strategy adopted here alleviates this difficulty
by confining the stiff term to a dedicated, sub-iterated acoustic
step, leaving the gravity step to advance at its own time step.

Accordingly, major efforts have recently been made to develop efficient
numerical methods for this class of models. For instance, the capabilities
of \textit{discontinuous Galerkin} ADER-type methods have been widely
investigated in a series of recent papers \citep{bassi_hyperbolic_2020,busto2021,escalante_efficient_2019},
in 1D and in 2D on Cartesian grids. We can also mention \citep{escalante2020},
using the 1D finite-volume path-conservative schemes developed in
\citep{castro2006}, or the FE methods used in \citep{guermond2019}
and \citep{guermond2022}. In these last two works, numerical stability
is ensured through the positivity of the water depth and the discrete
preservation of steady states, thanks to a meticulous calibration
of the viscosity terms. However, in the context of hyperbolic approximations
of the SGN equations, access to an exact energy equation also makes
it possible to assess the stability of numerical schemes through discrete
energy estimates. While the methods mentioned above prove very effective
in practice, their stability in the specific sense of a fully discrete
energy (entropy) control, in both time and space and in 2D on general
meshes, has so far received comparatively little attention. In this
respect, the numerical approach proposed in \citep{favrie_rapid_2017},
based on a splitting strategy that separates a hyperbolic part from
an ODE part solved exactly, provides a first way to achieve such a
control, albeit in 1D and on a flat bottom. The main objective of
the present work is to show that a suitable splitting strategy, applied
to the model of \citep{richard_extension_2021}, provides a fully
discrete energy estimate (i.e. energy dissipation at the discrete
level, in both time and space) on general unstructured meshes, and
that low numerical diffusion can be achieved within this framework.
It enables the design of finite-volume-type schemes that are efficient
in low-Froude regimes, with a good compromise between accuracy and
algorithmic complexity. More precisely, the contributions are threefold:
(i) a splitting of the model into a gravity part and an acoustic part,
each inheriting an exact energy budget, which constitutes the central
object of this paper; (ii) for each subsystem, a finite-volume scheme
with local energy estimates established at the fully discrete level
on general unstructured meshes, together with the associated CFL conditions
and admissible ranges for the stabilization constants (the gravity
step adapts an existing shallow-water scheme, while the acoustic scheme
is introduced here); for the acoustic step, these sufficient nonlinear
conditions coincide with the necessary and sufficient conditions for
linear energy contraction over all Fourier modes on uniform square
meshes at constant depth; (iii) complementary linear analyses in the
energy norm, which select the minimal transport weight and long-wave
stabilization, within the chosen family, of the gravity step used
in the computations, beyond the range covered by the nonlinear estimate,
and quantify the numerical dissipation, supported by validations ranging
from academic benchmarks to a realistic configuration.

The paper is organized as follows: in Section~\ref{sec:splitting},
after introducing the 2D model under consideration and several notations,
we detail the proposed splitting strategy. The gravity part of the
splitting is addressed in Section~\ref{subsec:splitting_sw}, while
the so-called \textit{acoustic} part, which accounts for the dispersive
effects, is detailed in Section~\ref{subsec:splitting_acoustic}.
Section~\ref{sec:implementation} specifies some elements concerning
the practical numerical implementation, including integration in time,
the treatment of wet/dry fronts, boundary conditions and indicators
of computational performance. We propose in Section~\ref{sec:validation}
a series of numerical tests covering two academic benchmarks, a comparison
with experimental measurements, and an application to a realistic
configuration based on real data, as a first step towards operational
use. The proofs of the energy estimates are deferred to the appendices:
\ref{app:gravity-global-balance} for the global balance of the gravity
step and \ref{sec:dissipation} for the acoustic step. \ref{sec:linear_waves}
details the dispersion relation and \ref{sec:cnoidal} the computation
of cnoidal wave solutions of the relaxation model.

\section{Splitting strategy}\label{sec:splitting}

\subsection{Introduction}\label{subsec:splitting_intro}

An asymptotic approximation of the Serre--Green--Naghdi model was
first introduced in \citep{escalante_efficient_2019} and justified
in \citep{richard_extension_2021} by depth-averaging the free-surface
Euler equations under a weak compressibility assumption. In this work,
we restrict ourselves to a mild-slope formulation without the additional
bathymetry and enhanced dispersion terms of the full model, a restriction
under which the system retains a well-defined mathematical entropy
through an exact energy equation. The rigorous extension to fully
variable topography and enhanced dispersive properties is not trivial,
especially in two dimensions and if one aims to preserve an exact
energy balance. This extension will be the subject of future work.
Note, however, that a first study including improved dispersive properties
and an exact energy equation has recently been proposed in one dimension
\citep{duran2024}. In the present work, the core equations are those
derived in \citep{richard_extension_2021} and correspond to the following
system of partial differential equations, hereafter referred to as
the Leucothea (LCT) model:

\begin{subequations}
\label{eq:model} 
\begin{equation}
{\displaystyle \frac{\partial h}{\partial t}+\mathrm{div}\left(h\boldsymbol{u}\right)=0}\,,
\end{equation}
\begin{equation}
{\displaystyle \frac{\partial\left(h\boldsymbol{u}\right)}{\partial t}+\mathrm{div}\left(h\boldsymbol{u}\otimes\boldsymbol{u}\right)=-gh\boldsymbol{\nabla}\left(b+h\right)-\boldsymbol{\nabla}\left(hp\right)}\,,
\end{equation}
\begin{equation}
{\displaystyle \frac{\partial\left(hw\right)}{\partial t}+\mathrm{div}\left(hw\boldsymbol{u}\right)=\frac{3}{2}p}\,,
\end{equation}
\begin{equation}
{\displaystyle \frac{\partial\left(hp\right)}{\partial t}+\mathrm{div}\left(hp\boldsymbol{u}\right)=-a^{2}\left(2w+h\mathrm{div}\left(\boldsymbol{u}\right)\right)}\,,
\end{equation}
\end{subequations}
 which admits an exact energy conservation equation, also defining
a mathematical entropy 
\begin{equation}
{\displaystyle \frac{\partial\left(hE\right)}{\partial t}+\mathrm{div}\left(hE\boldsymbol{u}+\left(\frac{gh^{2}}{2}+hp\right)\boldsymbol{u}\right)=0}\,,\label{eq:energy_budget}
\end{equation}
where the total energy $E$ is 
\begin{equation}
{\displaystyle E=\frac{1}{2}\left\Vert \boldsymbol{u}\right\Vert ^{2}+\frac{1}{2}gh+gb+\frac{2}{3}w^{2}+\frac{1}{2a^{2}}p^{2}}\,.\label{eq:energyE}
\end{equation}
In the above system, in addition to the classical flow variables,
namely the water depth $h=h(t,\boldsymbol{x})$ and the depth-averaged
horizontal velocity $\boldsymbol{u}=\boldsymbol{u}(t,\boldsymbol{x})$,
functions of time $t\geq0$ and space $\boldsymbol{x}\in\mathbb{R}^{2}$,
one has to consider two additional variables: $w$, which is the depth-averaged
vertical velocity of the fluid, and $p$, the depth-averaged non-hydrostatic
pressure. The scalar function $b=b(\boldsymbol{x})$ encodes the topography
variations, and is assumed to be smooth and time-independent. The
constant $a>0$ denotes the artificial sound speed associated with
the weak compressibility assumption; the accuracy of the hyperbolic
approximation to the SGN equations improves as $a$ increases (see
Section~\ref{sec:implementation} for its practical choice). For
$h>0$, the system is hyperbolic with a complete set of eigenvectors
provided that $a^{2}+gh+p>0$. In order to numerically solve system
\eqref{eq:model}, we employ a splitting method that separates what
we call the \textit{shallow-water} (or \textit{gravity}) part: 
\begin{subequations}
\label{eq:splitting_sw} 
\begin{equation}
{\displaystyle \frac{\partial h}{\partial t}+\mathrm{div}\left(h\boldsymbol{u}\right)=0}\,,
\end{equation}
\begin{equation}
{\displaystyle \frac{\partial\left(h\boldsymbol{u}\right)}{\partial t}+\mathrm{div}\left(h\boldsymbol{u}\otimes\boldsymbol{u}\right)=-gh\boldsymbol{\nabla}\left(b+h\right)}\,,
\end{equation}
\begin{equation}
{\displaystyle \frac{\partial\left(hw\right)}{\partial t}+\mathrm{div}\left(hw\boldsymbol{u}\right)=0}\,,
\end{equation}
\begin{equation}
{\displaystyle \frac{\partial\left(hp\right)}{\partial t}+\mathrm{div}\left(hp\boldsymbol{u}\right)=0}\,,
\end{equation}
\end{subequations}
 and the so-called \textit{acoustic} part, governing the dispersive
effects: 
\begin{subequations}
\label{eq:splitting_acoustic} 
\begin{equation}
{\displaystyle \frac{\partial h}{\partial t}=0}\,,
\end{equation}
\begin{equation}
{\displaystyle \frac{\partial\left(h\boldsymbol{u}\right)}{\partial t}+\boldsymbol{\nabla}\left(hp\right)=0}\,,
\end{equation}
\begin{equation}
{\displaystyle \frac{\partial\left(hw\right)}{\partial t}=\frac{3}{2}p}\,,
\end{equation}
\begin{equation}
{\displaystyle \frac{\partial\left(hp\right)}{\partial t}=-a^{2}\left(2w+h\mathrm{div}\left(\boldsymbol{u}\right)\right)}\,.
\end{equation}
\end{subequations}
 At this stage it is important to note that each of these two systems
admits its own energy equation, namely for \eqref{eq:splitting_sw}
\begin{equation}
{\displaystyle \frac{\partial\left(hE\right)}{\partial t}+\mathrm{div}\left(\left(E+\frac{gh}{2}\right)h\boldsymbol{u}\right)=0}\,,\label{eq:energy_SW}
\end{equation}
while system \eqref{eq:splitting_acoustic} satisfies
\begin{equation}
{\displaystyle \frac{\partial\left(hE\right)}{\partial t}+\mathrm{div}\left(hp\boldsymbol{u}\right)=0}\,.\label{eq:energy_acoustic}
\end{equation}

We wish to emphasize this structure, which underpins the whole strategy.
At the continuous level, each subsystem possesses its own exact energy
balance, \eqref{eq:energy_SW} for the gravity part and \eqref{eq:energy_acoustic}
for the acoustic part, and the two energy fluxes recombine exactly
into the flux of the complete budget \eqref{eq:energy_budget}. The
splitting is therefore exact at the level of the energy: it distributes
the total energy budget between the two steps without creating or
destroying energy at the continuous level. This clean separation is
what makes the approach tractable, since each subsystem can then be
discretized independently by a scheme that dissipates its own discrete
energy, the control of the total energy being built step by step.
The splitting does introduce a time-splitting error; the numerical
results below suggest that it is not the dominant source of error
in the regimes and at the resolutions considered here.

\subsection{Notations}\label{subsec:splitting_notation}

\begin{figure}[!tbh]
\centering{}\includegraphics[width=0.4\textwidth]{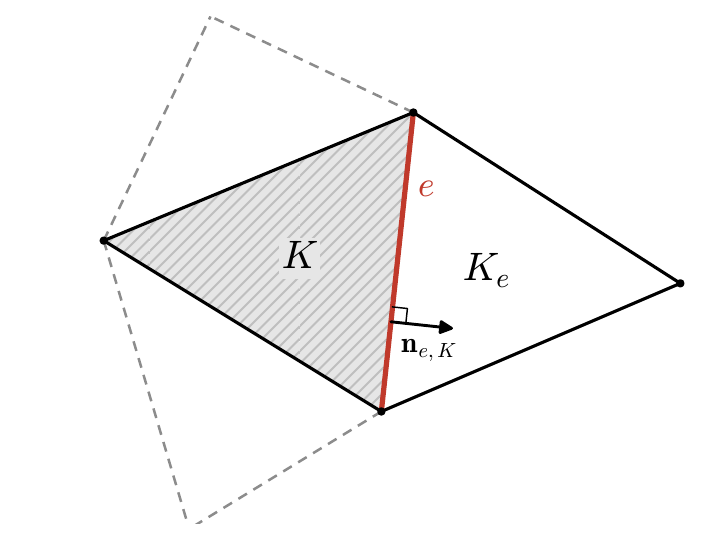}\caption{For any edge $e$ of a cell $K\in\mathbb{T}$, with area $m_{K}$
and perimeter $m_{\partial K}$, $K_{e}$ is the neighbouring cell
across $e$, $m_{e}$ its length and $\boldsymbol{n}_{e,K}$ the outward
normal. }\label{fig:geom}
\end{figure}

In order to introduce a finite-volume-type discretization of system
\eqref{eq:model}, a tessellation $\mathbb{T}$ of the computational
domain $\Omega\subset\mathbb{R}^{2}$ is considered. For a cell $K\in\mathbb{T}$,
we denote by $\partial K$ its boundary, $m_{K}$ its area and $m_{\partial K}$
its perimeter. For any edge $e\in\partial K$, we denote by $m_{e}$
its length and $\boldsymbol{n}_{e,K}$ the outward unit normal vector
pointing from cell $K$ to the neighbouring cell $K_{e}$ (see Fig.
\ref{fig:geom}). All the edges $e$ are collected in the set $\mathbb{E}$.

Let us now introduce some useful notations. For a piecewise constant
field $\omega$ over the cells $K$ of the tessellation $\mathbb{T}$,
scalar or vector valued, and for a given edge $e\in\mathbb{E}$ separating
two cells $K$ and $K_{e}$, we first define the two elementary interface
quantities,

\begin{equation}
{\displaystyle \overline{\omega}_{e}=\frac{\omega_{K}+\omega_{K_{e}}}{2}}\quad,\quad{\displaystyle \llbracket\omega\rrbracket_{e}=\frac{\omega_{K_{e}}-\omega_{K}}{2}}\,,\label{eq:notations_interface}
\end{equation}

the interface average and the interface half-jump, the latter oriented
from $K$ towards $K_{e}$. They are the only interface quantities
the transport fluxes require.

The stabilization terms, on the other hand, are discrete differential
operators, and they are built from the half-jump by inserting the
outward normal $\boldsymbol{n}_{e,K}$:

\begin{equation}
{\displaystyle \boldsymbol{\delta}_{e}\omega=\llbracket\omega\rrbracket_{e}\boldsymbol{n}_{e,K}}\quad,\quad{\displaystyle \delta_{e}\boldsymbol{\omega}=\llbracket\boldsymbol{\omega}\rrbracket_{e}\cdot\boldsymbol{n}_{e,K}}\,.\label{eq:notations_operators}
\end{equation}

The bold operator sends a scalar field to a vector and is the interface
analogue of the gradient; the plain one sends a vector field to a
scalar and is the analogue of the divergence.

Finally, the geometric factor

\[
r_{K}:=\frac{m_{\partial K}}{m_{K}}
\]

appears in every stabilization coefficient; on a uniform square grid
of step $\Delta$ it equals $4/\Delta$.

\subsection{Shallow-water subsystem}\label{subsec:splitting_sw}

\subsubsection{Numerical scheme}

The subsystem \eqref{eq:splitting_sw} corresponds to the two-dimensional
shallow-water system supplemented with two passive transport equations,
and can therefore be solved using any solver for the shallow-water
equations with topography. Bearing in mind the objectives of this
work, this method should be entropy-stable and readily extendable
to two-dimensional unstructured meshes. Among the possibilities available
in the literature, building on the general entropy-stability framework
introduced in \citep{tadmor1987,tadmor2003}, one may first consider
entropy-stable spatial discretizations. Fully discrete entropy-stable
constructions for general nonlinear hyperbolic systems on unstructured
meshes include \citep{guermond2016,kuzmin2022,martaud2023,gallice2022,gaburro2023}
and can be employed for the two-dimensional shallow-water equations.
Moreover, since the two energy balances associated with \eqref{eq:splitting_sw}
and \eqref{eq:splitting_acoustic} can be treated independently, schemes
based on different geometries can also be considered, as for example
the staggered approaches proposed in \citep{duran2020,herbin2023},
or the node-based finite-volume algorithm recently proposed in \citep{delgrosso2024}.
In the present work we chose a fully explicit version of the numerical
method arising from \citep{couderc2017}, and recently revisited and
extended in \citep{duran2026} (see also the initial works \citep{grenier2013}
and \citep{parisot2016} for a first semi-implicit version of this
scheme applicable to the multi-layer system).

Based on \citep{couderc2017} and the parametrized CPR formulation
recently analysed in \citep{duran2026}, the centred and fully upwind
CPR (Centered-Potential-Regularization) variants considered here for
\eqref{eq:splitting_sw} are the following: 
\begin{subequations}
\label{eq:scheme_sw} 
\begin{equation}
{\displaystyle h^{n+1}_{K}=h^{n}_{K}-\frac{\Delta t}{m_{K}}\sum_{e\in\partial K}\mathcal{F}^{\ast}_{e}\,m_{e}}\,,
\end{equation}

\begin{equation}
{\displaystyle h^{n+1}_{K}\,\boldsymbol{u}^{n+1}_{K}=h^{n}_{K}\,\boldsymbol{u}^{n}_{K}{\displaystyle -\frac{\Delta t}{m_{K}}\sum_{e\in\partial K}\left(\overline{\boldsymbol{u}}^{n}_{e}\,\mathcal{F}^{\ast}_{e}-\theta_{e}\llbracket\boldsymbol{u}^{n}\rrbracket_{e}\left|\mathcal{F}^{\ast}_{e}\right|\right)m_{e}{\displaystyle -{\displaystyle \frac{\Delta t}{m_{K}}h^{n}_{K}\sum_{e\in\partial K}\phi^{\ast}_{e}\,\boldsymbol{n}_{e,K}\,m_{e}}}}}\,,
\end{equation}

\begin{equation}
{\displaystyle h^{n+1}_{K}\,w^{n+1}_{K}=h^{n}_{K}\,w^{n}_{K}{\displaystyle -\frac{\Delta t}{m_{K}}\sum_{e\in\partial K}\left(\overline{w}^{n}_{e}\,\mathcal{F}^{\ast}_{e}-\theta_{e}\llbracket w^{n}\rrbracket_{e}\left|\mathcal{F}^{\ast}_{e}\right|\right)m_{e}}}\,,
\end{equation}

\begin{equation}
{\displaystyle h^{n+1}_{K}\,p^{n+1}_{K}=h^{n}_{K}\,p^{n}_{K}{\displaystyle -\frac{\Delta t}{m_{K}}\sum_{e\in\partial K}\left(\overline{p}^{n}_{e}\,\mathcal{F}^{\ast}_{e}-\theta_{e}\llbracket p^{n}\rrbracket_{e}\left|\mathcal{F}^{\ast}_{e}\right|\right)m_{e}}}\,,
\end{equation}

where the interface mass flux is

\begin{equation}
{\displaystyle \mathcal{F}^{*}_{e}=\left(\overline{h\boldsymbol{u}}^{n}_{e}-\boldsymbol{\Pi}_{e}\right)\cdot\boldsymbol{n}_{e,K}}\,,\label{eq:flux}
\end{equation}

and the weight $\theta_{e}\in[0,1]$ sets the amount of upwinding
in the transport by the interface mass flux $\mathcal{F}^{\ast}_{e}$.
It is common to the three transported quantities, the velocity $\boldsymbol{u}$
and the passive scalars $w$ and $p$, which are carried by the same
mass flux; its value is dictated by transport alone, whether of the
passive scalars or of the velocity within the gravity block, as shown
in Section~\ref{subsec:linear_sw}. The choice $\theta_{e}=0$ yields
the centred variant and $\theta_{e}=1$ the upwind variant. Indeed,
in the latter case, the second term reconstitutes the usual splitting
into positive and negative parts, since 
\[
\overline{\omega}_{e}\,\mathcal{F}^{\ast}_{e}-\llbracket\omega\rrbracket_{e}\left|\mathcal{F}^{\ast}_{e}\right|=\omega_{K}\left(\mathcal{F}^{\ast}_{e}\right)^{+}+\omega_{K_{e}}\left(\mathcal{F}^{\ast}_{e}\right)^{-}\,,
\]
with $\left(\mathcal{F}^{\ast}_{e}\right)^{\pm}=\left(\mathcal{F}^{\ast}_{e}\pm\vert\mathcal{F}^{\ast}_{e}\vert\right)/2$,
amounting to an upwinding of $\omega\in\{\boldsymbol{u},w,p\}$ with
respect to the sign of $\mathcal{F}^{\ast}_{e}$.

The two variants share the same energy stabilization terms and coincide
at rest; their respective stability properties are compared in the
linear analysis of Section~\ref{subsec:linear_sw}. The discrete
pressure relies on the following formula: 
\begin{equation}
{\displaystyle \phi^{\ast}_{e}=\overline{\phi}^{n}_{e}-\Lambda_{e}}\,,\label{eq:sw_scheme_phi_star}
\end{equation}
with $\phi^{n}_{K}=g(h^{n}_{K}+b_{K})$. The terms $\boldsymbol{\Pi}_{e}$
and $\Lambda_{e}$ govern the numerical viscosity of the scheme, and
admit the following form: 
\begin{equation}
{\displaystyle \boldsymbol{\Pi}_{e}=\gamma\,\Delta t\,\overline{rh^{n}\nu^{n}}_{e}\,\boldsymbol{\delta}_{e}\phi^{n}}\,,\label{eq:Pie}
\end{equation}
with the edge average $\overline{rh^{n}\nu^{n}}_{e}=\bigl(r_{K}h^{n}_{K}\nu^{n}_{K,e}+r_{K_{e}}h^{n}_{K_{e}}\nu^{n}_{K_{e},e}\bigr)/2$,
and the dimensionless viscosity amplification 
\begin{equation}
\nu^{n}_{K,e}=\left(1+\frac{\lvert\boldsymbol{u}^{n}_{K}\cdot\boldsymbol{n}_{e,K}\rvert}{c^{n}_{K}}\right)^{2},\qquad c^{n}_{K}=\sqrt{g\,h^{n}_{K}}\,,\label{eq:nu_edge}
\end{equation}
and 
\begin{equation}
{\displaystyle \Lambda_{e}=\alpha\,g\,\Delta t\,\overline{r}_{e}\,\delta_{e}(h^{n}\boldsymbol{u}^{n})}\,.\label{eq:Lambdae}
\end{equation}
\end{subequations}
 In \eqref{eq:nu_edge}, the amplification is attached to each side
of the edge. Since $(c^{n}_{K})^{2}=g\,h^{n}_{K}$, one has $g\,h^{n}_{K}\nu^{n}_{K,e}=(s^{n}_{K,e})^{2}$,
where $s^{n}_{K,e}=\lvert\boldsymbol{u}^{n}_{K}\cdot\boldsymbol{n}_{e,K}\rvert+c^{n}_{K}$
is the largest wave speed of cell $K$ normal to the edge $e$, so
that $\boldsymbol{\Pi}_{e}=\gamma\,\Delta t\,\overline{r\,(s^{n})^{2}}_{e}\,\boldsymbol{\delta}_{e}(h^{n}+b)$.
The mass-flux dissipation therefore acts on the free surface $h+b$,
not on the water depth $h$, and vanishes for the lake at rest. On
a uniform grid, this correction introduces a numerical diffusion coefficient
of order $\Delta t\,s^{2}$, rather than $s\,\Delta x$ as in first-order
Godunov-type schemes. Its time-step scaling is thus similar to that
of the Lax--Wendroff correction. At rest, $\nu_{K,e}=1$ and \eqref{eq:Pie}
reduces to the original CPR stabilization. The linear analysis below
examines the magnitude and directional dependence of this dissipation
on a uniform Cartesian grid. The two terms play complementary roles:
$\boldsymbol{\Pi}_{e}$ adds diffusion to the mass flux through the
jump of the free-surface potential $\phi$ (which gives the scheme
its name, Centered-Potential-Regularization), while $\Lambda_{e}$
regularizes the interface pressure through the jump of the discharge
$h\boldsymbol{u}$. Both vanish at the lake-at-rest equilibrium, ensuring
the well-balanced property, and are the terms responsible for the
discrete energy dissipation established in the following sections.
The momentum update is not strictly conservative: even over a flat
bottom with periodic boundaries, the stabilization through $\Lambda_{e}$
may produce a variation of the total discrete momentum. The constants
$\gamma>0,\alpha>0$ essentially govern the stability of the scheme,
and must be carefully calibrated, in relation to the CFL condition.
We first derive, from a linear stability analysis in the energy norm,
conditions for the dissipation of the linearized energy, necessary
and sufficient at rest, for the transport weight and in the long-wave
limit, before recalling in Section~\ref{subsec:nonlinear-sw} the
sufficient conditions ensuring local fully discrete energy estimates
in the nonlinear framework for the upwind version $\theta_{e}=1$,
based on \citep{couderc2017} and the fully upwind specialization
of the analysis developed in \citep[Remark 8]{duran2026}. The latter
fixes the structure of the operators, while the calibrated values
are those of the linear analysis, supplemented by a direct verification
of the global energy balance.

\subsubsection{Linear stability in the energy norm: gravity step}\label{subsec:linear_sw}

The stabilization operators of \eqref{eq:scheme_sw} have the structure
required by the fully discrete energy estimate recalled in Section~\ref{subsec:nonlinear-sw},
but that estimate only provides sufficient bounds on their intensity.
We calibrate here the three free quantities of the scheme, namely
the transport weight, the stabilization constant and the time step,
by a von Neumann analysis in the energy norm about a uniform moving
state. The analysis is organized in three steps, whose order matters:
the transport weight is fixed by transport alone (Step one); the stabilization
constant is then calibrated at fixed weight from the long-wave limit,
whose worst case is a flow and a wave along a mesh diagonal, which
yields a closed-form law (Step two); the resulting time-step restriction
is derived from the highest-frequency, checkerboard mode $(\xi_{x},\xi_{y})=(\pi,\pi)$,
and assessed numerically over frequencies and flow orientations (Step
three). The results are gathered in Proposition~\ref{prop:linear_sw}
at the end of this section. The analysis is carried out on a uniform
Cartesian grid, where Fourier techniques apply. It provides energy-dissipation
conditions for this problem, and certifies neither polygonal meshes
nor nonlinear states; for the latter, an exact global energy identity
is derived in Section~\ref{subsec:nonlinear-sw}.

\paragraph{Linearized setting}

The vector of primitive unknowns $\mathbf{U}=(h,u,v,w,p)^{T}$ is
linearized around a uniform state $(H>0,U,V,0,0)^{T}$ with a possibly
nonzero velocity, over a flat bottom $b=0$, on a grid with $\Delta x=\Delta y=\Delta$.
At this state all coefficients of the scheme are frozen. The transport
weight takes a constant value $\theta_{e}=\theta$. The two sides
of every edge carry the same base value, so that the edge average
in \eqref{eq:Pie} coincides with any other symmetric combination,
and the amplification \eqref{eq:nu_edge} takes, on the edges normal
to the coordinate axes, the values

\begin{equation}
\nu_{x}=\left(1+\frac{|U|}{c_{0}}\right)^{2}=\frac{\Lambda^{2}_{x}}{c^{2}_{0}},\qquad\nu_{y}=\left(1+\frac{|V|}{c_{0}}\right)^{2}=\frac{\Lambda^{2}_{y}}{c^{2}_{0}},\qquad\Lambda_{x}=|U|+c_{0},\quad\Lambda_{y}=|V|+c_{0},\quad c_{0}=\sqrt{gH},\label{eq:nu_lin}
\end{equation}

where $\Lambda_{x}$ and $\Lambda_{y}$ are the largest wave speeds
normal to the $x$- and $y$-edges. The analysis therefore applies
to the scheme as it is implemented, with no separate linear prescription.
We normalize the stabilization constants and restrict the analysis
to the one-parameter family specified by

\begin{equation}
\gamma_{\ell}:=4\gamma,\qquad\alpha_{\ell}:=2\alpha,\qquad\alpha_{\ell}=\gamma_{\ell}.\label{eq:linear_viscosity_choice}
\end{equation}

The first two equalities define the normalization; the third is a
choice of coefficients. Within this family, as shown below, the rest-state
threshold is $\gamma_{\ell}=1$.

\paragraph{Energy-norm criterion}

The energy associated with the linearized equations, given by the
quadratic part of the relative entropy with respect to the uniform
reference state, is

\begin{equation}
E_{\mathrm{lin}}=\frac{1}{2}gh^{2}+\frac{1}{2}H\left(u^{2}+v^{2}\right)+\frac{2}{3}Hw^{2}+\frac{1}{2a^{2}}Hp^{2}=\langle\mathbf{U},\mathbf{D}\mathbf{U}\rangle,\qquad\mathbf{D}=\mathrm{diag}\left(\frac{g}{2},\frac{H}{2},\frac{H}{2},\frac{2H}{3},\frac{H}{2a^{2}}\right).\label{eq:energy_lin}
\end{equation}

Looking for Fourier modes $\mathbf{U}^{n}_{\ell,m}=\hat{\mathbf{U}}^{n}e^{i(\ell\xi_{x}+m\xi_{y})}$,
where $(\ell,m)\in\mathbb{Z}^{2}$ indexes the cells of the uniform
grid and $(\xi_{x},\xi_{y})\in[-\pi,\pi]^{2}$ are the dimensionless
wavenumbers, $\xi_{x}=k_{x}\Delta$ and $\xi_{y}=k_{y}\Delta$, all
symbols below being $2\pi$-periodic in each variable, the scheme
reduces to $\hat{\mathbf{U}}^{n+1}=\mathbf{A}_{g}\hat{\mathbf{U}}^{n}$,
and the energy does not increase over one step if and only if

\begin{equation}
\left\langle \hat{\mathbf{U}}^{n},\left(\mathbf{D}-\mathbf{A}^{*}_{g}\mathbf{D}\mathbf{A}_{g}\right)\hat{\mathbf{U}}^{n}\right\rangle \geq0.\label{eq:linear_sw_dissip_en}
\end{equation}

Expressing $\mathbf{A}_{g}$ in an orthonormal basis for $\langle\cdot,\mathbf{D}\cdot\rangle$,
that is, introducing the symmetrized amplification matrix

\begin{equation}
\tilde{\mathbf{A}}_{g}:=\mathbf{D}^{1/2}\mathbf{A}_{g}\,\mathbf{D}^{-1/2},\label{eq:linear_sw_sym_matrix}
\end{equation}

one has the congruence

\begin{equation}
\mathbf{D}-\mathbf{A}^{*}_{g}\mathbf{D}\mathbf{A}_{g}=\mathbf{D}^{1/2}\,\mathbf{B}\,\mathbf{D}^{1/2},\qquad\mathbf{B}:=\mathbf{I}-\tilde{\mathbf{A}}^{*}_{g}\tilde{\mathbf{A}}_{g}.\label{eq:linear_sw_matrix_B}
\end{equation}

Since $\mathbf{D}^{1/2}$ is invertible, \eqref{eq:linear_sw_dissip_en}
holds for every $\hat{\mathbf{U}}^{n}$ if and only if $\mathbf{B}$
is positive semidefinite, $\mathbf{B}\geq0$, that is $\lVert\tilde{\mathbf{A}}_{g}\rVert_{2}\leq1$:
the energy decay rates are the eigenvalues of $\mathbf{B}$, one minus
the squared singular values of $\tilde{\mathbf{A}}_{g}$. The scheme
is thus energy-contractive if and only if $\mathbf{B}(\xi_{x},\xi_{y})\geq0$
for all $(\xi_{x},\xi_{y})$. Throughout the following analysis, $\mathbf{B}$
denotes $\mathbf{I}-\tilde{\mathbf{A}}^{*}_{g}\tilde{\mathbf{A}}_{g}$.

\paragraph{Block structure}

Around the state $\left(H,U,V,0,0\right)^{T}$ the auxiliary variables
$w$ and $p$ are advected by the frozen mass flux and decouple from
$(h,\boldsymbol{u})$, so the amplification matrix is block-diagonal,
\[
\mathbf{A}_{g}=\mathrm{diag}\left(\mathbf{A}^{h\boldsymbol{u}}_{g},\,G^{wp},\,G^{wp}\right),
\]
with a $3\times3$ block $\mathbf{A}^{h\boldsymbol{u}}_{g}$ acting
on $(h,\boldsymbol{u})$ and a single scalar factor $G^{wp}$ (identical
for $w$ and $p$, which obey the same transport) acting on each passive
scalar. As $\mathbf{D}=\mathrm{diag}(\mathbf{D}^{h\boldsymbol{u}},2H/3,H/(2a^{2}))$
is diagonal, with $\mathbf{D}^{h\boldsymbol{u}}=\mathrm{diag}(g/2,H/2,H/2)$,
the energy estimate \eqref{eq:linear_sw_dissip_en} decouples accordingly,
and holds if and only if both $\mathbf{A}^{h\boldsymbol{u}}_{g}$
and $G^{wp}$ are contractive in the energy norm. The analysis therefore
splits into a scalar transport problem for $G^{wp}$, to which the
transverse velocity of grid-aligned modes reduces as well (Step one),
and a coupled problem for $\mathbf{A}^{h\boldsymbol{u}}_{g}$ (rest
state, Steps two and three).

\paragraph{Energy dissipation, $L^{2}$ stability and the Kreiss matrix theorem}

The condition $\mathbf{B}(\xi)\geq0$ for all $\xi$ is necessary
and sufficient for the discrete energy not to increase over one time
step, for every initial datum. This is the property sought here: the
decay of the discrete energy is the discrete counterpart of the entropy
inequality that the splitting is designed to inherit, and the property
targeted by the nonlinear estimate of Section~\ref{subsec:nonlinear-sw}.
Below we prove the rest-state conditions over all frequencies, the
scalar transport criterion and the long-wave threshold uniform over
flow orientations. Step three gives an exact restriction for the checkerboard
mode $(\xi_{x},\xi_{y})=(\pi,\pi)$; its sufficiency over the full
Fourier domain is supported by a numerical sweep, rather than established
analytically.

Energy dissipation also implies $L^{2}$ stability: since $\mathbf{D}$
is constant, with $d_{\min}\mathbf{I}\leq\mathbf{D}\leq d_{\max}\mathbf{I}$,
the energy norm and the Euclidean norm are equivalent and $\lVert\mathbf{A}^{n}_{g}\rVert_{2}\leq(d_{\max}/d_{\min})^{1/2}$
for all $n$. The converse does not hold in general. By the Kreiss
matrix theorem, power boundedness of $\mathbf{A}_{g}$ is equivalent
to the existence of some possibly frequency-dependent Hermitian symmetrizer
$\mathbf{P}(\xi)$, satisfying $c\mathbf{I}\leq\mathbf{P}(\xi)\leq C\mathbf{I}$
for positive constants $c,C$ independent of $\xi$, with $\mathbf{A}^{*}_{g}\mathbf{P}\mathbf{A}_{g}\leq\mathbf{P}$
\citep[Chap. 4]{richtmyer1967}, and nothing forces an admissible
symmetrizer to be the matrix $\mathbf{D}$ dictated by the physical
energy: a scheme may let the energy increase transiently while remaining
$L^{2}$-stable. The two notions coincide, however, for scalar amplification
factors, since $|G(\xi)|^{n}$ is bounded if and only if $|G(\xi)|\leq1$.
This applies to the passive scalars, whose factor $G^{wp}$ is a diagonal
block of $\mathbf{A}_{g}$, and to the transverse velocity of grid-aligned
modes, which decouples from $h$ (Step one). For these scalar modes,
the transport condition on $\theta$ derived in Step one is therefore
necessary and sufficient for $L^{2}$ stability as well, and the centred
transport is $L^{2}$-unstable. Since these modes are part of the
full system, the transport condition is a necessary condition for
the stability of the whole gravity step, but not a sufficient one:
the coupled block acting on $(h,\boldsymbol{u})$ is studied separately.
At rest, the conditions of Lemma~\ref{lem:rest} ensure energy dissipation,
hence $L^{2}$ stability, over all frequencies. About a moving state,
the condition of Step two is necessary; together with the restriction
of Step three, it ensures energy dissipation, hence $L^{2}$ stability,
over the sampled Fourier modes, the sufficiency for all frequencies
resting on the numerical sweep.

\paragraph{The amplification matrix}

On a uniform Cartesian grid with $\Delta x=\Delta y=\Delta$, let
$z_{\ell,m}$ denote the value of a scalar grid function in cell $(\ell,m)$.
We define the centred first and three-point second differences by
\begin{equation}
\begin{aligned}(D^{1}_{x}z)_{\ell,m} & =\frac{z_{\ell+1,m}-z_{\ell-1,m}}{2\Delta}, & (D^{1}_{y}z)_{\ell,m} & =\frac{z_{\ell,m+1}-z_{\ell,m-1}}{2\Delta},\\
(D^{2}_{x}z)_{\ell,m} & =\frac{z_{\ell+1,m}-2z_{\ell,m}+z_{\ell-1,m}}{\Delta^{2}}, & (D^{2}_{y}z)_{\ell,m} & =\frac{z_{\ell,m+1}-2z_{\ell,m}+z_{\ell,m-1}}{\Delta^{2}}.
\end{aligned}
\label{eq:diffop}
\end{equation}
In the updates below, the cell indices are omitted and the operators
act componentwise. We also introduce 
\begin{equation}
\mathcal{T}=UD^{1}_{x}+VD^{1}_{y},\qquad\mathcal{D}=\frac{\Delta t\,\Delta}{2}\left(|U|D^{2}_{x}+|V|D^{2}_{y}\right).\label{eq:transport_op}
\end{equation}
Here $\mathcal{T}$ is the discrete advection operator and $\mathcal{D}$
is the dimensionless upwind diffusion in one time step.

With these notations and the directional speeds of \eqref{eq:nu_lin},
the linearized scheme reads

\begin{subequations}
\label{eq:scheme_sw_lin} 
\begin{equation}
h^{n+1}=h^{n}-\Delta t\,\mathcal{T}h^{n}-H\Delta t\left(D^{1}_{x}u^{n}+D^{1}_{y}v^{n}\right)+\frac{\gamma_{\ell}}{2}\Delta t^{2}\left(\Lambda^{2}_{x}D^{2}_{x}+\Lambda^{2}_{y}D^{2}_{y}\right)h^{n}\,,
\end{equation}

\begin{equation}
u^{n+1}=u^{n}-\Delta t\,\mathcal{T}u^{n}-g\Delta t\,D^{1}_{x}h^{n}+\gamma_{\ell}g\Delta t^{2}D^{2}_{x}\left(Uh^{n}+Hu^{n}\right)+\theta\,\mathcal{D}u^{n}\,,
\end{equation}

\begin{equation}
v^{n+1}=v^{n}-\Delta t\,\mathcal{T}v^{n}-g\Delta t\,D^{1}_{y}h^{n}+\gamma_{\ell}g\Delta t^{2}D^{2}_{y}\left(Vh^{n}+Hv^{n}\right)+\theta\,\mathcal{D}v^{n}\,,
\end{equation}

\begin{equation}
\omega^{n+1}=\omega^{n}-\Delta t\,\mathcal{T}\omega^{n}+\theta\,\mathcal{D}\omega^{n}\quad,\qquad\omega\in\{w,p\}\,.
\end{equation}
\end{subequations}
 The operator $\mathcal{T}$ is the full advection operator, not the
one-dimensional advection along each axis: the vector-transport flux
$\boldsymbol{G}_{e}=\overline{\boldsymbol{u}}_{e}\mathcal{F}^{\ast}_{e}$
carries both velocity components through every edge, so $u$ picks
up the transverse term $VD^{1}_{y}u$ and $v$ the term $UD^{1}_{x}v$.
Likewise the upwind operator $\mathcal{D}$ acts on both components
and on both passive scalars with the same two directional pieces;
setting $\theta=0$ recovers the centred variant, shown below to be
unstable about any moving state.

The cross terms $\gamma_{\ell}gU\Delta t^{2}D^{2}_{x}h$ and $\gamma_{\ell}gV\Delta t^{2}D^{2}_{y}h$
originate from $\Lambda_{e}$, defined in \eqref{eq:Lambdae}, diffusing
the discharge $h\boldsymbol{u}$ about a moving state, which couples
$h$ to the momentum equations; they vanish at rest.

Denoting $s_{\bullet}=i\sin\xi_{\bullet}$, $c_{\bullet}=1-\cos\xi_{\bullet}\geq0$,
the Courant number $\sigma_{0}=c_{0}\Delta t/\Delta$ based on the
celerity, the gravity-step Courant number $\sigma_{g}=(\|\boldsymbol{U}\|+c_{0})\Delta t/\Delta=(1+Fr)\,\sigma_{0}$
with $Fr=\|\boldsymbol{U}\|/c_{0}$, the viscous Courant numbers $\sigma_{x}=\Lambda_{x}\Delta t/\Delta$
and $\sigma_{y}=\Lambda_{y}\Delta t/\Delta$ (both equal to $\sigma_{0}$
at rest), and the advective ones $\tilde{\sigma}_{x}=U\Delta t/\Delta$
and $\tilde{\sigma}_{y}=V\Delta t/\Delta$ (both zero at rest), two
scalars recur on every diagonal entry:

\begin{equation}
\mathfrak{t}:=\tilde{\sigma}_{x}s_{x}+\tilde{\sigma}_{y}s_{y}\quad\text{(advection)}\,,\qquad\mathfrak{d}:=\theta\left(\left|\tilde{\sigma}_{x}\right|c_{x}+\left|\tilde{\sigma}_{y}\right|c_{y}\right)\quad\text{(upwinding)}\,.\label{eq:symbols_ad}
\end{equation}

The first is purely imaginary and carries no dissipation; the second
is real, non-negative, and is the only place where $\theta$ appears.
With them the amplification matrix reads

\[
\mathbf{A}^{h\boldsymbol{u}}_{g}=\left(\begin{array}{ccc}
1-\gamma_{\ell}\left(\sigma^{2}_{x}c_{x}+\sigma^{2}_{y}c_{y}\right)-\mathfrak{t} & -\sqrt{\dfrac{H}{g}}\,\sigma_{0}s_{x} & -\sqrt{\dfrac{H}{g}}\,\sigma_{0}s_{y}\\[6pt]
-\sqrt{\dfrac{g}{H}}\,\sigma_{0}\left(s_{x}+2\gamma_{\ell}\tilde{\sigma}_{x}c_{x}\right) & 1-2\gamma_{\ell}\sigma^{2}_{0}c_{x}-\mathfrak{t}-\mathfrak{d} & 0\\[6pt]
-\sqrt{\dfrac{g}{H}}\,\sigma_{0}\left(s_{y}+2\gamma_{\ell}\tilde{\sigma}_{y}c_{y}\right) & 0 & 1-2\gamma_{\ell}\sigma^{2}_{0}c_{y}-\mathfrak{t}-\mathfrak{d}
\end{array}\right)\,.
\]

\paragraph{The rest state}

In full generality, $\mathbf{B}^{h\boldsymbol{u}}(\xi_{x},\xi_{y}):=\mathbf{I}-\tilde{\mathbf{A}}^{h\boldsymbol{u},*}_{g}\tilde{\mathbf{A}}^{h\boldsymbol{u}}_{g}$
has no closed-form spectrum, so that neither the admissible values
of $\gamma_{\ell}$ nor the CFL-type condition on $\Delta t$ ensuring
$\mathbf{B}^{h\boldsymbol{u}}\geq0$ for all $(\xi_{x},\xi_{y})$
can be read off directly. The rest state is the natural starting point,
for three reasons: it fixes the normalization $\gamma_{\ell}=1$ of
the stabilization constant; it is the limit $Fr\to0$ of the laws
$\Gamma(Fr)$ and $\sigma_{\max}(Fr)$ derived in Steps two and three;
and it is the only configuration where the contractivity of the coupled
block is proved over the whole Fourier domain (Lemma~\ref{lem:rest}).
Moreover, at rest all variants coincide, the cross-coupling terms
vanish, and the advective and viscous diagonals symmetrize. Along
the diagonal $\xi_{x}=\xi_{y}$, the spectrum is explicit and yields
necessary conditions, which Lemma~\ref{lem:rest} below shows to
be sufficient over the whole Fourier domain. Setting $\xi_{x}=\xi_{y}=\xi$,
the three diagonal entries of $\mathbf{A}^{h\boldsymbol{u}}_{g}$
become equal and the symmetrization is computed by hand. Setting 
\[
\delta:=1-2\gamma_{\ell}\sigma^{2}_{0}\left(1-\cos\xi\right),\quad\zeta:=i\sigma_{0}\sin\xi\,,
\]
the square root $\mathbf{D}^{1/2}$ homogenizes the two coupling terms
and the symmetrized matrix \eqref{eq:linear_sw_sym_matrix} takes
the form 
\[
\tilde{\mathbf{A}}^{h\boldsymbol{u}}_{g}=(\mathbf{D}^{h\boldsymbol{u}})^{1/2}\,\mathbf{A}^{h\boldsymbol{u}}_{g}\,(\mathbf{D}^{h\boldsymbol{u}})^{-1/2}=\left(\begin{array}{ccc}
\delta & -\zeta & -\zeta\\
-\zeta & \delta & 0\\
-\zeta & 0 & \delta
\end{array}\right).
\]
The product $\tilde{\mathbf{A}}^{h\boldsymbol{u},*}_{g}\tilde{\mathbf{A}}^{h\boldsymbol{u}}_{g}$
is then computed directly. Setting $\tau:=\sigma^{2}_{0}\sin^{2}\xi\geq0$,
\[
\tilde{\mathbf{A}}^{h\boldsymbol{u},*}_{g}\,\tilde{\mathbf{A}}^{h\boldsymbol{u}}_{g}=\left(\begin{array}{ccc}
\delta^{2}+2\tau & 0 & 0\\
0 & \delta^{2}+\tau & \tau\\
0 & \tau & \delta^{2}+\tau
\end{array}\right).
\]
The block-diagonal structure of $\mathbf{B}^{h\boldsymbol{u}}$ makes
its spectrum immediate. Expanding $\delta$ and $\tau$, we find the
three eigenvalues of $\mathbf{B}^{h\boldsymbol{u}}$: 
\[
\begin{array}{lllll}
\lambda_{1} & = & 1-\delta^{2} & = & 4\gamma_{\ell}\sigma^{2}_{0}\left(1-\cos\xi\right)\left[\gamma_{\ell}\sigma^{2}_{0}\left(\cos\xi-1\right)+1\right]\,,\\
\\\lambda_{2,3} & = & 1-\delta^{2}-2\tau & = & 2\sigma^{2}_{0}\left(1-\cos\xi\right)\left[2\gamma^{2}_{\ell}\sigma^{2}_{0}\left(\cos\xi-1\right)+2\gamma_{\ell}-(1+\cos\xi)\right]\,,
\end{array}
\]
the eigenvalue $\lambda_{2,3}$ being of multiplicity two. Stability
in the energy sense is equivalent to $\lambda_{1}\geq0$ and $\lambda_{2,3}\geq0$
for all $\xi$. As $\lambda_{1}-\lambda_{2,3}=2\tau\geq0$, only the
condition $\lambda_{2,3}\geq0$ is constraining: 
\[
\lambda_{2,3}=2\sigma^{2}_{0}\left(1-\cos\xi\right)f\left(\cos\xi\right)\quad,\qquad f(x):=2\gamma_{\ell}-(1+x)-2\gamma^{2}_{\ell}\sigma^{2}_{0}(1-x)\,.
\]
The prefactor $2\sigma^{2}_{0}(1-\cos\xi)$ being strictly positive
for $\sigma_{0}>0$ and nonzero frequencies, the condition $\lambda_{2,3}\geq0$
for all $\xi$ is equivalent, by continuity, to $f(x)\geq0$ for all
$x=\cos\xi\in[-1,1]$. As $f$ is affine in $x$, it suffices to control
its two values at the endpoints, which decouples the two necessary
conditions on the single constant $\gamma_{\ell}$ of \eqref{eq:linear_viscosity_choice}:
\begin{equation}
\gamma_{\ell}\geq1\quad\text{and}\quad\gamma_{\ell}\,\sigma^{2}_{0}\leq\frac{1}{2}\,.\label{eq:rest_conditions}
\end{equation}
The first constraint, $\gamma_{\ell}\geq1$, stems from the large-wavelength
limit $\xi\rightarrow0$ and is necessary independently of the time
step; the second is the CFL-type condition proper, arising from the
limit $\xi\rightarrow\pi$.
\begin{lem}
\label{lem:rest} Let $\sigma_{0}>0$. At rest, $\mathbf{B}^{h\boldsymbol{u}}(\xi_{x},\xi_{y})\geq0$
for every $(\xi_{x},\xi_{y})\in[-\pi,\pi]^{2}$ if and only if the
conditions \eqref{eq:rest_conditions} hold.
\end{lem}
\begin{proof}
Necessity follows from the diagonal $\xi_{x}=\xi_{y}$ treated above.
For sufficiency, the symmetrized amplification matrix at rest reads
\[
\tilde{\mathbf{A}}^{h\boldsymbol{u}}_{g}=\mathbf{I}-\gamma_{\ell}\sigma^{2}_{0}\,\mathbf{C}-i\sigma_{0}\,\mathbf{S},\qquad\mathbf{C}=\mathrm{diag}\left(c_{x}+c_{y},\,2c_{x},\,2c_{y}\right),\qquad\mathbf{S}=\left(\begin{array}{ccc}
0 & \sin\xi_{x} & \sin\xi_{y}\\
\sin\xi_{x} & 0 & 0\\
\sin\xi_{y} & 0 & 0
\end{array}\right),
\]
where $\mathbf{C}$ and $\mathbf{S}$ are real symmetric. Consider
first $\gamma_{\ell}=1$, and set $\mathbf{A}_{1}(\sigma):=\mathbf{I}-\sigma^{2}\mathbf{C}-i\sigma\mathbf{S}$
for $\sigma\in\mathbb{R}$. Expanding the product gives 
\[
\mathbf{I}-\mathbf{A}_{1}(\sigma)^{*}\mathbf{A}_{1}(\sigma)=\sigma^{2}\,\mathbf{G}(\sigma),\qquad\mathbf{G}(\sigma):=2\mathbf{C}-\mathbf{S}^{2}+i\sigma\left(\mathbf{S}\mathbf{C}-\mathbf{C}\mathbf{S}\right)-\sigma^{2}\mathbf{C}^{2},
\]
where $\mathbf{G}(\sigma)$ is Hermitian, $\mathbf{S}\mathbf{C}-\mathbf{C}\mathbf{S}$
being real and skew-symmetric. At $\sigma=1/\sqrt{2}$, using $\sin^{2}\xi_{j}=2c_{j}-c^{2}_{j}$
and writing $d=\cos\xi_{x}-\cos\xi_{y}$, we obtain 
\[
\mathbf{G}(1/\sqrt{2})=\begin{pmatrix}d^{2}/2 & -id\sin\xi_{x}/\sqrt{2} & id\sin\xi_{y}/\sqrt{2}\\
id\sin\xi_{x}/\sqrt{2} & \sin^{2}\xi_{x} & -\sin\xi_{x}\sin\xi_{y}\\
-id\sin\xi_{y}/\sqrt{2} & -\sin\xi_{x}\sin\xi_{y} & \sin^{2}\xi_{y}
\end{pmatrix}=\boldsymbol{\psi}\boldsymbol{\psi}^{*},\qquad\boldsymbol{\psi}=\begin{pmatrix}d/\sqrt{2}\\
i\sin\xi_{x}\\
-i\sin\xi_{y}
\end{pmatrix}.
\]
Thus $z^{*}\mathbf{G}(1/\sqrt{2})z=|\boldsymbol{\psi}^{*}z|^{2}\geq0$
for every $z\in\mathbb{C}^{3}$, proving $\mathbf{G}(1/\sqrt{2})\geq0$.
Since $\mathbf{A}_{1}(-\sigma)$ is the complex conjugate of $\mathbf{A}_{1}(\sigma)$,
$\mathbf{G}(-1/\sqrt{2})=\overline{\mathbf{G}(1/\sqrt{2})}\geq0$
as well. To extend this positivity to $|\sigma|\leq1/\sqrt{2}$, write
\[
\mathbf{G}(\sigma)=\mathbf{G}_{0}+\sigma\mathbf{G}_{1}-\sigma^{2}\mathbf{C}^{2},\qquad\mathbf{G}_{0}=2\mathbf{C}-\mathbf{S}^{2},\qquad\mathbf{G}_{1}=i(\mathbf{S}\mathbf{C}-\mathbf{C}\mathbf{S}).
\]
The endpoint values are therefore 
\[
\mathbf{G}(\pm1/\sqrt{2})=\mathbf{G}_{0}\pm\frac{1}{\sqrt{2}}\mathbf{G}_{1}-\frac{1}{2}\mathbf{C}^{2}.
\]
Multiplying these two identities by $(1+\sqrt{2}\sigma)/2$ and $(1-\sqrt{2}\sigma)/2$,
respectively, and adding them gives 
\[
\frac{1+\sqrt{2}\sigma}{2}\mathbf{G}(1/\sqrt{2})+\frac{1-\sqrt{2}\sigma}{2}\mathbf{G}(-1/\sqrt{2})=\mathbf{G}_{0}+\sigma\mathbf{G}_{1}-\frac{1}{2}\mathbf{C}^{2}.
\]
Consequently,
\[
\mathbf{G}(\sigma)=\frac{1+\sqrt{2}\,\sigma}{2}\,\mathbf{G}(1/\sqrt{2})+\frac{1-\sqrt{2}\,\sigma}{2}\,\mathbf{G}(-1/\sqrt{2})+\left(\frac{1}{2}-\sigma^{2}\right)\mathbf{C}^{2}\geq0\,.
\]
Indeed, all three coefficients are nonnegative when $|\sigma|\leq1/\sqrt{2}$,
both endpoint matrices are positive semidefinite, and $z^{*}\mathbf{C}^{2}z=\|\mathbf{C}z\|^{2}\geq0$.
Hence $\lVert\mathbf{A}_{1}(\sigma)\rVert_{2}\leq1$ whenever $\sigma^{2}\leq1/2$.
For $\gamma_{\ell}\geq1$, set $t=\sqrt{\gamma_{\ell}}\,\sigma_{0}$.
Then 
\[
\tilde{\mathbf{A}}^{h\boldsymbol{u}}_{g}=\frac{1+\gamma^{-1/2}_{\ell}}{2}\,\mathbf{A}_{1}(t)+\frac{1-\gamma^{-1/2}_{\ell}}{2}\,\mathbf{A}_{1}(-t)\,,
\]
a convex combination of two contractions as soon as $t^{2}=\gamma_{\ell}\sigma^{2}_{0}\leq1/2$.
Therefore $\lVert\tilde{\mathbf{A}}^{h\boldsymbol{u}}_{g}\rVert_{2}\leq1$,
that is $\mathbf{B}^{h\boldsymbol{u}}\geq0$, for every $(\xi_{x},\xi_{y})$.
\end{proof}

We retain the smallest admissible value $\gamma_{\ell}=1$: at rest,
it maximizes the admissible time step, $\sigma_{0}\leq1/\sqrt{2}$,
and gives the smallest stabilization coefficients; then 
\[
\left.\lambda_{2,3}\right|_{\gamma_{\ell}=1}=2\sigma^{2}_{0}\left(1-2\sigma^{2}_{0}\right)\left(1-\cos\xi\right)^{2}\geq0\quad\Longleftrightarrow\quad\sigma^{2}_{0}\leq\frac{1}{2}\,.
\]
 This is the normalization we adopt: $\gamma_{\ell}=1$ is the rest-state
linear threshold.

The rest state is, however, not representative. Two distinct mechanisms
are lost when $\boldsymbol{U}=0$, and both must be recovered before
the constants can be calibrated for an application. We treat them
in turn, and the order matters: the transport condition fixes $\theta$
without reference to $\gamma_{\ell}$, after which $\gamma_{\ell}$
can be calibrated at fixed $\theta$.

\paragraph{Step one: the transport weight}

For $\xi_{y}=0$, $s_{y}=c_{y}=0$ and the transverse velocity $v$
decouples from $(h,u)$ in $\mathbf{A}^{h\boldsymbol{u}}_{g}$: it
is a shear mode, on which $\boldsymbol{\Pi}_{e}$ and $\Lambda_{e}$
have no effect, with amplification factor

\begin{equation}
G_{\perp}=1-\tilde{\sigma}_{x}s_{x}-\theta\left|\tilde{\sigma}_{x}\right|c_{x}\,.\label{eq:Gperp}
\end{equation}

For $\theta=0$,

\begin{equation}
\left|G_{\perp}\right|^{2}=1+\tilde{\sigma}^{2}_{x}\sin^{2}\xi_{x}>1\label{eq:centred_growth}
\end{equation}

as soon as $\tilde{\sigma}_{x}\sin\xi_{x}\neq0$, whatever $\gamma_{\ell}$.
Upwinding is therefore required by the gravity block itself: the centred
variant is unstable about any moving state.

The passive scalars are purely advected by the frozen mass flux, with
the factor

\begin{equation}
G^{wp}=1-\left(\tilde{\sigma}_{x}s_{x}+\tilde{\sigma}_{y}s_{y}\right)-\theta\left(\left|\tilde{\sigma}_{x}\right|c_{x}+\left|\tilde{\sigma}_{y}\right|c_{y}\right),\label{eq:awp}
\end{equation}

which reduces to \eqref{eq:Gperp} for $\xi_{y}=0$. Let $\lVert\tilde{\sigma}\rVert_{1}=|\tilde{\sigma}_{x}|+|\tilde{\sigma}_{y}|$.
\begin{lem}
\label{lem:awp} Assume $\lVert\tilde{\sigma}\rVert_{1}\leq1$ and
$0\leq\theta\leq1$. Then $|G^{wp}|\leq1$ for every Fourier mode
if and only if $\theta\geq\lVert\tilde{\sigma}\rVert_{1}$.
\end{lem}
\begin{proof}
Assume first $\tilde{\sigma}_{x},\tilde{\sigma}_{y}\geq0$. If $\lVert\tilde{\sigma}\rVert_{1}=0$,
then $G^{wp}=1$ and the claim is immediate. Otherwise, let 
\[
G^{LW}(\xi)=1-\theta^{2}(1-\cos\xi)-i\,\theta\sin\xi,\qquad|G^{LW}(\xi)|^{2}=1-\theta^{2}(1-\theta^{2})(1-\cos\xi)^{2},
\]
the amplification factor of the one-dimensional Lax--Wendroff scheme
at Courant number $\theta$, of modulus at most one for $\theta\leq1$.
\emph{Sufficiency.} If $\theta\geq\lVert\tilde{\sigma}\rVert_{1}$,
a direct computation from \eqref{eq:awp} gives 
\[
G^{wp}=\Bigl(1-\frac{\lVert\tilde{\sigma}\rVert_{1}}{\theta}\Bigr)+\frac{\tilde{\sigma}_{x}}{\theta}\,G^{LW}(\xi_{x})+\frac{\tilde{\sigma}_{y}}{\theta}\,G^{LW}(\xi_{y}),
\]
a convex combination of $1$, $G^{LW}(\xi_{x})$ and $G^{LW}(\xi_{y})$,
hence $|G^{wp}|\leq1$. \emph{Necessity.} For $\xi_{x}=\xi_{y}=\varepsilon$,
$G^{wp}=1-\theta\lVert\tilde{\sigma}\rVert_{1}(1-\cos\varepsilon)-i\,\lVert\tilde{\sigma}\rVert_{1}\sin\varepsilon$,
hence 
\[
1-|G^{wp}|^{2}=\lVert\tilde{\sigma}\rVert_{1}\bigl(\theta-\lVert\tilde{\sigma}\rVert_{1}\bigr)\varepsilon^{2}+O(\varepsilon^{4}),
\]
which is negative for all sufficiently small $\varepsilon\neq0$ if
$\theta<\lVert\tilde{\sigma}\rVert_{1}$. Finally, the case of arbitrary
signs reduces to this one through $\xi_{i}\mapsto-\xi_{i}$, which
flips the sign of $s_{i}$ and leaves $c_{i}$ unchanged.
\end{proof}

The transport step thus blends the identity with one-dimensional Lax--Wendroff
steps at Courant number $\theta$, one per direction; the threshold
$\theta=\lVert\tilde{\sigma}\rVert_{1}$ removes the identity part
and corresponds to the minimal stable viscosity within this family,
as Lax--Wendroff does in one dimension. This minimal weight ensures
$L^{2}$ stability, not monotonicity: for $\theta<1$ the transport
step does not satisfy a discrete maximum principle.

The condition $\theta\geq\lVert\tilde{\sigma}\rVert_{1}$ is thus
necessary and sufficient for the passive scalars and for the transverse
modes of the gravity block, and fixes the single weight $\theta_{e}$
of \eqref{eq:scheme_sw}; it says nothing about the coupled modes
of $(h,\boldsymbol{u})$, whose stability is treated at fixed $\theta$
in Step two. Since $\lVert\tilde{\sigma}\rVert_{1}=\sigma_{0}(|F_{x}|+|F_{y}|)\leq\sqrt{2}\,\sigma_{0}Fr$,
with $(F_{x},F_{y})=\boldsymbol{U}/c_{0}$ and equality for a diagonal
flow, the smallest isotropic law satisfying it for every orientation
is

\begin{equation}
\theta=\min\!\left(1,\ \sqrt{2}\,\sigma_{0}\,Fr\right)=\min\!\left(1,\ \sqrt{2}\,\frac{\lVert\boldsymbol{U}\rVert\Delta t}{\Delta}\right)\,.\label{eq:theta_law}
\end{equation}

The cap in \eqref{eq:theta_law} is inactive whenever the scheme is
energy-contractive, as shown in Step three; we therefore use the reduced
weight $q:=\theta/\sigma_{0}=\sqrt{2}\,Fr$, independent of the time
step, in Steps two and three.

 On a general mesh, the metric $(r_{K}+r_{K_{e}})/8$ recovers $1/\Delta$
on square cells, which gives the edge prescription

\begin{equation}
\theta_{e}=\min\!\left(1,\ \frac{1}{2\sqrt{2}}\,\lVert\overline{\boldsymbol{u}}_{e}\rVert\,\Delta t\,\overline{r}_{e}\right)\,,\label{eq:theta_mesh}
\end{equation}

where $\overline{\boldsymbol{u}}_{e}$ is the edge average of the
velocity (any symmetric combination, a maximum for instance, is equally
admissible).

\paragraph{Step two: the stabilization constant}

With $\theta$ fixed by \eqref{eq:theta_law}, it remains to calibrate
$\gamma_{\ell}$ on the coupled block $\mathbf{A}^{h\boldsymbol{u}}_{g}$.
Unlike at rest, its spectrum is not available in closed form, but
the long-wave limit $(\xi_{x},\xi_{y})\to(0,0)$ along a fixed propagation
direction yields a necessary condition in closed form. The argument
proceeds in three stages: (a) at leading order, the dissipation matrix
reduces to $\sigma^{2}_{0}$ times a real symmetric matrix $\mathbf{L}$,
which does not depend on the time step, so that its positivity can
only constrain $\gamma_{\ell}$; (b) among all flow and wave directions,
the worst case is a flow and a wave aligned with the same mesh diagonal,
where $\mathbf{L}$ reduces to a $2\times2$ block whose determinant
gives a quadratic inequality in $\gamma_{\ell}$; (c) this inequality
yields the threshold $\gamma_{\ell}\geq\Gamma(Fr)$, and Lemma~\ref{lem:diag}
shows that it suffices for every orientation.

In this limit, with $\lVert\xi\rVert=\sqrt{\xi^{2}_{x}+\xi^{2}_{y}}$,
the restriction $\mathbf{B}^{h\boldsymbol{u}}$ of $\mathbf{B}$ to
the coupled block satisfies

\begin{equation}
\mathbf{B}^{h\boldsymbol{u}}=\sigma^{2}_{0}\,\mathbf{L}(\xi_{x},\xi_{y})+O\bigl(\lVert\xi\rVert^{3}\bigr)\,.\label{eq:longwave}
\end{equation}

The first-order term in $\mathbf{B}^{h\boldsymbol{u}}$ vanishes because
the first-order coefficient in the expansion of $\tilde{\mathbf{A}}^{h\boldsymbol{u}}_{g}$
is skew-Hermitian, and all second-order terms scale as $\sigma^{2}_{0}$
since $\theta=q\,\sigma_{0}$. Energy contractivity therefore requires
$\mathbf{L}(\xi_{x},\xi_{y})$ to be positive semidefinite in every
propagation direction. We give this matrix below, then evaluate it
in the diagonal configuration.

Here $\mathbf{L}$ is a real symmetric matrix, homogeneous of degree
two in the Fourier variables, so that its sign depends only on the
propagation direction. With $F_{\xi}=F_{x}\xi_{x}+F_{y}\xi_{y}$ and
$d_{F}=|F_{x}|\xi^{2}_{x}+|F_{y}|\xi^{2}_{y}$, its entries read

\begin{equation}
\mathbf{L}(\xi_{x},\xi_{y})=\begin{pmatrix}\begin{gathered}\gamma_{\ell}\left[(1+|F_{x}|)^{2}\xi^{2}_{x}+(1+|F_{y}|)^{2}\xi^{2}_{y}\right]\\
-F^{2}_{\xi}-(\xi^{2}_{x}+\xi^{2}_{y})
\end{gathered}
 & \gamma_{\ell}F_{x}\xi^{2}_{x}-2\xi_{x}F_{\xi} & \gamma_{\ell}F_{y}\xi^{2}_{y}-2\xi_{y}F_{\xi}\\[3pt]
\gamma_{\ell}F_{x}\xi^{2}_{x}-2\xi_{x}F_{\xi} & (2\gamma_{\ell}-1)\xi^{2}_{x}+q\,d_{F}-F^{2}_{\xi} & -\xi_{x}\xi_{y}\\[3pt]
\gamma_{\ell}F_{y}\xi^{2}_{y}-2\xi_{y}F_{\xi} & -\xi_{x}\xi_{y} & (2\gamma_{\ell}-1)\xi^{2}_{y}+q\,d_{F}-F^{2}_{\xi}
\end{pmatrix}\,.\label{eq:Lmat}
\end{equation}

Since $\mathbf{L}$ does not depend on the time step, a negative eigenvalue
cannot be removed by reducing $\Delta t$. The long-wave limit thus
gives a necessary lower bound on $\gamma_{\ell}$, which depends only
on $Fr$ once all flow and wave directions are considered; the time
step is constrained at finite frequencies, in Step three.

Parametrizing the flow and wave directions by the angles $\varphi_{\text{flow}}$
and $\varphi_{\text{wave}}$,
\[
(F_{x},F_{y})=Fr(\cos\varphi_{\text{flow}},\sin\varphi_{\text{flow}}),\qquad(\xi_{x},\xi_{y})=\lVert\xi\rVert(\cos\varphi_{\text{wave}},\sin\varphi_{\text{wave}}),
\]
the worst configuration announced above reads

\begin{equation}
\varphi_{\text{flow}}=\varphi_{\text{wave}}=\frac{\pi}{4}\quad\Longleftrightarrow\quad F_{x}=F_{y}=\frac{Fr}{\sqrt{2}}\,,\quad\xi_{x}=\xi_{y}\,.\label{eq:worst_case}
\end{equation}

In this configuration the calculation closes. We evaluate $\mathbf{L}$
at the unit wave vector $\xi_{x}=\xi_{y}=1/\sqrt{2}$. In the variables
$(h,u_{\parallel},u_{\perp})$, where $u_{\parallel}=(u+v)/\sqrt{2}$
and $u_{\perp}=(u-v)/\sqrt{2}$ are the velocity components along
and across the flow, $\mathbf{L}$ is block diagonal:

\begin{equation}
\mathbf{L}=\begin{pmatrix}\mathbf{L}_{\parallel} & 0\\
0 & \gamma_{\ell}
\end{pmatrix},\qquad\mathbf{L}_{\parallel}=\begin{pmatrix}\gamma_{\ell}\bigl(1+\dfrac{Fr}{\sqrt{2}}\bigr)^{2}-(1+Fr^{2}) & \dfrac{1}{2}Fr(\gamma_{\ell}-4)\\[3pt]
\dfrac{1}{2}Fr(\gamma_{\ell}-4) & \gamma_{\ell}-1
\end{pmatrix}.\label{eq:Lpar}
\end{equation}

The transverse velocity is decoupled and contributes the eigenvalue
$\gamma_{\ell}>0$; $\mathbf{L}_{\parallel}$ couples the surface
and the longitudinal velocity. For $Fr>0$, positivity requires $\gamma_{\ell}>1$
together with a non-negative determinant, which gives a quadratic
inequality in $\gamma_{\ell}$ with no time step in it,

\begin{equation}
\bigl(4+4\sqrt{2}\,Fr+Fr^{2}\bigr)\gamma^{2}_{\ell}-\bigl(8+4\sqrt{2}\,Fr-2Fr^{2}\bigr)\gamma_{\ell}+4-12Fr^{2}\geq0\,.\label{eq:gamma_quadratic}
\end{equation}

At $\gamma_{\ell}=1$, the left-hand side equals $-9Fr^{2}<0$: the
value $\gamma_{\ell}=1$ lies strictly between the two roots, and
\eqref{eq:gamma_quadratic} together with $\gamma_{\ell}>1$ reduces
to $\gamma_{\ell}\geq\Gamma(Fr)$, the upper root,

\begin{equation}
\Gamma(Fr)=\frac{4+2\sqrt{2}\,Fr-Fr^{2}+Fr\sqrt{44+44\sqrt{2}\,Fr+13Fr^{2}}}{4+4\sqrt{2}\,Fr+Fr^{2}}\,.\label{eq:gamma_law}
\end{equation}

On a general mesh, the law is evaluated at the edge Froude number
$Fr_{e}=\lVert\overline{\boldsymbol{u}}_{e}\rVert/\sqrt{g\,\overline{h}_{e}}$,
built from the same edge average as $\theta_{e}$ in \eqref{eq:theta_mesh}.
The law is increasing, and its two limits are explicit,

\begin{equation}
\Gamma(0)=1\,,\qquad\Gamma(Fr)=1+\frac{\sqrt{11}-\sqrt{2}}{2}Fr+O(Fr^{2})\,,\qquad\Gamma(\infty)=\sqrt{13}-1\simeq2.6056\,.\label{eq:gamma_limits}
\end{equation}

\begin{figure}[htbp]
\begin{centering}
\includegraphics[width=0.95\textwidth]{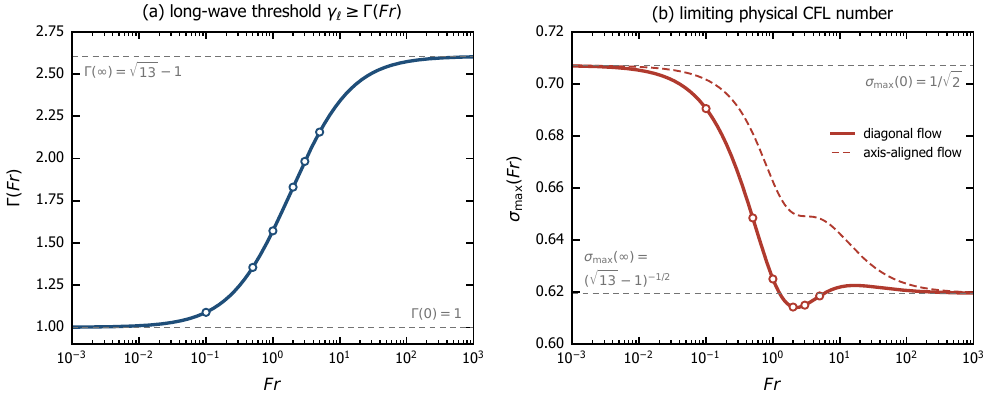}
\par\end{centering}
\caption{Linear calibration of the gravity step. (a) Long-wave threshold $\Gamma(Fr)$
of \eqref{eq:gamma_law}. (b) Limiting Courant number $\sigma_{\max}(Fr)$
of \eqref{eq:numax}, for a diagonal flow (solid) and an axis-aligned
flow (dashed); over the tested Froude range, the diagonal flow is
the most restrictive and gives $\sigma_{\max}$. Circles: tabulated
values; dashed horizontal lines: limits $Fr\to0$ and $Fr\to\infty$.}\label{fig:gamma_numax}
\end{figure}

The law, plotted in Fig.~\ref{fig:gamma_numax}(a), extends the rest-state
condition \eqref{eq:rest_conditions}, recovered at $Fr=0$. In particular,
as shown above, $\gamma_{\ell}=1$ is not admissible at any nonzero
Froude number with the calibrated weight. The fully upwind variant,
$\theta=1$, also used for comparison in Section~\ref{sec:validation},
has a larger transport diffusion, $q=1/\sigma_{0}$, which depends
on the time step and lowers the long-wave threshold below $\Gamma(Fr)$:
reducing $\theta$ to its minimum is what requires raising $\gamma_{\ell}$
to $\Gamma(Fr)$.
\begin{lem}
\label{lem:diag} Let $Fr>0$, $q=\sqrt{2}\,Fr$ and $F^{2}_{x}+F^{2}_{y}=Fr^{2}$.
Then $\mathbf{L}(\xi_{x},\xi_{y})\geq0$ for every flow orientation
satisfying the above constraint and every propagation direction if
and only if $\gamma_{\ell}\geq\Gamma(Fr)$.
\end{lem}
\begin{proof}
Necessity follows from the diagonal configuration \eqref{eq:worst_case},
which belongs to the configurations considered. For sufficiency, let
$\gamma_{\ell}\geq\Gamma(Fr)>1$ and, by homogeneity, $\lVert\xi\rVert=1$.
For $\boldsymbol{Z}=(h,u,v)\in\mathbb{R}^{3}$, let $u_{\parallel}=\xi_{x}u+\xi_{y}v$
be the velocity component along the wave. We seek a lower bound involving
only the height and the longitudinal velocity. Expanding \eqref{eq:Lmat}
gives
\[
\boldsymbol{Z}^{T}\mathbf{L}\,\boldsymbol{Z}=\gamma_{\ell}\Psi-(1+F^{2}_{\xi})h^{2}-4F_{\xi}hu_{\parallel}-u^{2}_{\parallel}+\bigl(q\,d_{F}-F^{2}_{\xi}\bigr)(u^{2}+v^{2}),
\]
\[
\Psi=(1+2d_{F})\,h^{2}+\xi^{2}_{x}(F_{x}h+u)^{2}+\xi^{2}_{y}(F_{y}h+v)^{2}+\xi^{2}_{x}u^{2}+\xi^{2}_{y}v^{2},
\]
where $\Psi$ gathers the directional diffusions. We bound them below
by their isotropic counterparts. By the Cauchy--Schwarz inequality
and $|F_{x}|+|F_{y}|\leq\sqrt{2}\,Fr$,
\[
F^{2}_{\xi}\leq\Bigl(\sum_{i}|F_{i}||\xi_{i}|\Bigr)^{2}\leq d_{F}\sum_{i}|F_{i}|\leq\sqrt{2}\,Fr\,d_{F},
\]
and
\[
\xi^{2}_{x}u^{2}+\xi^{2}_{y}v^{2}\geq\frac{1}{2}u^{2}_{\parallel},\qquad\xi^{2}_{x}(F_{x}h+u)^{2}+\xi^{2}_{y}(F_{y}h+v)^{2}\geq\frac{1}{2}(F_{\xi}h+u_{\parallel})^{2}.
\]
The first inequality gives $q\,d_{F}\geq F^{2}_{\xi}$ and $2d_{F}\geq\sqrt{2}\,F^{2}_{\xi}/Fr$.
Hence $\Psi\geq\Psi_{\parallel}:=\bigl(1+\sqrt{2}\,F^{2}_{\xi}/Fr\bigr)h^{2}+(F_{\xi}h+u_{\parallel})^{2}/2+u^{2}_{\parallel}/2$,
and, since $\gamma_{\ell}>0$ and $q\,d_{F}\geq F^{2}_{\xi}$,
\[
\boldsymbol{Z}^{T}\mathbf{L}\,\boldsymbol{Z}\geq(h,u_{\parallel})\,\mathbf{L}_{\parallel}(F_{\xi})\,(h,u_{\parallel})^{T},\qquad\mathbf{L}_{\parallel}(F_{\xi})=\begin{pmatrix}\gamma_{\ell}\Bigl(1+\Bigl(\dfrac{\sqrt{2}}{Fr}+\dfrac{1}{2}\Bigr)F^{2}_{\xi}\Bigr)-1-F^{2}_{\xi} & \dfrac{1}{2}F_{\xi}(\gamma_{\ell}-4)\\[3pt]
\dfrac{1}{2}F_{\xi}(\gamma_{\ell}-4) & \gamma_{\ell}-1
\end{pmatrix},
\]
which depends on the orientation only through $F_{\xi}$, with $|F_{\xi}|\leq Fr$
since $\lVert\xi\rVert=1$. Since $\gamma_{\ell}-1>0$, $\mathbf{L}_{\parallel}(F_{\xi})\geq0$
if and only if its determinant is non-negative, and this determinant
is affine in $F^{2}_{\xi}\in[0,Fr^{2}]$. It equals $(\gamma_{\ell}-1)^{2}$
at $F_{\xi}=0$, and at $F^{2}_{\xi}=Fr^{2}$, $\mathbf{L}_{\parallel}(Fr)$
is the block \eqref{eq:Lpar}, since in the diagonal configuration
the lower bound differs from the original quadratic form only by the
non-negative transverse contribution $\gamma_{\ell}u^{2}_{\perp}$;
its determinant is non-negative for $\gamma_{\ell}\geq\Gamma(Fr)$.
Hence $\mathbf{L}_{\parallel}(F_{\xi})\geq0$ for all admissible $F_{\xi}$,
and $\mathbf{L}\geq0$.
\end{proof}

The three inequalities are simultaneously saturated when the flow
and the wave are aligned with the same mesh diagonal and the velocity
perturbation is longitudinal: this is where the directional diffusions,
built axis by axis, are weakest compared with their isotropic counterparts.
The first one also explains the transport law: $q=\sqrt{2}\,Fr$ is
precisely the value for which the upwind diffusion $q\,d_{F}$ dominates
the advective term $F^{2}_{\xi}$ in every orientation, so that transport
and gravity share the same worst case. The dependence of the threshold
on the orientation is thus a genuine defect of the Cartesian discretization,
which is not rotation invariant. Calibrating on the worst orientation
makes the law independent of orientation, at the price of some excess
dissipation for the other orientations. This motivates its extension
to unstructured meshes, although the present proof is restricted to
Cartesian grids. As at rest, we retain the smallest admissible value,
$\gamma_{\ell}=\Gamma(Fr)$, since a larger constant increases the
stabilization coefficients and lowers the admissible time step (Step
three).

\paragraph{Step three: the resulting time-step restriction}

The long-wave condition says nothing about finite frequencies, which
is where the time step enters. Over the sampled Fourier modes, flow
orientations and Froude numbers, the numerical sweep shows that with
$\gamma_{\ell}=\Gamma(Fr)$ and $\theta$ from \eqref{eq:theta_law}
the most restrictive mode is the checkerboard mode $(\pi,\pi)$, at
which the restriction can be computed exactly. There $s_{x}=s_{y}=0$
and $c_{x}=c_{y}=2$, so that $\tilde{\mathbf{A}}^{h\boldsymbol{u}}_{g}$
is real and lower triangular,

\[
\tilde{\mathbf{A}}^{h\boldsymbol{u}}_{g}(\pi,\pi)=\begin{pmatrix}\lambda_{0} & 0 & 0\\
-4\gamma_{\ell}\sigma^{2}_{0}F_{x} & \lambda_{1} & 0\\
-4\gamma_{\ell}\sigma^{2}_{0}F_{y} & 0 & \lambda_{1}
\end{pmatrix},\qquad\lambda_{0}=1-2\gamma_{\ell}(\sigma^{2}_{x}+\sigma^{2}_{y}),\qquad\lambda_{1}=1-4\gamma_{\ell}\sigma^{2}_{0}-2\theta\lVert\tilde{\sigma}\rVert_{1}.
\]

Here $\lambda_{0},\lambda_{1}$ denote the diagonal entries, and hence
the eigenvalues, of the amplification matrix. The velocity block is
invariant under rotation, and the coupling is parallel to the flow.
In coordinates aligned with the flow, $\mathbf{B}^{h\boldsymbol{u}}(\pi,\pi)$
therefore splits into the transverse eigenvalue $1-\lambda^{2}_{1}$
and the longitudinal block

\[
\begin{pmatrix}1-\lambda^{2}_{0}-16\gamma^{2}_{\ell}\sigma^{4}_{0}Fr^{2} & 4\gamma_{\ell}\sigma^{2}_{0}Fr\,\lambda_{1}\\
4\gamma_{\ell}\sigma^{2}_{0}Fr\,\lambda_{1} & 1-\lambda^{2}_{1}
\end{pmatrix}.
\]

Its positive semidefiniteness is equivalent to

\begin{equation}
|\lambda_{0}|\leq1,\qquad|\lambda_{1}|\leq1,\qquad(1-\lambda^{2}_{0})(1-\lambda^{2}_{1})\geq16\gamma^{2}_{\ell}\sigma^{4}_{0}Fr^{2}.\label{eq:checker_cond}
\end{equation}

The orientation enters only through $\lVert\boldsymbol{F}\rVert_{1}=|F_{x}|+|F_{y}|\in[Fr,\sqrt{2}\,Fr]$:
using the uncapped weight, whose consistency is verified below, $\theta=\sqrt{2}\,\sigma_{0}Fr$,
one has $\lambda_{0}=1-4\rho_{0}\sigma^{2}_{0}$ and $\lambda_{1}=1-4\rho_{1}\sigma^{2}_{0}$,
where $\rho_{0}=\gamma_{\ell}\bigl(1+\lVert\boldsymbol{F}\rVert_{1}+Fr^{2}/2\bigr)$
and $\rho_{1}=\gamma_{\ell}+Fr\,\lVert\boldsymbol{F}\rVert_{1}/\sqrt{2}$.
For $\sigma_{0}>0$, the determinant condition in \eqref{eq:checker_cond}
reads
\[
4\rho_{0}\rho_{1}\bigl(1-2\rho_{0}\sigma^{2}_{0}\bigr)\bigl(1-2\rho_{1}\sigma^{2}_{0}\bigr)\geq\gamma^{2}_{\ell}Fr^{2},
\]
a quadratic inequality in $\sigma^{2}_{0}$. The largest range of
time steps from $\Delta t=0$ on which \eqref{eq:checker_cond} holds
is

\begin{equation}
\sigma^{2}_{0}\leq\sigma^{2}_{*}:=\frac{\rho_{0}+\rho_{1}-\sqrt{(\rho_{0}-\rho_{1})^{2}+\gamma^{2}_{\ell}Fr^{2}}}{4\,\rho_{0}\rho_{1}}\,.\label{eq:sigma_p}
\end{equation}

At fixed $\sigma_{0}$, $\lambda_{0}$ and $\lambda_{1}$ are affine
in $\lVert\boldsymbol{F}\rVert_{1}$, so that $1-\lambda^{2}_{0}$
and $1-\lambda^{2}_{1}$ are concave; non-negative at both ends of
$[Fr,\sqrt{2}\,Fr]$, they remain so in between, and their product,
being log-concave, attains its minimum at an end. Condition \eqref{eq:checker_cond}
therefore holds for every flow orientation if and only if it holds
for the axis-aligned flow, $\lVert\boldsymbol{F}\rVert_{1}=Fr$, and
for the diagonal flow, $\lVert\boldsymbol{F}\rVert_{1}=\sqrt{2}\,Fr$,
that is, if and only if $\sigma_{0}\leq\min(\sigma_{\text{axis}},\sigma_{\text{diag}})$,
where $\sigma_{\text{axis}}$ and $\sigma_{\text{diag}}$ denote the
values of $\sigma_{*}$ at these two ends. Numerically, the diagonal
bound is the smaller one over the tested Froude range; we keep the
minimum so as not to rely on this observation. In terms of $\sigma_{g}$,
the condition reads $\sigma_{g}\leq\sigma_{\max}(Fr)$ with

\begin{equation}
\begin{gathered}\sigma_{\max}(Fr)=(1+Fr)\min\!\left(\sigma_{\text{axis}},\sigma_{\text{diag}}\right),\\
\rho_{0,\text{diag}}=\Gamma(Fr)\Bigl(1+\frac{Fr}{\sqrt{2}}\Bigr)^{2},\qquad\rho_{1,\text{diag}}=\Gamma(Fr)+Fr^{2}\,.
\end{gathered}
\label{eq:numax}
\end{equation}

This bound keeps the transport weight below its cap. For a diagonal
flow, $\lambda_{1}=1-4(\gamma_{\ell}+Fr^{2})\sigma^{2}_{0}$, and
$|\lambda_{1}|\leq1$ gives $2Fr^{2}\sigma^{2}_{0}\leq Fr^{2}/(\gamma_{\ell}+Fr^{2})<1$,
hence $\lVert\tilde{\sigma}\rVert_{1}\leq\theta=\sqrt{2}\,Fr\,\sigma_{0}<1$:
the uncapped weight is consistent, and the assumptions of Lemma~\ref{lem:awp}
hold.

The resulting values are

\[
\begin{array}{c|ccccccc}
Fr & 0 & 0.1 & 0.5 & 1 & 2 & 3 & 5\\
\hline \Gamma(Fr) & 1.0000 & 1.0889 & 1.3549 & 1.5715 & 1.8304 & 1.9827 & 2.1566\\
\sigma_{\max} & 0.7071 & 0.6905 & 0.6485 & 0.6251 & 0.6143 & 0.6150 & 0.6185
\end{array}\,.
\]

The limits are explicit, $\sigma_{\max}(0)=1/\sqrt{2}$, which is
the rest-state condition, sufficient over all frequencies by Lemma~\ref{lem:rest},
and $\sigma_{\max}(\infty)=(\sqrt{13}-1)^{-1/2}\simeq0.6195$; the
minimum, $\sigma_{\max}\simeq0.614$, is reached near $Fr\simeq2.2$.
The law therefore removes the time-step dependence of the long-wave
threshold without a severe CFL penalty. Both laws are plotted in Fig.~\ref{fig:gamma_numax}.

These results hold on a Cartesian grid. On a general mesh, the prescriptions
\eqref{eq:theta_mesh} and $\gamma_{\ell,e}=\Gamma(Fr_{e})$ are extensions,
to be checked on the mesh actually used through the fully discrete
energy balance.

\begin{figure}[htbp]
\begin{centering}
\includegraphics[width=0.5\textwidth]{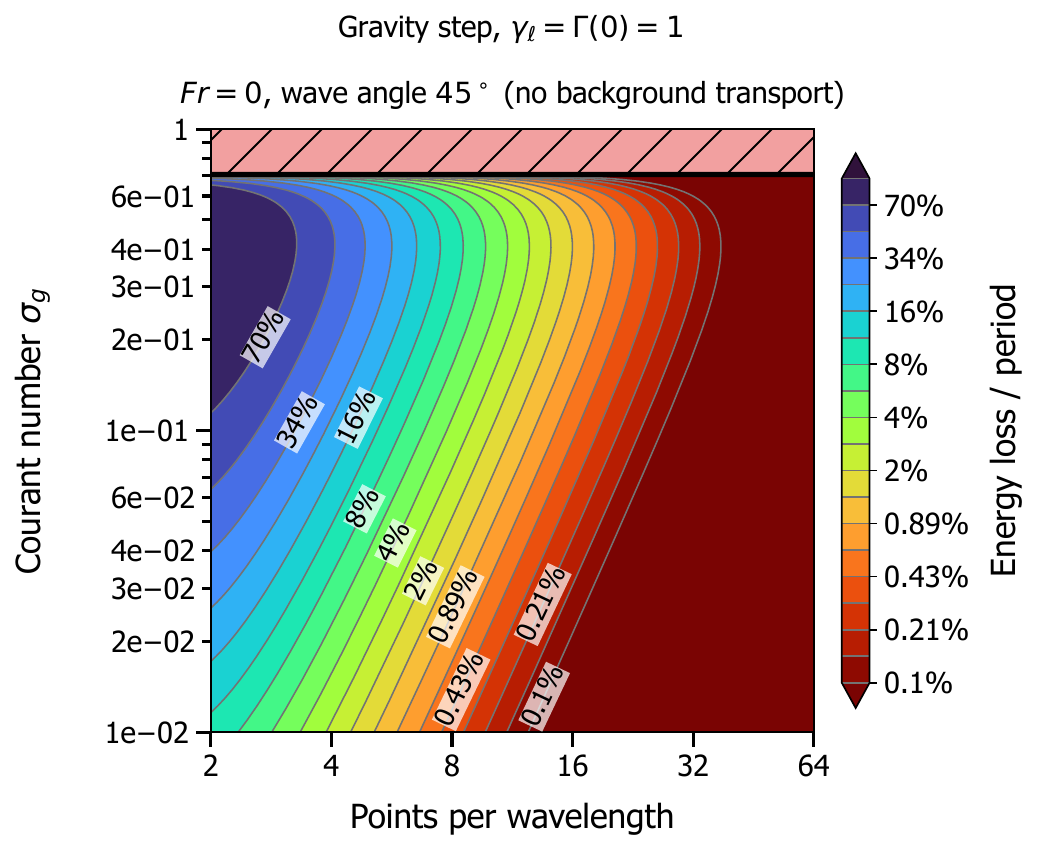}\includegraphics[width=0.5\textwidth]{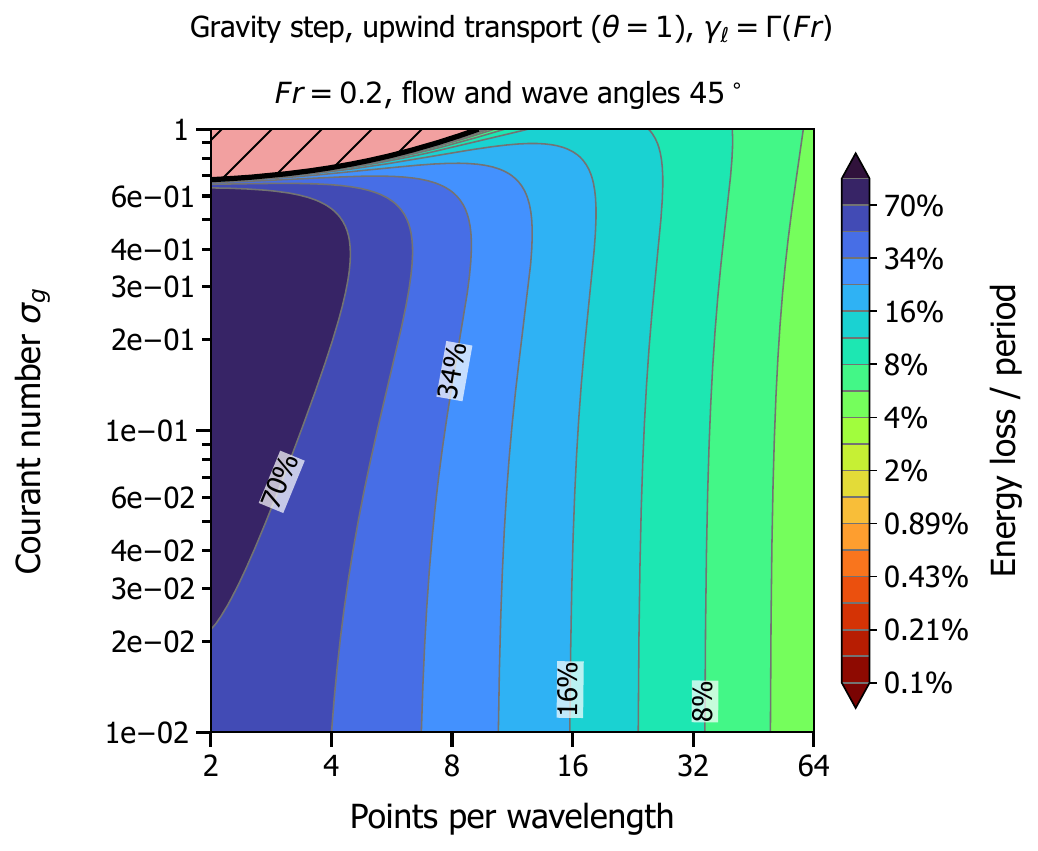}
\par\end{centering}
\begin{centering}
\includegraphics[width=0.5\textwidth]{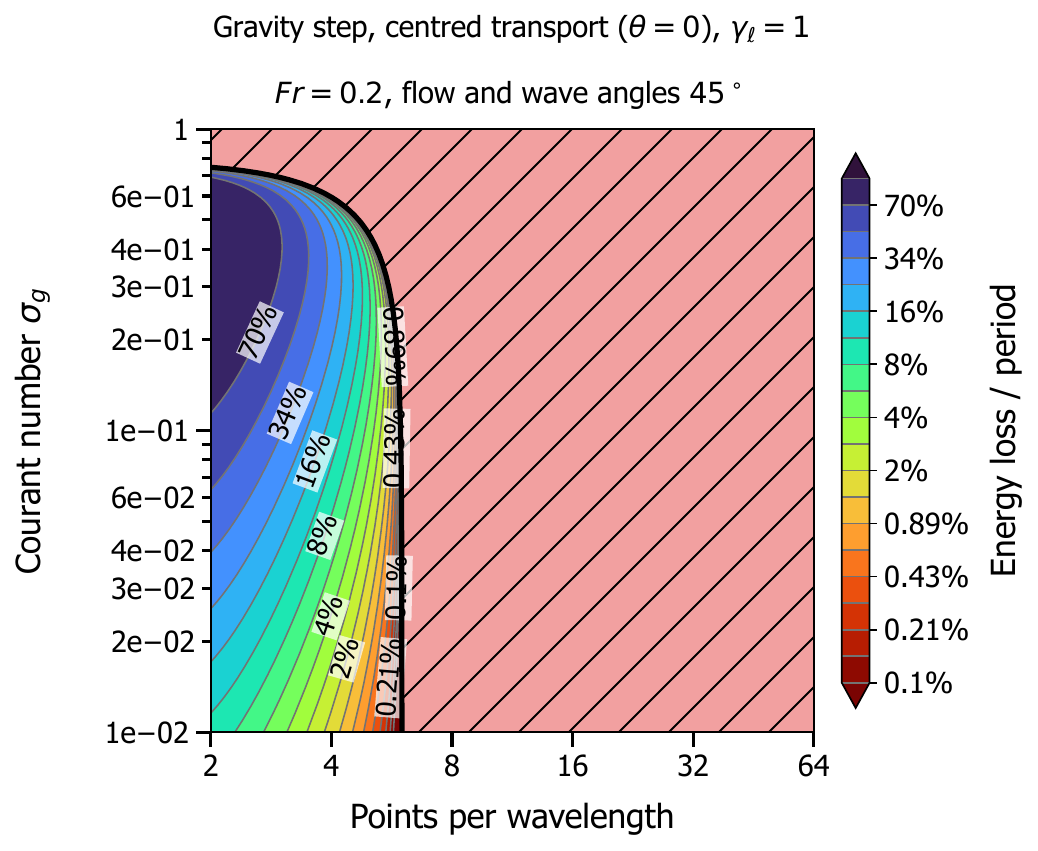}\includegraphics[width=0.5\textwidth]{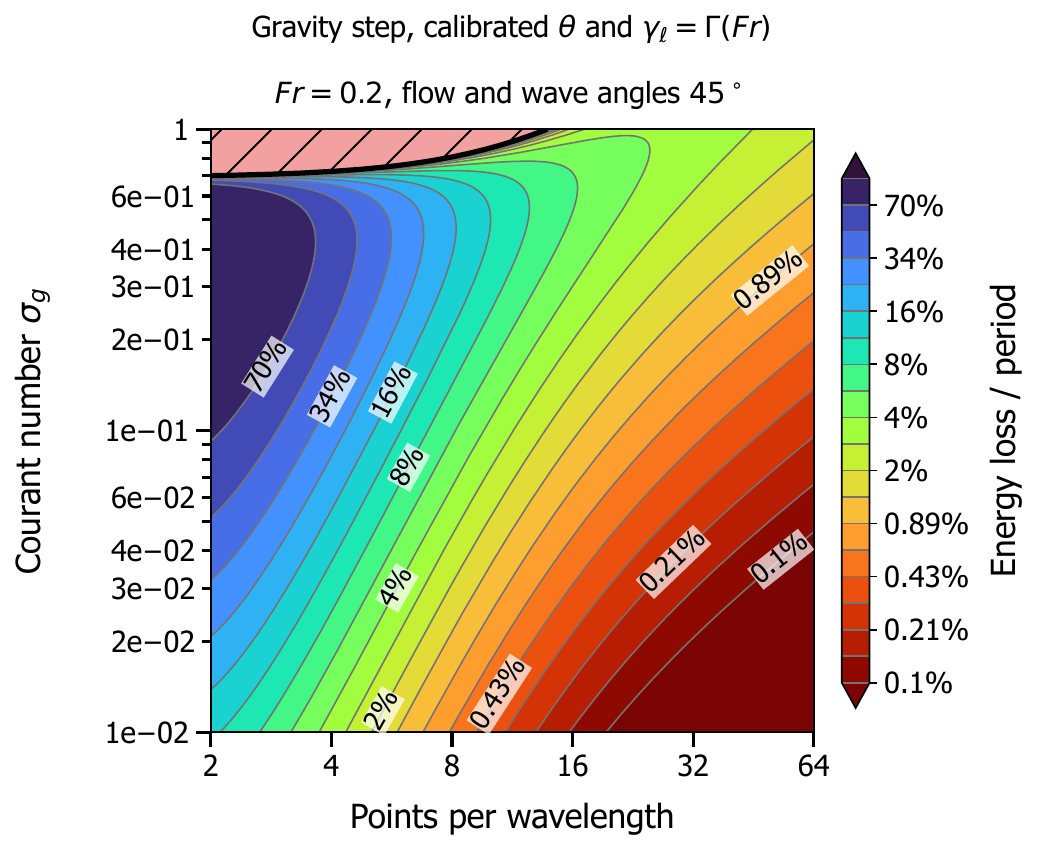}
\par\end{centering}
\begin{centering}
\includegraphics[width=0.5\textwidth]{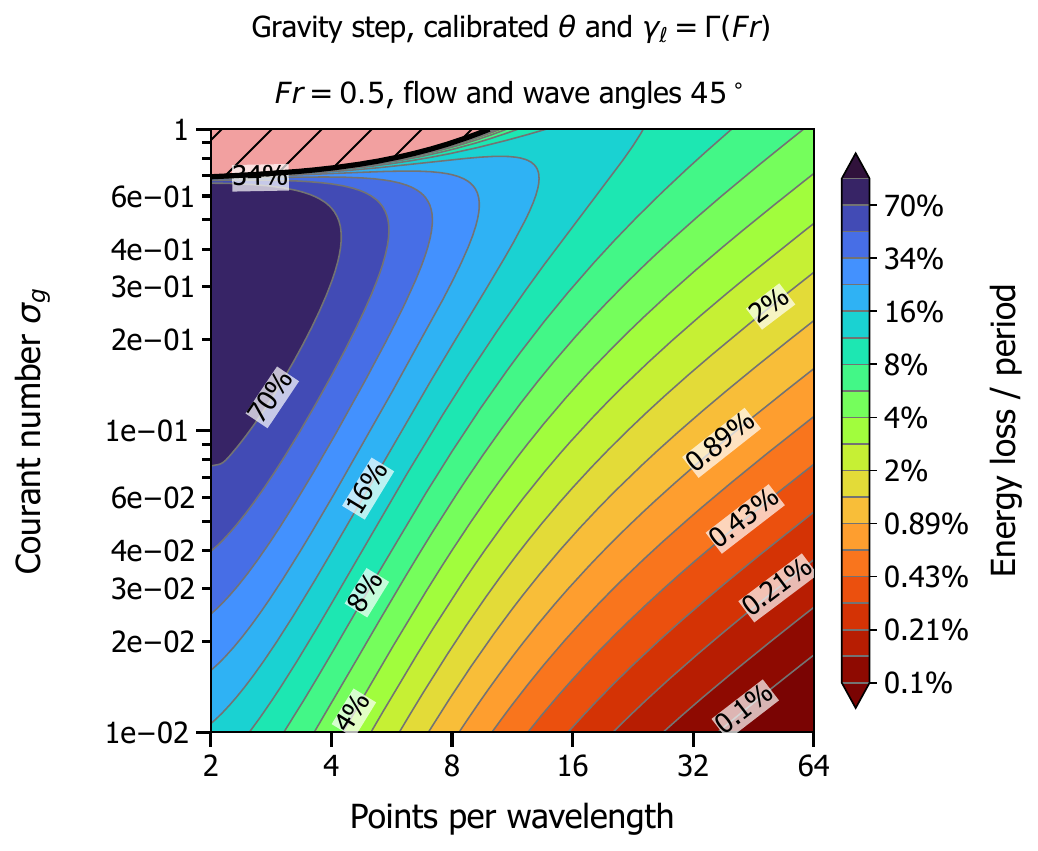}\includegraphics[width=0.5\textwidth]{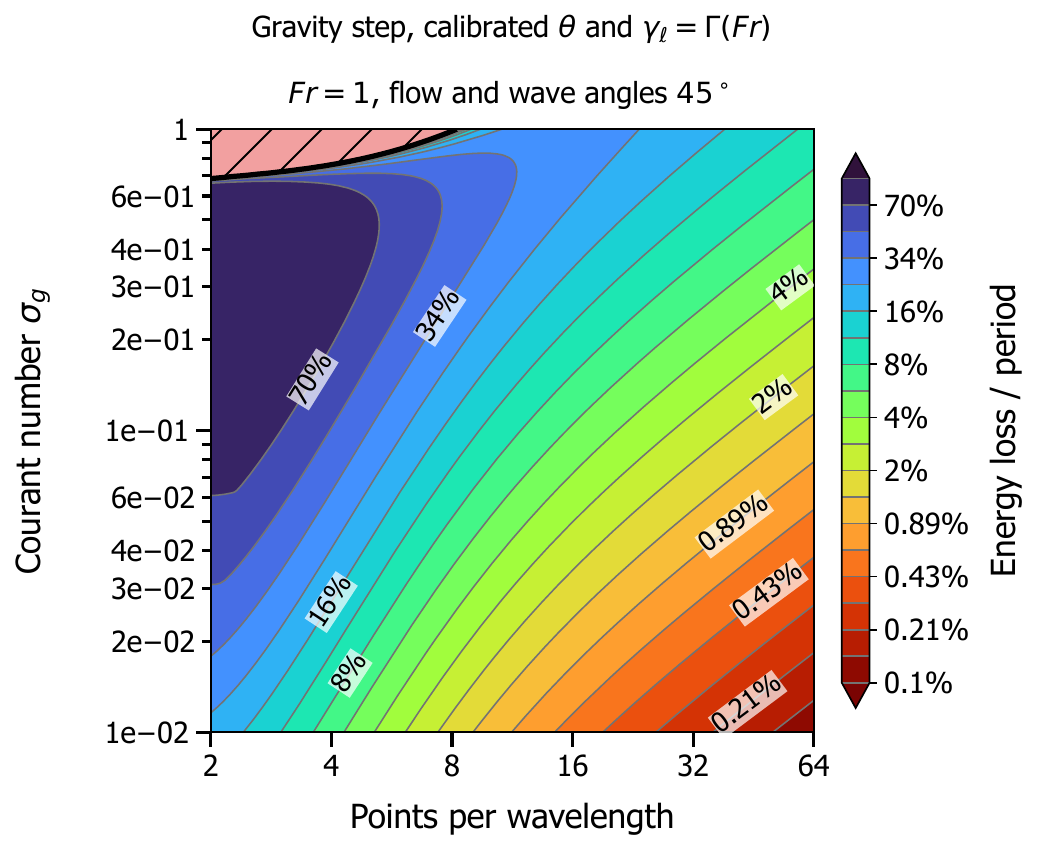}
\par\end{centering}
\caption{Energy loss per intrinsic period for the gravity-wave mode, as a function
of $\sigma_{g}$ and $L/\Delta$, with flow and wave angles of $45^{\circ}$.
Top row: rest state (left), fully upwind transport at $Fr=0.2$ with
$\gamma_{\ell}=\Gamma(Fr)$ (right). Middle row: centred transport
at $Fr=0.2$ with $\gamma_{\ell}=1$ (left), calibrated transport
and viscosity at $Fr=0.2$ (right). Bottom row: calibrated transport
and viscosity at $Fr=0.5$ (left) and $Fr=1$ (right). Hatched regions
denote energy growth.}\label{fig:dissip-sw_base-Fr-0}
\end{figure}

To assess the numerical dissipation, Fig.~\ref{fig:dissip-sw_base-Fr-0}
uses the gravity-step Courant number $\sigma_{g}$ as its vertical
coordinate and the number of points per wavelength $L/\Delta$ as
its horizontal coordinate. This convention differs from the prescribed
parameters used in the implementation; their relation is given in
Remark \ref{rem:courant_conventions}. The plotted quantity is the
geometric energy loss per intrinsic wave period, evaluated over approximately
$100$ periods. At fixed $\sigma_{g}$, increasing $Fr$ decreases
the time step.

The rest-state panel provides a common reference. At $Fr=0.2$, comparing
the fully upwind and calibrated panels at the same $\gamma_{\ell}=\Gamma(Fr)$
isolates the reduction in dissipation obtained by adapting the transport
weight. The centred case with $\gamma_{\ell}=1$ instead exhibits
energy growth over a broad part of the map; since both coefficients
differ, it is a diagnostic rather than a comparison of the transport
weight alone. The last two panels illustrate the calibrated scheme
at $Fr=0.5$ and $Fr=1$. The thick black curve separates energy loss
from energy growth for the selected wave over the simulated duration;
it is not the uniform contractivity bound $\sigma_{g}\leq\sigma_{\max}(Fr)$
of \eqref{eq:numax}.

\begin{prop}[Linear energy stability of the gravity step]
\label{prop:linear_sw} Consider the linearized scheme \eqref{eq:scheme_sw_lin}
on a uniform Cartesian grid, about a uniform state $(H,U,V,0,0)$
with $c_{0}=\sqrt{gH}$, $\boldsymbol{F}=(U,V)/c_{0}$, $Fr=\lVert\boldsymbol{F}\rVert$,
$\sigma_{0}=c_{0}\Delta t/\Delta>0$, a constant weight $\theta\in[0,1]$
and the family \eqref{eq:linear_viscosity_choice}. The scheme is
energy-contractive when $\mathbf{B}(\xi_{x},\xi_{y})\geq0$ for every
$(\xi_{x},\xi_{y})\in[-\pi,\pi]^{2}$, and then $L^{2}$-stable. Then:

(i)\emph{ Rest state} ($Fr=0$): the scheme is energy-contractive
if and only if $\gamma_{\ell}\geq1$ and $\gamma_{\ell}\sigma^{2}_{0}\leq1/2$.

(ii)\emph{ Transport}: if $\lVert\tilde{\sigma}\rVert_{1}=\sigma_{0}(|F_{x}|+|F_{y}|)\leq1$,
the passive scalars and the transverse modes are contractive at every
frequency if and only if $\theta\geq\lVert\tilde{\sigma}\rVert_{1}$.
In particular, the centred variant $\theta=0$ is unstable about any
moving state, and $\theta=\sqrt{2}\,Fr\,\sigma_{0}$ is the smallest
weight depending only on $Fr$ that is admissible for every flow orientation.

(iii)\emph{ Long waves}: for $Fr>0$ and $\theta=\sqrt{2}\,Fr\,\sigma_{0}$,
the long-wave matrix $\mathbf{L}$ of \eqref{eq:Lmat} is nonnegative
for every flow orientation and every propagation direction if and
only if $\gamma_{\ell}\geq\Gamma(Fr)$, with $\Gamma$ given by \eqref{eq:gamma_law}.
This condition, independent of the time step, is necessary for energy
contractivity.

(iv)\emph{ Checkerboard mode}: for $\theta=\sqrt{2}\,Fr\,\sigma_{0}$
and $\gamma_{\ell}=\Gamma(Fr)$, $\mathbf{B}(\pi,\pi)\geq0$ for every
flow orientation if and only if $\sigma_{g}\leq\sigma_{\max}(Fr)$,
with $\sigma_{\max}$ given by \eqref{eq:numax}. This condition is
necessary for energy contractivity.

The sufficiency of (iii)--(iv) over the whole Fourier domain in the
presence of a current is not proved; it is supported by the numerical
sweep of Step three.
\end{prop}
\begin{proof}
(i) is Lemma~\ref{lem:rest}; (ii) is Lemma~\ref{lem:awp}, together
with \eqref{eq:Gperp} for the transverse modes and the bound $|F_{x}|+|F_{y}|\leq\sqrt{2}\,Fr$,
attained for a diagonal flow; (iii) is Lemma~\ref{lem:diag}; (iv)
is the checkerboard computation of Step three.
\end{proof}

These results concern the calibrated scheme, not the sufficient nonlinear
coefficient interval \eqref{eq:rbound}.

\subsubsection{Nonlinear stability: gravity step}

\paragraph{Global energy balance}

The linear analysis calibrates the stabilization around uniform states.
We now examine the total discrete energy of the nonlinear gravity
step, 
\begin{equation}
\mathcal{E}^{n}:=\sum_{K}m_{K}h^{n}_{K}E^{n}_{K},\label{eq:gravity-total-energy}
\end{equation}
where $E^{n}_{K}$ is the specific energy defined above. We retain
the edge notations introduced previously. Assume periodic boundaries
and positive depths at both time levels. In the sums below, each edge
is counted once, and $\Delta z_{K}=z^{n+1}_{K}-z^{n}_{K}$. In particular,
$\Delta(hz)_{K}=h^{n+1}_{K}z^{n+1}_{K}-h^{n}_{K}z^{n}_{K}$.

The scheme satisfies the exact identity (proved in~\ref{app:gravity-global-balance})
\begin{equation}
\begin{gathered}\mathcal{E}^{n+1}-\mathcal{E}^{n}=-2\Delta t\sum_{e}m_{e}\Bigl[(\boldsymbol{\Pi}^{n}_{e}\cdot\boldsymbol{n}_{e})\llbracket\phi^{n}\rrbracket_{e}+\Lambda^{n}_{e}\llbracket h^{n}\boldsymbol{u}^{n}\rrbracket_{e}\cdot\boldsymbol{n}_{e}\Bigr]\\
-2\Delta t\sum_{e}m_{e}\theta^{n}_{e}|\mathcal{F}^{*,n}_{e}|\left(\|\llbracket\boldsymbol{u}^{n}\rrbracket_{e}\|^{2}+\frac{4}{3}\llbracket w^{n}\rrbracket^{2}_{e}+\frac{1}{a^{2}}\llbracket p^{n}\rrbracket^{2}_{e}\right)+\mathcal{R}_{\Delta t},
\end{gathered}
\label{eq:gravity-global-balance}
\end{equation}
where all edge quantities are evaluated at $t^{n}$ and 
\begin{equation}
\mathcal{R}_{\Delta t}=\sum_{K}m_{K}\left[\frac{g}{2}(\Delta h_{K})^{2}+\frac{h^{n+1}_{K}}{2}\left(\|\Delta\boldsymbol{u}_{K}\|^{2}+\frac{4}{3}(\Delta w_{K})^{2}+\frac{1}{a^{2}}(\Delta p_{K})^{2}\right)\right]\geq0.\label{eq:gravity-time-remainder}
\end{equation}
The prescribed $\boldsymbol{\Pi}^{n}_{e}\cdot\boldsymbol{n}_{e}$
and $\Lambda^{n}_{e}$ are proportional, with positive coefficients,
to $\llbracket\phi^{n}\rrbracket_{e}$ and $\llbracket h^{n}\boldsymbol{u}^{n}\rrbracket_{e}\cdot\boldsymbol{n}_{e}$,
respectively. Their contributions in \eqref{eq:gravity-global-balance}
are therefore non-positive, as is the transport term for $\theta^{n}_{e}\geq0$.
Their factors $\Delta t\,r_{K}$ give contributions of order $\Delta t^{2}/\Delta^{2}$
per unit volume on a regular mesh, matching the quadratic temporal
remainder. Energy decay requires these contributions, together with
transport dissipation, to dominate $\mathcal{R}_{\Delta t}$. The
following sufficient estimate ensures this for a fully upwind, cellwise
prescription.

\paragraph{A sufficient nonlinear estimate}

\label{subsec:nonlinear-sw}  We now address the stability conditions
of the first part of the splitting in the nonlinear framework with
upwind fluxes, meaning $\theta_{e}=1$ in \eqref{eq:scheme_sw}. Here
the sufficient nonlinear estimate uses the cellwise prescription $\nu_{K,e}=\nu_{K}$
in the mass-flux stabilization.  The fully upwind specialization of
the mechanical-energy estimate established in \citep[Remark 8]{duran2026},
together with the compatible passive-scalar estimate of \citep[Section 5.2]{duran2026}
applied to $w$ and $p$, yields the following result. Denoting $E^{n}_{K}$
the local discrete energy arising from formula \eqref{eq:energyE},
the scheme \eqref{eq:scheme_sw} satisfies a local discrete entropy
inequality of the form 
\begin{equation}
{\displaystyle (hE)^{n+1}_{K}\leq(hE)^{n}_{K}-\dfrac{\Delta t}{m_{K}}\sum_{e\in\partial K}\mathcal{G}_{g}(W^{n}_{K},W^{n}_{K_{e}},\boldsymbol{n}_{e,K})m_{e}\,,}\label{eq:Energy_estimSW_final}
\end{equation}
with $W^{n}_{K}=(h,h\boldsymbol{u},hw,hp)^{n}_{K}$ and a conservative
numerical entropy flux function $\mathcal{G}_{g}$, consistent with
the energy flux given by \eqref{eq:energy_SW}; $\mathcal{G}_{g}$
also depends on the bottom, on the mesh geometry and, for the local
prescription below, on coefficients computed on the neighbourhood
of the edge. Accordingly, for the uniform choice $\xi_{K}=\xi$, the
fully explicit CPR scheme is stable provided that, for every cell
$K$ and for every edge $e\in\partial K$, 
\begin{align}
\Delta t\,\max\!\left(\frac{m_{\partial K}}{m_{K}},\frac{m_{\partial K_{e}}}{m_{K_{e}}}\right)\frac{\lvert\mathcal{F}^{\ast}_{e,K}\rvert}{\widecheck{h}_{e}} & \leq\frac{\xi_{K}}{1+\xi_{K}},\label{eq:cpr-cfl-transport}\\
\Delta t\,\max\!\left(\frac{m_{\partial K}}{m_{K}},\frac{m_{\partial K_{e}}}{m_{K_{e}}}\right)\sqrt{g\overline{h}_{e}} & \leq\frac{1}{2\sqrt{\nu_{K}}},\label{eq:cpr-cfl-gravity}
\end{align}
where $\widecheck{h}_{e}=\min(h^{n}_{K},h^{n}_{K_{e}})$, $0<\xi_{K}<1$
and $\nu_{K}=(1+\xi_{K})/(1-\xi_{K})$ in the stabilization term \eqref{eq:Pie}.
The constants $\alpha,\gamma$ have to satisfy the following local
bounds: 
\begin{equation}
\gamma,\alpha\in[r^{-},r^{+}]\quad,\quad\text{with }r^{\pm}(\vartheta)=\dfrac{1\pm\sqrt{1-\vartheta}}{\vartheta/2}\,,\label{eq:rbound}
\end{equation}
where $\vartheta=4g\,\Delta t^{2}\,\overline{r}_{e}\,\overline{rh\nu}_{e}$.
The two CFL conditions separately constrain the transport and gravity-wave
contributions.

\paragraph{Local amplification and explicit time-step selection}

 Following \citet{duran2026}, we introduce a shallow-water-type gravity
restriction, supplemented below by an explicit transport restriction,
fixing the local parameters $\xi_{K}$ and $\nu_{K}$. For an edge
$e=K|K_{e}$, set 
\[
\dfrac{1}{\Delta_{e}}:=\max\!\left(\frac{m_{\partial K}}{m_{K}},\frac{m_{\partial K_{e}}}{m_{K_{e}}}\right),\qquad\widehat{\mathcal{F}}_{e,K}:=\frac{1}{2}\left(h^{n}_{K}\boldsymbol{u}^{n}_{K}+h^{n}_{K_{e}}\boldsymbol{u}^{n}_{K_{e}}\right)\cdot\boldsymbol{n}_{e,K}.
\]
Define 
\begin{align}
u^{\Delta}_{K} & :=\max_{e\in\partial K}\frac{|\mathcal{F}^{\ast}_{e,K}|}{\Delta_{e}\widecheck h_{e}}, & c^{\Delta}_{K} & :=\max_{e\in\partial K}\dfrac{\sqrt{g\overline{h}_{e}}}{\Delta_{e}},\nonumber \\
\widehat{u}^{\Delta}_{K} & :=\max_{e\in\partial K}\frac{|\widehat{\mathcal{F}}_{e,K}|}{\Delta_{e}\widecheck h_{e}}, & \mathscr{G}_{K} & :=\max_{e=K|K_{e}}\frac{\gamma}{4\Delta_{e}\widecheck h_{e}}\left(\frac{m_{\partial K}}{m_{K}}h^{n}_{K}+\frac{m_{\partial K_{e}}}{m_{K_{e}}}h^{n}_{K_{e}}\right)|\phi^{n}_{K_{e}}-\phi^{n}_{K}|.
\end{align}
Finally, let 
\[
T_{0}:=\min_{K}\frac{1}{2c^{\Delta}_{K}},\qquad\mathscr{H}_{K}:=\sqrt{c^{\Delta}_{K}T_{0}\mathscr{G}_{K}}.
\]
Note that these quantities are computed from the state at time $t^{n}$,
before selecting $\Delta t$. From this, introduce the local Froude
indicator 
\begin{equation}
Fr_{K}:=\frac{\mathfrak{u}^{\Delta}_{K}}{c^{\Delta}_{K}},\qquad\mathfrak{u}^{\Delta}_{K}:=\widehat{u}^{\Delta}_{K}+\mathscr{H}_{K}\,,\label{eq:frk}
\end{equation}
and choose 
\begin{equation}
\nu_{K}=(1+Fr_{K})^{2}.\label{eq:nuk}
\end{equation}
With the local choices above and $\gamma=\alpha=2$, set $\xi_{K}=(\nu_{K}-1)/(\nu_{K}+1)$.
The gravity restriction is ensured under the shallow-water-type CFL
condition
\begin{equation}
\Delta t\,\bigl(c^{\Delta}_{K}+\mathfrak{u}^{\Delta}_{K}\bigr)\leq\frac{1}{2}\qquad\text{for every cell }K.\label{eq:CFL_hyperbolic}
\end{equation}
This corresponds to the first branch of the explicit CFL condition
in \citep[Section~4.4]{duran2026}.

\begin{figure}[htbp]
\centering{}\includegraphics[width=0.6\textwidth]{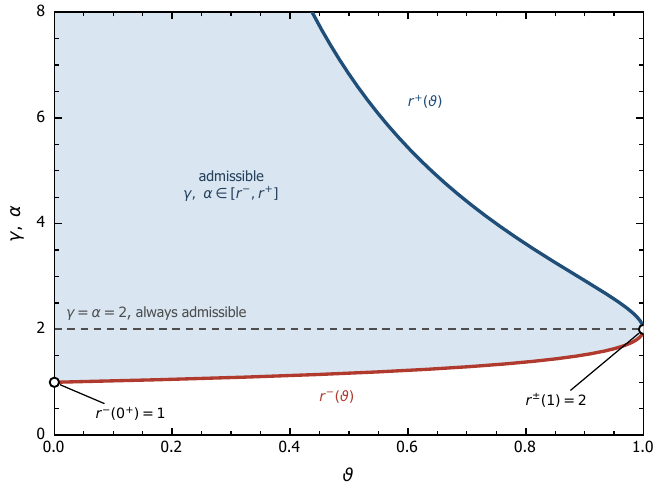} \caption{Nonlinear stability. Admissibility region for the viscosity constants
$\gamma,\alpha$ as a function of $\vartheta=4g\,\Delta t^{2}\,\overline{r}_{e}\,\overline{rh\nu}_{e}$,
arising from \eqref{eq:rbound}.}\label{fig:carto}
\end{figure}

To ensure the transport restriction as well, impose 
\begin{equation}
\Delta t\widehat{u}^{\Delta}_{K}+\Delta t^{2}\mathscr{P}_{K}\leq\frac{\xi_{K}}{1+\xi_{K}},\label{eq:cpr-transport-explicit}
\end{equation}
where 
\[
\mathscr{P}_{K}=\max_{e=K|K_{e}}\frac{\gamma(\nu_{K}r_{K}h^{n}_{K}+\nu_{K_{e}}r_{K_{e}}h^{n}_{K_{e}})}{4\Delta_{e}\widecheck h_{e}}|\phi^{n}_{K_{e}}-\phi^{n}_{K}|.
\]
Indeed, $u^{\Delta}_{K}\leq\widehat{u}^{\Delta}_{K}+\Delta t\mathscr{P}_{K}$.
This local construction uses the cited local result, separately from
the uniform-parameter statement above. 
\begin{rem}
If $Fr_{K}=0$, then $\widehat{u}^{\Delta}_{K}=\mathscr{G}_{K}=0$.
The centred normal flux and the potential jump therefore vanish on
every edge of $K$. Hence $\mathcal{F}^{\ast}_{e,K}=0$ on these edges,
and the transport condition \eqref{eq:cpr-cfl-transport} holds directly,
although $\xi_{K}=0$. The CFL condition \eqref{eq:CFL_hyperbolic}
remains well defined and reduces locally to $\Delta t\,c^{\Delta}_{K}\leq1/2$. 
\end{rem}
\begin{rem}
\label{rem:constants} We have $\lim\limits_{\vartheta\rightarrow0}r^{-}(\vartheta)=1$,
$r^{\pm}(1)=2$ and $\lim\limits_{\vartheta\rightarrow0}r^{+}(\vartheta)=+\infty$.
The function $\vartheta\mapsto r^{-}(\vartheta)$ being increasing
and $\vartheta\mapsto r^{+}(\vartheta)$ decreasing on the interval
$[0,1]$, we get that the value $\gamma=\alpha=2$ is always admissible,
independently of the edge $e$ and of the neighbouring cells (see
Fig. \ref{fig:carto}). Note also the possibility of taking $\gamma,\alpha$
close to $1$ by diminishing $\vartheta$ (that is, the time step).
The common admissibility interval $[r^{-}(\vartheta),r^{+}(\vartheta)]$
results from a symmetric treatment of the energy production terms,
as presented in \citep{duran2026}. By allowing independent coefficients
in the corresponding Young inequalities and pushing the estimates
to their limiting values (which implies in particular the limit $\vartheta\to0$,
that is, $\Delta t/\Delta x\to0$), one obtains the necessary conditions
\[
\gamma\geq\frac{1}{2},\qquad\alpha\geq1
\]
for the nonlinear discrete energy argument to close. This refinement
only provides necessary conditions, and is not stated explicitly in
\citep{duran2026}, but follows from the same discrete energy estimates.
These bounds concern the closure of this refined energy argument;
they are not asserted as necessary conditions for stability of the
scheme itself.
\end{rem}

\paragraph{Relation with the linear calibration}

The nonlinear estimate and the linear analysis serve complementary
purposes. The former provides sufficient conditions for the fully
upwind, cellwise prescription. The latter selects a less dissipative
directional calibration around uniform states, with $\gamma_{\ell}=\alpha_{\ell}=\Gamma(Fr)$,
where $\Gamma$ is defined by \eqref{eq:gamma_law}. The sufficient
choice $\gamma=\alpha=2$ corresponds to $\gamma_{\ell}=8$ and $\alpha_{\ell}=4$;
the comparison therefore involves different prescriptions as well
as different coefficient bounds. Both prescriptions share the global
energy balance derived above, but the sufficient nonlinear estimate
does not certify the linear calibration for arbitrary nonlinear states.

\subsection{Acoustic subsystem}\label{subsec:splitting_acoustic}

We now turn to the discretization of the acoustic subsystem. The discrete
approach can be seen as a two-dimensional extension of what has recently
been proposed in \citep{duran2024}. We first propose the numerical
scheme and then discuss the main results leading to its stability.

\subsubsection{Numerical scheme}

Denoting $\Delta\tau$ the time step associated with the acoustic
part of the splitting, the numerical scheme we consider for \eqref{eq:splitting_acoustic}
is the following: 
\begin{subequations}
\label{eq:scheme_Acoustic} 
\begin{equation}
{\displaystyle h^{n+1}_{K}=h^{n}_{K}}\,,\label{eq:acoustic_mass}
\end{equation}

\begin{equation}
h^{n+1}_{K}\boldsymbol{u}^{n+1}_{K}=h^{n}_{K}\boldsymbol{u}^{n}_{K}-\frac{\Delta\tau}{m_{K}}\sum_{e\in\partial K}\left(\theta_{u}\overline{h^{n+1}p^{n+1}}_{e}+(1-\theta_{u})\overline{h^{n}p^{n}}_{e}\right)\boldsymbol{n}_{e,K}m_{e},\label{eq:acoustic_hu}
\end{equation}

\begin{equation}
{\displaystyle h^{n+1}_{K}w^{n+1}_{K}=h^{n}_{K}w^{n}_{K}+\frac{3\Delta\tau}{2}\left(\theta_{w}p^{n+1}_{K}+(1-\theta_{w})p^{n}_{K}\right)}\,,\label{eq:acoustic_hw}
\end{equation}

\begin{equation}
\begin{array}{l}
{\displaystyle h^{n+1}_{K}p^{n+1}_{K}=h^{n}_{K}p^{n}_{K}-a^{2}\Delta\tau\left(2\left(\theta_{p}w^{n+1}_{K}+(1-\theta_{p})w^{n}_{K}\right)+{\displaystyle \frac{h^{n}_{K}}{m_{K}}\sum_{e\in\partial K}\boldsymbol{u}^{*}_{e}\cdot\boldsymbol{n}_{e,K}m_{e}}\right)}\end{array}\,.\label{eq:acoustic_hp}
\end{equation}

Above, we set 
\begin{equation}
\boldsymbol{u}^{*}_{e}=\overline{\boldsymbol{u}^{n}}_{e}-\boldsymbol{\Gamma}_{e},\qquad\boldsymbol{\Gamma}_{e}=\beta\Delta\tau\,\overline{\left(\frac{r}{h^{n}}\right)}_{e}\,\boldsymbol{\delta}_{e}(h^{n}p^{n}),\label{eq:ustar}
\end{equation}
\end{subequations}
 where $\beta>0$ is a stabilization constant. As with $\boldsymbol{\Pi}_{e}$
and $\Lambda_{e}$ in the gravity step, $\boldsymbol{\Gamma}_{e}$
is a diffusion term: it regularizes the interface velocity through
the jump of the pressure $hp$, vanishes at the lake-at-rest equilibrium,
where $p=0$, and is responsible for the discrete energy dissipation
of the acoustic step established below. The weights $\theta_{u},\theta_{w},\theta_{p}\in[0,1]$
introduced in \eqref{eq:acoustic_hu}--\eqref{eq:acoustic_hp} set
the implicitness of the pressure gradient in the momentum equation,
of the “$p$” source term in the “$w$” equation \eqref{eq:acoustic_hw}
and of the “$w$” term in the “$p$” equation \eqref{eq:acoustic_hp},
respectively. The choice used in practice is $\theta_{u}=\theta_{w}=1$
and $\theta_{p}=1/2$. The linear stability analysis below establishes
an exact energy-contractivity condition for this choice and a sufficient
condition for $\theta_{p}\in[1/2,1]$. The conditions ensuring the
stability of the scheme are discussed in the following sections.

\subsubsection{Linear stability in the energy norm: acoustic step}\label{subsec:linear_acoustic}

As in Section~\ref{subsec:linear_sw}, we study the linear stability
of the acoustic scheme \eqref{eq:scheme_Acoustic} in the energy norm
\eqref{eq:energy_lin}, about the same uniform state. The acoustic
step carries no transport of the background flow, which belongs to
the gravity step, so the analysis holds for any background velocity;
and since $h$ is unchanged during the step, its contribution cancels
from the energy increment. We use the Fourier variables and the symbols
$s_{\bullet}=i\sin\xi_{\bullet}$, $c_{\bullet}=1-\cos\xi_{\bullet}$
of Section~\ref{subsec:linear_sw}, and take $\theta_{u}=\theta_{w}=1$
and $\theta_{p}\in[1/2,1]$.

\paragraph{Amplification matrix}

On a square mesh $\Delta=\Delta x=\Delta y$, set

\begin{equation}
{\displaystyle \sigma_{a}=\frac{a\Delta\tau}{\Delta},\qquad\kappa=\frac{\sqrt{3}\,a\Delta\tau}{H},\qquad\delta_{p}=1+\theta_{p}\kappa^{2},\qquad\mathcal{B}_{\beta}=1-4\beta\sigma^{2}_{a}(c_{x}+c_{y}).}\label{eq:acoustic_param}
\end{equation}

The full acoustic amplification matrix is $\mathbf{A}_{a}=\mathrm{diag}(1,\mathbf{A}^{\boldsymbol{u}wp}_{a})$.
Its block on $(u,v,w,p)$, the height being decoupled, reads

\begin{equation}
{\displaystyle \mathbf{A}^{\boldsymbol{u}wp}_{a}=\frac{1}{\delta_{p}}\begin{pmatrix}{\displaystyle \delta_{p}+\sigma^{2}_{a}s^{2}_{x}} & {\displaystyle \sigma^{2}_{a}s_{x}s_{y}} & {\displaystyle \frac{2}{\sqrt{3}}\kappa\sigma_{a}s_{x}} & {\displaystyle -\frac{\sigma_{a}}{a}s_{x}\mathcal{B}_{\beta}}\\[4pt]
{\displaystyle \sigma^{2}_{a}s_{x}s_{y}} & {\displaystyle \delta_{p}+\sigma^{2}_{a}s^{2}_{y}} & {\displaystyle \frac{2}{\sqrt{3}}\kappa\sigma_{a}s_{y}} & {\displaystyle -\frac{\sigma_{a}}{a}s_{y}\mathcal{B}_{\beta}}\\[4pt]
{\displaystyle -\frac{\sqrt{3}}{2}\kappa\sigma_{a}s_{x}} & {\displaystyle -\frac{\sqrt{3}}{2}\kappa\sigma_{a}s_{y}} & {\displaystyle 1-(1-\theta_{p})\kappa^{2}} & {\displaystyle \frac{\sqrt{3}}{2a}\kappa\mathcal{B}_{\beta}}\\[4pt]
{\displaystyle -a\sigma_{a}s_{x}} & {\displaystyle -a\sigma_{a}s_{y}} & {\displaystyle -\frac{2a}{\sqrt{3}}\kappa} & {\displaystyle \mathcal{B}_{\beta}}
\end{pmatrix}.}\label{eq:A_acoustic}
\end{equation}

Let $\mathbf{N}=\mathrm{diag}(i,i,2/\sqrt{3},1/a)$ and $\boldsymbol{Y}^{n}=\mathbf{N}(\hat{u}^{n},\hat{v}^{n},\hat{w}^{n},\hat{p}^{n})^{T}$,
where the hats denote the Fourier amplitudes of Section~\ref{subsec:linear_sw}.
The acoustic part of $E_{\mathrm{lin}}$ is $(H/2)\lVert\boldsymbol{Y}^{n}\rVert^{2}$,
and $\boldsymbol{Y}^{n+1}=\tilde{\mathbf{A}}_{a}\boldsymbol{Y}^{n}$
with

\begin{equation}
{\displaystyle \tilde{\mathbf{A}}_{a}=\mathbf{N}\mathbf{A}^{\boldsymbol{u}wp}_{a}\mathbf{N}^{-1}=\frac{1}{\delta_{p}}\begin{pmatrix}{\displaystyle \delta_{p}\mathbf{I}_{3}-\boldsymbol{b}\boldsymbol{b}^{T}} & {\displaystyle \mathcal{B}_{\beta}\boldsymbol{b}}\\
{\displaystyle -\boldsymbol{b}^{T}} & {\displaystyle \mathcal{B}_{\beta}}
\end{pmatrix},\qquad\boldsymbol{b}=\bigl(\sigma_{a}\sin\xi_{x},\ \sigma_{a}\sin\xi_{y},\ \kappa\bigr)^{T}.}\label{eq:A_acoustic_sym}
\end{equation}

The energy does not increase over one acoustic sub-step, for every
initial amplitude, if and only if 
\[
{\displaystyle \mathbf{I}-\tilde{\mathbf{A}}^{*}_{a}\tilde{\mathbf{A}}_{a}\geq0.}
\]

\paragraph{Reduction to a $2\times2$ problem}

Since $\boldsymbol{b}$ is real, $\boldsymbol{b}^{T}=\boldsymbol{b}^{*}$,
and $\lVert\boldsymbol{b}\rVert\geq\kappa>0$ for a positive time
step. Consider the two-dimensional subspace

\[
{\displaystyle \mathcal{K}=\left\{ (\boldsymbol{z},0)^{T}\,:\ \boldsymbol{z}\in\mathbb{C}^{3},\ \boldsymbol{b}^{*}\boldsymbol{z}=0\right\} .}
\]

Every vector of $\mathcal{K}$ is a fixed point of the acoustic step,

\[
{\displaystyle \tilde{\mathbf{A}}_{a}\begin{pmatrix}{\displaystyle \boldsymbol{z}}\\
{\displaystyle 0}
\end{pmatrix}=\frac{1}{\delta_{p}}\begin{pmatrix}{\displaystyle \delta_{p}\boldsymbol{z}-\boldsymbol{b}(\boldsymbol{b}^{T}\boldsymbol{z})}\\
{\displaystyle -\boldsymbol{b}^{T}\boldsymbol{z}}
\end{pmatrix}=\begin{pmatrix}{\displaystyle \boldsymbol{z}}\\
{\displaystyle 0}
\end{pmatrix},}
\]
so that 1 is an eigenvalue of $\tilde{\mathbf{A}}_{a}$ of multiplicity
at least two, whose eigenspace contains $\mathcal{K}$.

The orthogonal complement of $\mathcal{K}$ in $\mathbb{C}^{4}$ is
$\mathcal{K}^{\perp}=\left\{ (z_{b}\boldsymbol{b},z_{p})^{T}:\ z_{b},z_{p}\in\mathbb{C}\right\} $,
with orthonormal basis 
\[
{\displaystyle \boldsymbol{f}_{1}=\left(\frac{\boldsymbol{b}}{\lVert\boldsymbol{b}\rVert},0\right)^{T},\qquad\boldsymbol{f}_{2}=(0,0,0,1)^{T}.}
\]
Applying $\tilde{\mathbf{A}}_{a}$ to these two vectors gives

\[
{\displaystyle \tilde{\mathbf{A}}_{a}\boldsymbol{f}_{1}=\Bigl(1-\frac{\lVert\boldsymbol{b}\rVert^{2}}{\delta_{p}}\Bigr)\boldsymbol{f}_{1}-\frac{\lVert\boldsymbol{b}\rVert}{\delta_{p}}\boldsymbol{f}_{2},\qquad\tilde{\mathbf{A}}_{a}\boldsymbol{f}_{2}=\frac{\lVert\boldsymbol{b}\rVert\,\mathcal{B}_{\beta}}{\delta_{p}}\boldsymbol{f}_{1}+\frac{\mathcal{B}_{\beta}}{\delta_{p}}\boldsymbol{f}_{2},}
\]

so that $\mathcal{K}^{\perp}$ is invariant as well. Taking first
an orthonormal basis of $\mathcal{K}$ and then $(\boldsymbol{f}_{1},\boldsymbol{f}_{2})$
as columns yields a unitary matrix $\mathbf{Q}$ such that

\begin{equation}
{\displaystyle \mathbf{Q}^{*}\tilde{\mathbf{A}}_{a}\mathbf{Q}=\begin{pmatrix}{\displaystyle \mathbf{I}_{2}} & {\displaystyle 0}\\
{\displaystyle 0} & {\displaystyle \mathbf{R}}
\end{pmatrix},\qquad\mathbf{R}=\frac{1}{\delta_{p}}\begin{pmatrix}{\displaystyle \delta_{p}-\lVert\boldsymbol{b}\rVert^{2}} & {\displaystyle \lVert\boldsymbol{b}\rVert\,\mathcal{B}_{\beta}}\\
{\displaystyle -\lVert\boldsymbol{b}\rVert} & {\displaystyle \mathcal{B}_{\beta}}
\end{pmatrix}.}\label{eq:R_acoustic}
\end{equation}

Since $\det(\mathbf{I}_{2}-\mathbf{R})=\lVert\boldsymbol{b}\rVert^{2}/\delta_{p}\neq0$,
1 is not an eigenvalue of $\mathbf{R}$: it is thus a double eigenvalue
of $\tilde{\mathbf{A}}_{a}$, with eigenspace exactly $\mathcal{K}$,
and the action of the acoustic step on the energy is entirely carried
by the $2\times2$ block $\mathbf{R}$. Writing $\mathbf{Q}^{*}\boldsymbol{Y}^{n}=(\boldsymbol{Z}^{n},\boldsymbol{V}^{n})^{T}$
with $\boldsymbol{Z}^{n},\boldsymbol{V}^{n}\in\mathbb{C}^{2}$, one
has $\boldsymbol{Z}^{n+1}=\boldsymbol{Z}^{n}$ and $\boldsymbol{V}^{n+1}=\mathbf{R}\boldsymbol{V}^{n}$,
and since $\mathbf{Q}$ is unitary,

\[
{\displaystyle E^{n+1}_{\mathrm{lin}}-E^{n}_{\mathrm{lin}}=\frac{H}{2}\,(\boldsymbol{V}^{n})^{*}\left(\mathbf{R}^{*}\mathbf{R}-\mathbf{I}_{2}\right)\boldsymbol{V}^{n}.}
\]

As $\boldsymbol{V}^{n}$ is arbitrary when the initial amplitude is,
the energy does not increase for every initial amplitude if and only
if $\mathbf{I}_{2}-\mathbf{R}^{*}\mathbf{R}\geq0$. Writing $q^{2}_{a}=\sigma^{2}_{a}(\sin^{2}\xi_{x}+\sin^{2}\xi_{y})$,
so that $\lVert\boldsymbol{b}\rVert^{2}=q^{2}_{a}+\kappa^{2}$, a
direct calculation gives

\[
{\displaystyle (\mathbf{I}-\mathbf{R}^{*}\mathbf{R})_{11}=\frac{\lVert\boldsymbol{b}\rVert^{2}}{\delta^{2}_{p}}\left(1-q^{2}_{a}+(2\theta_{p}-1)\kappa^{2}\right),\qquad\det(\mathbf{I}-\mathbf{R}^{*}\mathbf{R})=\frac{\lVert\boldsymbol{b}\rVert^{2}}{\delta^{2}_{p}}\left(1-q^{2}_{a}-\mathcal{B}^{2}_{\beta}+(2\theta_{p}-1)\kappa^{2}\right).}
\]

For a positive time step, energy contractivity is therefore equivalent
to

\begin{equation}
{\displaystyle q^{2}_{a}+\mathcal{B}^{2}_{\beta}\leq1+(2\theta_{p}-1)\kappa^{2}.}\label{eq:acoustic_cond}
\end{equation}

\begin{rem}[Long-wave limit]
In the long-wave limit $kH\to0$, obtained by letting $k\to0$ at
fixed $\Delta,H,\Delta\tau$, one has $q_{a}\to0$ and $\boldsymbol{b}\to(0,0,\kappa)^{T}$.
The acoustic step therefore tends to the identity on the horizontal
velocity: only the $(w,p)$ block remains active. For the gravity-wave
branch, the amplitudes of these auxiliary variables also vanish in
this limit, at fixed height amplitude (\ref{sec:linear_waves}). The
acoustic correction thus becomes negligible on long gravity waves.
\end{rem}

\paragraph{Stability condition}
\begin{prop}
\label{prop:acoustic_lin} Let $\theta_{u}=\theta_{w}=1$, $\beta>0$
and $\beta_{\ell}:=4\beta$. On a square mesh, the acoustic step is
energy-contractive at every Fourier frequency if
\begin{equation}
{\displaystyle \beta_{\ell}\geq1,\qquad\beta_{\ell}\,\sigma^{2}_{a}\leq\frac{1}{2},}\label{eq:acoustic_lin_cond}
\end{equation}
for every $\theta_{p}\in[1/2,1]$. For $\theta_{p}=1/2$, these conditions
are also necessary. 
\end{prop}
\begin{proof}
Since $(2\theta_{p}-1)\kappa^{2}\geq0$, it suffices to establish
\eqref{eq:acoustic_cond} for $\theta_{p}=1/2$. Let $c_{+}=c_{x}+c_{y}$
and $c_{2}=c^{2}_{x}+c^{2}_{y}$. Since $\sin^{2}\xi_{j}=2c_{j}-c^{2}_{j}$,
\[
{\displaystyle 1-q^{2}_{a}-\mathcal{B}^{2}_{\beta}=\sigma^{2}_{a}\left[2(\beta_{\ell}-1)c_{+}+c_{2}-\beta^{2}_{\ell}\sigma^{2}_{a}c^{2}_{+}\right].}
\]
Under \eqref{eq:acoustic_lin_cond}, $c_{2}\geq c^{2}_{+}/2$ and
$\beta^{2}_{\ell}\sigma^{2}_{a}\leq\beta_{\ell}/2$ bound the bracket
from below by 
\[
{\displaystyle (\beta_{\ell}-1)c_{+}(4-c_{+})/2\geq0,}
\]
 since $0\leq c_{+}\leq4$. For $\theta_{p}=1/2$, necessity follows
from two modes. For $\xi_{y}=0$ and $\xi_{x}\to0$, the bracket is
$2(\beta_{\ell}-1)c_{x}+O(c^{2}_{x})$, which requires $\beta_{\ell}\geq1$.
At the checkerboard mode $(\pi,\pi)$, $c_{+}=4$ and $c_{2}=8$,
and the bracket equals $8\beta_{\ell}(1-2\beta_{\ell}\sigma^{2}_{a})$,
which requires $\beta_{\ell}\sigma^{2}_{a}\leq1/2$.
\end{proof}

For $\theta_{p}=1/2$, the least admissible stabilization is thus
$\beta=1/4$, with $\sigma^{2}_{a}\leq1/2$. The structure is the
same as for the gravity step: the long-wave limit fixes the constant,
and the checkerboard mode the time step. For $\theta_{p}>1/2$ and
$\beta_{\ell}\geq1$, the same argument, with $c_{2}\geq c^{2}_{+}/2$
attained for $c_{x}=c_{y}$, shows that contraction at every frequency
is equivalent to the relaxed checkerboard condition $\beta_{\ell}\sigma^{2}_{a}\leq1/2+(6\theta_{p}-3)(\Delta/H)^{2}/(16\beta_{\ell})$,
which reduces to the upper bound in \eqref{eq:acoustic_lin_cond}
in the fine-mesh limit; \eqref{eq:acoustic_lin_cond} is then sufficient
but not optimal. We retain $\theta_{p}=1/2$, for which \eqref{eq:acoustic_lin_cond}
is exact. The Courant number $\sigma_{a}$ refers to one acoustic
sub-step; its relation with the implementation parameter $C_{a}$
is given in Remark~\ref{rem:courant_conventions}.

\paragraph{Splitting and dissipation maps}

For a gravity time step $\Delta t$, the acoustic step is subdivided
into $n_{\tau}$ equal sub-steps of length $\Delta\tau=\Delta t/n_{\tau}$.
Given an admissible acoustic Courant number $\sigma_{a,\max}>0$,
we choose
\[
{\displaystyle n_{\tau}=\max\left(1,\left\lceil \frac{a\Delta t}{\sigma_{a,\max}\Delta}\right\rceil \right),}
\]
so that $\sigma_{a}=a\Delta\tau/\Delta\leq\sigma_{a,\max}$. The full
splitting then has amplification matrix $\mathbf{A}^{n_{\tau}}_{a}\mathbf{A}_{g}$.
Since both steps are analysed in the same energy norm, it is energy-contractive
as soon as each factor is; the acoustic condition alone says nothing
about the gravity step, and conversely.

\begin{figure}[htbp]
\begin{centering}
\includegraphics[width=0.5\textwidth]{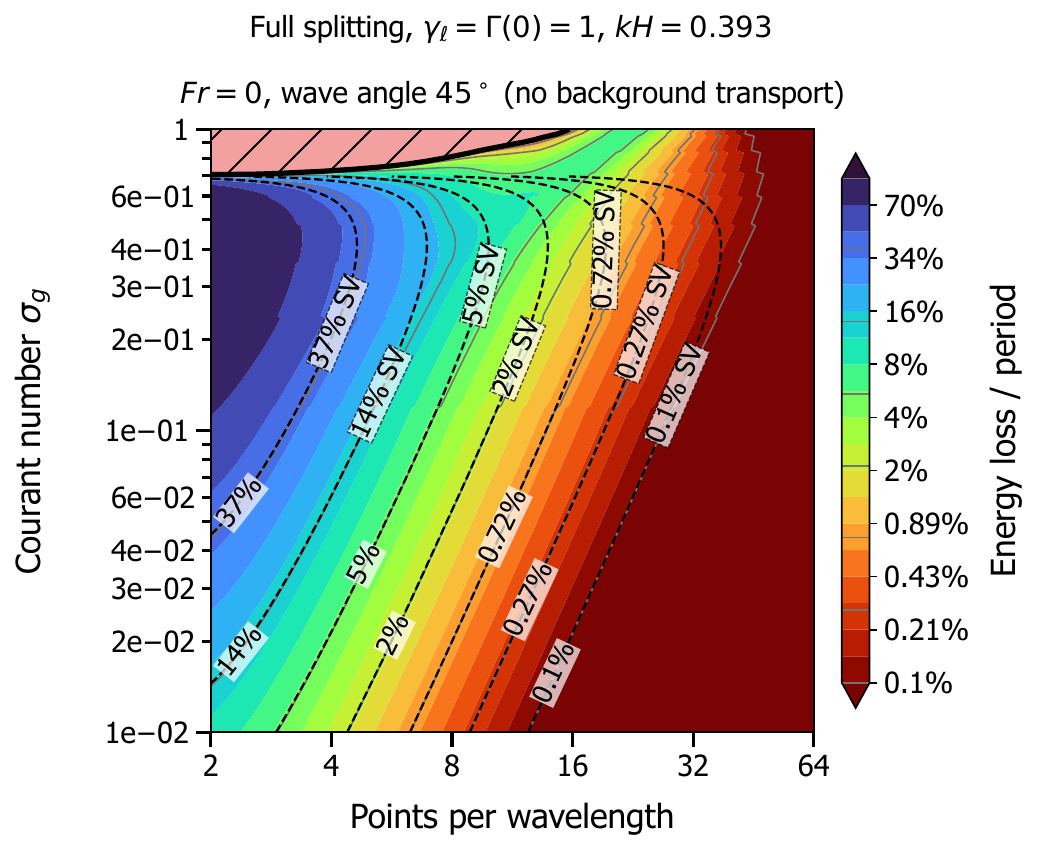}\includegraphics[width=0.5\textwidth]{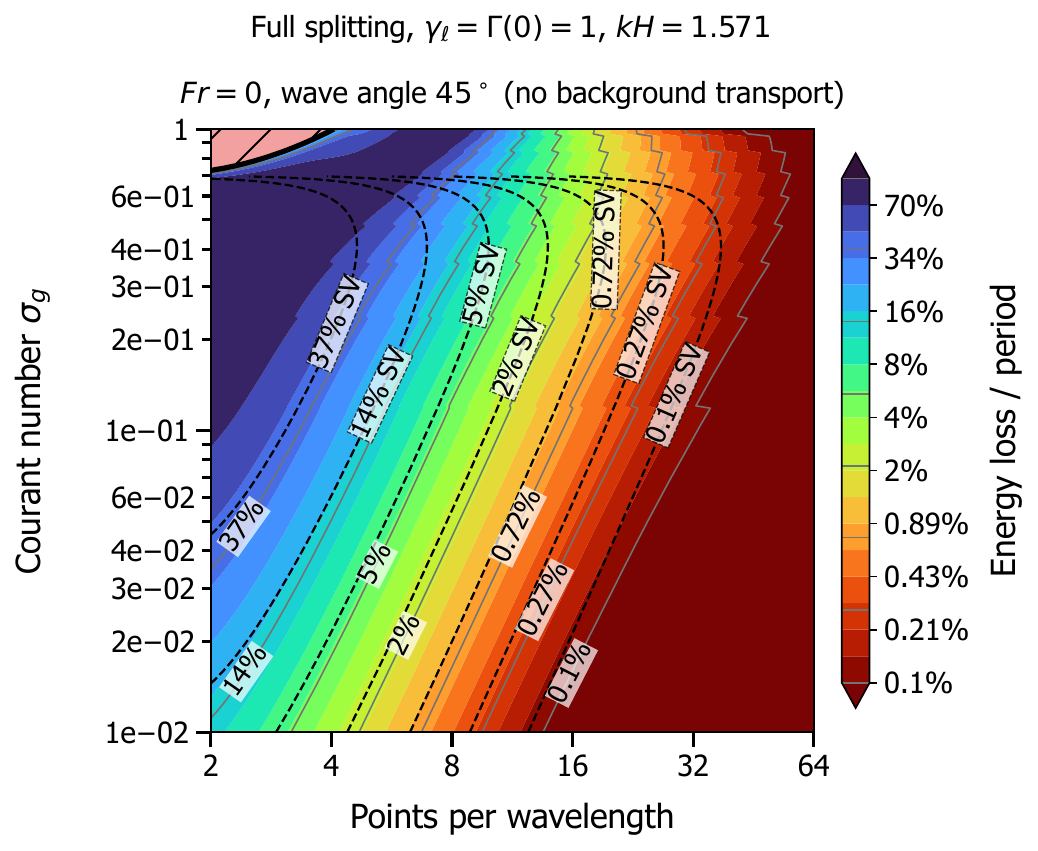}
\par\end{centering}
\caption{Energy loss per intrinsic period at $Fr=0$ for the full splitting
(colours and solid contours) and Saint-Venant alone (dashed contours,
labelled SV), as a function of $\sigma_{g}$ and $L/\Delta$. The
wave angle is $45^{\circ}$, $H=1\,\mathrm{m}$, $M=0.2$, $\gamma_{\ell}=\alpha_{\ell}=1$,
$\beta=1/4$, $\theta_{p}=1/2$, and the acoustic sub-step satisfies
$\sigma_{a}\protect\leq0.6$. Each model uses its own wave eigenmode
and period. Left: $L=16\,\mathrm{m}$ ($kH\approx0.39$); right: $L=4\,\mathrm{m}$
($kH\approx1.57$).}\label{fig:dissip-acoustic-Fr-0}
\end{figure}

Figure \ref{fig:dissip-acoustic-Fr-0} compares the full splitting
with the Saint-Venant reference at rest for two values of $kH$. The
coloured maps and solid contours represent the full splitting; dashed
contours represent Saint-Venant alone. Each model is initialized with
its own wave eigenmode, and the energy loss is measured per intrinsic
period of that model. The close agreement with the Saint-Venant reference
at small $kH$ indicates that attenuation is primarily governed by
the gravity step, while the acoustic splitting introduces only a small
additional effect. This is consistent with its vanishing action on
long gravity waves. A more pronounced departure appears at larger
$kH$, where the shallow-water approximation itself becomes less accurate.
Small breaks in the contours arise from changes in the integer number
of acoustic sub-steps.

\subsubsection{Nonlinear stability: acoustic step}

We finally turn to the derivation of fully discrete energy estimates
in the nonlinear case. During the acoustic step the heights are frozen,
so that the sub-problem is linear in $(\boldsymbol{u},w,p)$, with
spatially variable coefficients. The stability of the scheme relies
on an exact local energy balance and on the sufficient conditions
of Proposition~\ref{prop:acoustic_estim2}, both established in \ref{sec:dissipation}
(Sections~\ref{subsec:proof_acoustic_estim1} and \ref{subsec:proof_acoustic_estim2}).

For $\theta_{u}=\theta_{w}=1$ and $\theta_{p}\in[0,1]$, a direct
computation shows that the scheme \eqref{eq:acoustic_mass}, \eqref{eq:acoustic_hu},
\eqref{eq:acoustic_hw}, \eqref{eq:acoustic_hp} satisfies the exact
local energy balance
\begin{equation}
{\displaystyle (hE)^{n+1}_{K}-}{\displaystyle (hE)^{n}_{K}+\dfrac{\Delta\tau}{m_{K}}\sum_{e\in\partial K}\widetilde{(hp\boldsymbol{u})}_{e}\cdot\boldsymbol{n}_{e,K}m_{e}=\mathfrak{R}_{K}}\,,\label{eq:acoustic_estim1}
\end{equation}
\begin{equation}
\text{with}\quad{\displaystyle \mathfrak{R}_{K}=\frac{1}{2}h^{n}_{K}\left\Vert \boldsymbol{u}^{n+1}_{K}-\boldsymbol{u}^{n}_{K}\right\Vert ^{2}-\mathfrak{T}_{K}(\theta_{p})-\frac{h^{n}_{K}}{2a^{2}}\left(p^{n+1}_{K}-p^{n}_{K}\right)^{2}+\dfrac{\Delta\tau}{m_{K}}h^{n}_{K}p^{n+1}_{K}\sum_{e\in\partial K}\boldsymbol{\Gamma}_{e}\cdot\boldsymbol{n}_{e,K}m_{e}}\,,
\end{equation}
and $\widetilde{(hp\boldsymbol{u})}_{e}=\overline{(hp)}^{n+1}_{e}\,\overline{\boldsymbol{u}}^{n}_{e}-\boldsymbol{\delta}_{e}(hp^{n+1})(\delta_{e}\boldsymbol{u})$.
The $w$--$p$ exchange term $\mathfrak{T}_{K}(\theta_{p})$, a dissipation
for $\theta_{p}\geq1/2$, is given by 
\begin{equation}
\mathfrak{T}_{K}(\theta_{p})=\frac{2}{3}h^{n}_{K}\left(w^{n+1}_{K}-w^{n}_{K}\right)^{2}-2\Delta\tau(1-\theta_{p})p^{n+1}_{K}\left(w^{n+1}_{K}-w^{n}_{K}\right)\,.\label{eq:W-upwind-residual}
\end{equation}

Note that $\mathfrak{R}_{K}$ carries the positive (a priori anti-dissipative)
contribution $h^{n}_{K}\Vert\boldsymbol{u}^{n+1}_{K}-\boldsymbol{u}^{n}_{K}\Vert^{2}/2$
from the momentum increment; controlling it under the acoustic CFL
condition and the bounds on $\beta$ is precisely the object of the
next proposition.

\begin{prop}
\label{prop:acoustic_estim2}Assume $\theta_{u}=\theta_{w}=1$, $\theta_{p}\in[1/2,1]$
and $h_{K}>0$, with constant coefficients $a>0$ and $\beta>0$,
on a conforming polygonal mesh with straight edges. With $\mu_{K}=r_{K}/h_{K}$
and $\overline{\mu}_{e}=(\mu_{K}+\mu_{K_{e}})/2$, define for every
cell the symmetric matrix and the scalar

\begin{equation}
\mathsf{M}_{K}=\frac{1}{m_{K}h_{K}}\sum_{e\in\partial K}\frac{m_{e}}{\overline{\mu}_{e}}\,\boldsymbol{n}_{e,K}\boldsymbol{n}^{T}_{e,K},\qquad\mathsf{T}_{K}=\frac{h_{K}}{m_{K}}\sum_{e\in\partial K}m_{e}\,\overline{\mu}_{e},\label{eq:acoustic_geometry}
\end{equation}

and assume that, for every cell $K$,

\begin{equation}
\lambda_{\max}(\mathsf{M}_{K})\leq2\beta,\qquad\beta a^{2}\Delta\tau^{2}\,\mathsf{T}_{K}\leq2,\label{eq:CFL_acoustic}
\end{equation}

where $\lambda_{\max}$ denotes the largest eigenvalue. Then

\begin{equation}
(hE)^{n+1}_{K}\leq(hE)^{n}_{K}-\frac{\Delta\tau}{m_{K}}\sum_{e\in\partial K}\mathcal{G}_{a}(W_{K},W_{K_{e}},\boldsymbol{n}_{e,K})\,m_{e},\label{eq:acoustic_estim2}
\end{equation}

where $\mathcal{G}_{a}$ is a conservative energy flux, consistent
with $hp\boldsymbol{u}\cdot\boldsymbol{n}$, made of the leading part
$\widetilde{(hp\boldsymbol{u})}_{e}$ of \eqref{eq:acoustic_estim1}
and of two corrections of order $\beta\Delta\tau$, given in \eqref{eq:Ga};
$W_{K}$ gathers the states of the cell $K$ at both time levels.
\end{prop}

The first condition only involves the geometry of $K$ and of its
neighbours and the local depth ratios, and bounds $\beta$ from below;
the second one is the acoustic CFL condition. No small-amplitude assumption
is made. Under periodic boundary conditions, the fluxes cancel in
the sum over the cells and the total energy does not increase.

\begin{rem}
\label{rem:acoustic_meshes}On a uniform square mesh of side $\Delta$
at constant depth, $\mathsf{M}_{K}=I/2$ and $\mathsf{T}_{K}=16/\Delta^{2}$.
With $\beta_{\ell}=4\beta$ and $\sigma_{a}=a\Delta\tau/\Delta$,
the conditions \eqref{eq:CFL_acoustic} read $\beta_{\ell}\geq1$
and $\beta_{\ell}\sigma^{2}_{a}\leq1/2$, exactly those of the linear
analysis for $\theta_{p}=1/2$ (Section \ref{subsec:linear_acoustic}):
the nonlinear estimate is sharp in this configuration. More generally,
since $\overline{\mu}_{e}\geq\mu_{K}/2=r_{K}/(2h_{K})$, one has $\mathsf{M}_{K}\preceq(2/m_{\partial K})\sum_{e}m_{e}\boldsymbol{n}_{e,K}\boldsymbol{n}^{T}_{e,K}$
for any positive depths, so that $\lambda_{\max}(\mathsf{M}_{K})\leq2$
on any mesh; at constant depth, $\lambda_{\max}(\mathsf{M}_{K})=1/2$
for meshes of congruent regular polygons. For the interior cells tested
at constant depth on the Frontal-Delaunay meshes used in Section \ref{sec:validation},
$\lambda_{\max}(\mathsf{M}_{K})$ was observed not to exceed $0.85$,
so that the value $\beta_{\ell}=2$ used in the computations satisfies
the first condition for these cells at constant depth. This geometric
check alone does not certify the nonlinear computations: both local
conditions must be checked with the actual depths, and the boundary
treatment must be considered separately.
\end{rem}

\section{Implementation issues}\label{sec:implementation}

\subsection{Time stepping and acoustic sub-iterations}\label{subsec:time_stepping}

The two subsystems are advanced with two different time steps, reflecting
the strict separation between the gravity wave speed $\sqrt{gH}$
and the acoustic speed $a=\sqrt{gH}/M$, with $M\ll1$ the Mach number.
The \emph{effective} time step of the global scheme is that of the
shallow-water subsystem, set by its CFL-like condition, defined on
the unstructured mesh using the prescribed Courant parameter $C_{g}$
by 
\begin{equation}
\Delta t\,\leq\,C_{g}\,\left[\,\max_{K\in\mathbb{T}}\frac{m_{\partial K}\left(\left\Vert \boldsymbol{u}^{n}_{K}\right\Vert +\sqrt{gh^{n}_{K}}\right)}{2m_{K}}\,\right]^{-1}\,,\label{eq:cfl_sw_impl}
\end{equation}
with $C_{g}\leq1$. This is the practical wave-speed restriction used
in the computations. It does not by itself imply the sufficient nonlinear
conditions \eqref{eq:CFL_hyperbolic} and \eqref{eq:cpr-transport-explicit},
which must be checked together with the associated viscosity and upwinding
assumptions to invoke that result.

The acoustic subsystem carries the faster acoustic waves, of speed
$a$; we therefore \emph{sub-iterate} the acoustic step within each
global step. The interval $\Delta t$ is split into $n_{\tau}$ equal
acoustic sub-steps $\Delta\tau=\Delta t/n_{\tau}$, where $n_{\tau}$
is the smallest integer ensuring the practical acoustic CFL condition
with prescribed parameter $C_{a}$, 
\begin{equation}
n_{\tau}=\left\lceil \,\frac{\Delta t}{\Delta\tau_{\max}}\,\right\rceil \,,\qquad\Delta\tau_{\max}\,=\,C_{a}\,\left(a\:\max_{K\in\mathbb{T}}\frac{m_{\partial K}}{2m_{K}}\right)^{-1}\,,\label{eq:nsub}
\end{equation}
so that $a\,\Delta\tau\,m_{\partial K}/(2m_{K})\leq C_{a}$ on every
cell $K$.
\begin{rem}
\label{rem:courant_conventions}The prescribed Courant parameters
$C_{g}$ and $C_{a}$ use the geometric factor $m_{\partial K}/(2m_{K})$.
On a square mesh of side $\Delta$, this factor equals $2/\Delta$.
At a uniform state, the practical restrictions therefore read $2\sigma_{g}\leq C_{g}$
and $2\sigma_{a}\leq C_{a}$, where $\sigma_{g}=(\|\boldsymbol{U}\|+\sqrt{gH})\Delta t/\Delta$
and $\sigma_{a}=a\Delta\tau/\Delta$ are the Courant numbers used
in the linear analyses. The acoustic inequality may be strict because
the sub-iteration count is rounded up. This factor of two is a convention,
not a change in the time steps. On a uniform square mesh at constant
depth, the sufficient acoustic conditions of Proposition \ref{prop:acoustic_estim2}
reduce, through the practical time-step rule, to $\beta_{\ell}\ge1$
and $\beta_{\ell}C^{2}_{a}\le2$. In particular, $C_{a}\le1$ is sufficient
for $\beta_{\ell}=2$.
\end{rem}
Note that even with the implicit weights $\theta_{u}=\theta_{w}=1$
and $\theta_{p}=1/2$, the discrete equations for the auxiliary variables
$w$ \eqref{eq:acoustic_hw} and $p$ \eqref{eq:acoustic_hp} can
actually be solved in a fully explicit way. After basic calculations,
we indeed obtain 
\begin{equation}
p^{n+1}_{K}\left(1+\dfrac{3a^{2}\Delta\tau^{2}}{2(h^{n}_{K})^{2}}\right)=p^{n}_{K}-a^{2}\Delta\tau\left(\dfrac{2}{h^{n}_{K}}w^{n}_{K}+{\displaystyle \frac{1}{m_{K}}\sum_{e\in\partial K}\boldsymbol{u}^{*}_{e}\cdot\boldsymbol{n}_{e,K}m_{e}}\right)\,,\label{eq:acoustic_explicit_p}
\end{equation}
\begin{equation}
{\displaystyle w^{n+1}_{K}=w^{n}_{K}+\frac{3\Delta\tau}{2h^{n}_{K}}p^{n+1}_{K}},\label{eq:acoustic_explicit_w}
\end{equation}
so that the scheme $\eqref{eq:scheme_Acoustic}$ for the acoustic
part can be explicitly solved.

Each shallow-water subsystem step is thus followed by $n_{\tau}$
acoustic sub-steps, so that the cost of the acoustic part is set by
the chosen Mach number $M$; the linear analysis of the previous section
shows that these sub-iterations add little damping to the resolved
gravity wave.

The parameter values used in the computations are given in Section~\ref{sec:validation}.

\subsection{Wet/dry fronts}

The scheme of Section~\ref{subsec:splitting_sw} in its calibrated
low-dissipation form, is not monotone; near discontinuities it allows
oscillations that may generate negative water depths at wet/dry fronts.
Also, in such regions, dispersive effects are no longer physically
meaningful and need not be resolved by the acoustic step. We therefore
switch, in a continuous manner, to a more robust discretization for
both subsystems. The switch is driven by a mesh-based coefficient
$\chi_{K}\in[0,1]$ which equals one on the dry cells, where the water
depth falls below a cut-off $h<h_{\mathrm{cut}}$, equals zero in
the fully wet region, and varies smoothly in between through a $C^{2}$
Hermite profile for example, spreading from one cell to its neighbours,
layer by layer, over a fixed number of layers (typically five). At
an edge $e$ we use $\chi_{e}=\max(\chi_{K},\chi_{K_{e}})$.

For the shallow-water (gravity) step, the flux near the front is blended
towards a robust Rusanov flux with hydrostatic reconstruction \citep{audusse2004}.
Denoting by $\eta=h+b$ the free-surface elevation, the water depths
are first reconstructed: 
\[
\eta^{*}_{e}=\max\!\left(\eta_{K}-h_{K},\,\eta_{K_{e}}-h_{K_{e}}\right)\,,\quad h^{\star}_{K}=\max\!\left(0,\eta_{K}-\eta^{*}_{e}\right)\,,\quad h^{\star}_{K_{e}}=\max\!\left(0,\eta_{K_{e}}-\eta^{*}_{e}\right)\,,
\]
which ensures non-negative reconstructed water depths. The corresponding
Rusanov flux $\mathcal{F}^{\mathrm{Rus}}_{e}$ (mass and momentum
components) uses the local wave speed $s_{e}=\max\big(\vert\boldsymbol{u}_{K}\cdot\boldsymbol{n}\vert+\sqrt{gh^{\star}_{K}},\,\vert\boldsymbol{u}_{K_{e}}\cdot\boldsymbol{n}\vert+\sqrt{gh^{\star}_{K_{e}}}\big)$;
its hydrostatic pressure part is removed and carried by the potential
$\phi^{\chi}_{e}$ below, so that it is not counted twice. Both the
mass and momentum fluxes and the pressure term actually applied at
the edge are then taken as the convex combination 
\[
\mathcal{F}^{\chi}_{e}=\chi_{e}\,\mathcal{F}^{\mathrm{Rus}}_{e}+\left(1-\chi_{e}\right)\mathcal{F}^{*}_{e}\,,\qquad\phi^{\chi}_{e}=\chi_{e}\,g\left(\dfrac{\eta_{K}+\eta_{K_{e}}}{2}\right)+\left(1-\chi_{e}\right)\phi^{*}_{e}\,,
\]
so that the scheme uses the entropy-stable flux and pressure $(\mathcal{F}^{*}_{e},\phi^{*}_{e})$
in the wet region $(\chi_{e}=0)$ and the robust hydrostatic Rusanov
ones at the front $(\chi_{e}=1)$, with a smooth transition.

In the transition region, the acoustic speed is reduced according
to $a_{K}=a(1-\chi_{K})$. This reduction weakens the pressure response
to the velocity divergence and vertical motion, but does not by itself
switch off all acoustic corrections. On the dry cells, where $h<h_{\mathrm{cut}}$,
the auxiliary variables $p_{K}$ and $w_{K}$ are set to zero. This
near-dry treatment is a practical robustness measure. The robustness
of this treatment is illustrated by the realistic Saint-Malo configuration
of Section~\ref{subsec:saint-malo}, whose macro-tidal setting involves
extensive wet/dry areas.

In this blending layer the scheme prioritizes positivity and robustness
over the discrete energy estimate: the convex combination of fluxes
and the damped acoustic speed are not covered by the energy analysis
of the previous sections, which applies to the fully wet region $\chi_{e}=0$.
In practice this is confined to a few cells around the moving shoreline
and did not affect the stability of any of the computations below.

\subsection{Boundary conditions}\label{subsec:bc}

Boundary conditions are imposed weakly through ghost cells, consistently
for the gravity and the acoustic step. For each boundary edge $e$
separating an interior cell $K$ from an exterior ghost cell $K_{e}$,
the interior state is rotated into the local edge frame $(\boldsymbol{n},\boldsymbol{t})$,
a ghost state $(h,u_{n},u_{t},w,p)$ is prescribed according to the
boundary type, and the state is rotated back before the edge fluxes
and the acoustic interface velocity are evaluated. Both steps share
the same ghost values, so the splitting structure is preserved up
to the boundary.

Solid walls are imposed by reflecting the normal velocity, $u^{K_{e}}_{n}=-u^{K}_{n}$,
and copying the remaining variables. The interface normal velocity
then vanishes, $\bigl((\boldsymbol{u}_{K}+\boldsymbol{u}_{K_{e}})/2\bigr)\cdot\boldsymbol{n}=0$,
so that the mass flux $\mathcal{F}^{*}_{e}$ and the acoustic transport
velocity $\boldsymbol{u}^{*}_{e}\cdot\boldsymbol{n}$ vanish, and
with them the physical energy fluxes $\left(E+gh/2\right)h\boldsymbol{u}\cdot\boldsymbol{n}$
and $hp\boldsymbol{u}\cdot\boldsymbol{n}$. Thus, no physical energy
flux crosses the wall, although the gravity-step stabilization may
introduce additional numerical dissipation at the boundary.

Open (or forcing) boundaries, where a time-dependent water level $\eta_{b}(t)$
(tide, wave, or a combination) or a height $h_{b}(t)$ is prescribed,
require a specific treatment to avoid artificial reflections of outgoing
waves. The prescribed surface fixes the ghost height, while the normal
velocity is reconstructed from the interior state so as to preserve
the outgoing Riemann invariant $u_{n}+2\sqrt{gh}=\mathrm{cst}$ along
the outgoing characteristic: 
\begin{equation}
u^{K_{e}}_{n}=u^{K}_{n}+2\sqrt{g}\left(\sqrt{h_{K}}-\sqrt{h_{K_{e}}}\right)\,.\label{eq:riemann_bc}
\end{equation}
The auxiliary variables $w$ and $p$ are simply copied from the adjacent
interior cell, which amounts to a homogeneous Neumann (zero-gradient)
condition. This characteristic prescription provides a practical boundary
treatment by imposing the water level while retaining the outgoing
shallow-water Riemann invariant. It does not, however, guarantee complete
transparency: reflected waves may remain, even in the long-wave regime.
In the tests below, it is used for wave generation, while outgoing
waves are absorbed using a relaxation layer. Improving this boundary
treatment, particularly for dispersive waves, remains a direction
for future work, drawing on the boundary-condition analyses of \citet{lannes_general_2024}.

Where outgoing waves carry significant energy, a relaxation layer
is added along the open boundary. In a band of width $L_{r}$, the
solution is driven after each time step towards a prescribed sea surface
height $\eta_{b}$ at rest:

\begin{equation}
h+b\leftarrow\varrho\,(h+b)+(1-\varrho)\,\eta_{b}\,,\qquad(u,v,w,p)\leftarrow\varrho\,(u,v,w,p)\,,\label{eq:relaxation}
\end{equation}

with $\varrho(s)=1-\bigl(\exp\bigl((1-s)^{6}\bigr)-1\bigr)/(e-1)$,
$s=d/L_{r}\in[0,1]$, where $d$ is the distance to the boundary.
The weight $\varrho$ vanishes at the boundary and equals one at the
inner edge of the layer. Acting on the free surface $h+b$, the relaxation
preserves the lake at rest, and the height is clipped at zero to preserve
positivity.

\subsection{Computational cost}

The Mach number $M=\sqrt{gH}/a$ sets the artificial acoustic speed
$a$ and is the main tuning parameter of the method. The approximation
of the target Serre--Green--Naghdi dynamics improves as $M\to0$
(that is $a\to\infty$), the relaxation error on the dispersion relation
scaling like $M^{2}/2$ (see \ref{sec:linear_waves}), while the number
of acoustic sub-iterations grows as $n_{\tau}\sim1/M$ (see \eqref{eq:nsub}).
In practice we recommend $M\in[0.2,0.4]$, based on the dispersion
comparison in \ref{sec:linear_waves}, Fig.~\ref{fig:dispersion-M}:
this provides a practical compromise between relaxation error and
computational overhead.

\begin{figure}[htbp]
\centering{}\includegraphics[width=0.7\textwidth]{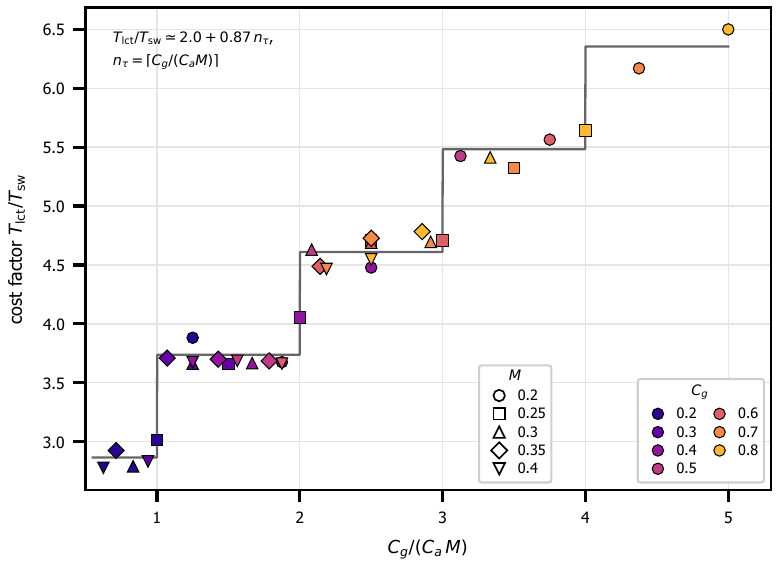}
\caption{Measured cost factor $T_{\mathrm{lct}}/T_{\mathrm{sw}}$ (full Leucothea
solver over the pure shallow-water solver, same time steps) versus
$C_{g}/(C_{a}M)$, whose ceiling is the number of acoustic sub-steps
per time step $n_{\tau}$, for $C_{g}\in[0.2,0.8]$ (colours), $M\in[0.2,0.4]$
(markers) and $C_{a}=0.8$, on the linear-wave test case (about $1.5\times10^{5}$
cells, serial runs, minimum over three runs). The staircase is the
least-squares fit $T_{\mathrm{lct}}/T_{\mathrm{sw}}\simeq2.0+0.87\,\lceil C_{g}/(C_{a}M)\rceil$.}\label{fig:cost}
\end{figure}

The overhead compared to the original shallow-water solver comes from
the two additional transported scalars $w$ and $p$ in the gravity
step and from the $n_{\tau}$ acoustic sub-steps per global step.
It was measured on the linear-wave test case of Section~\ref{subsec:linear_regime},
on an unstructured mesh of about $1.5\times10^{5}$ cells, as the
ratio $T_{\mathrm{lct}}/T_{\mathrm{sw}}$ of the wall-clock times
of the full solver and of the pure shallow-water solver run at the
same gravity Courant number, hence with the same time steps; each
time is the minimum over three serial runs without output. Seven Courant
numbers $C_{g}\in[0.2,0.8]$ and five Mach numbers $M\in[0.2,0.4]$
were tested, with $C_{a}=0.8$. As shown in Fig.~\ref{fig:cost},
the cost factor depends on $(C_{g},M)$ only through $n_{\tau}=\lceil C_{g}/(C_{a}M)\rceil$,
and is accurately fitted by $T_{\mathrm{lct}}/T_{\mathrm{sw}}\simeq2.0+0.87\,n_{\tau}$.
The gravity step with its two additional scalars thus costs about
twice a shallow-water step, and each acoustic sub-step slightly less
than one: although it involves no Riemann solver, it performs an edge
loop of comparable cost. In the recommended range $M\in[0.2,0.4]$,
$n_{\tau}$ varies from $1$ to $5$, and the full solver costs three
to six times the shallow-water one, the lower values corresponding
to small Courant numbers or large Mach numbers.

This overhead remains moderate for a dispersive model. In the classical
splitting strategies for the SGN equations, the dispersive correction
requires the solution of an elliptic problem at every time step (see
e.g. \citep{bonneton2011}), which is the most expensive part of the
algorithm. By replacing this elliptic solver with a fully explicit,
sub-iterated hyperbolic relaxation, the present approach keeps the
cost at a small multiple of the shallow-water solver while retaining
an explicit, local and easily parallelizable algorithm. Moreover,
each building block of the method is as cheap as a first-order finite-volume
scheme: the two subsystems, and all the variants derived from them,
use a single flux evaluation per edge, with no MUSCL-type reconstruction
(which, in our implementation, multiplies the cost by about four with
second-order RK time stepping) and no discontinuous Galerkin quadrature
with several points per cell. No global linear system is ever assembled
or solved either, so the whole solver reduces to explicit, local edge
updates that scale linearly with the number of cells and parallelize
naturally on unstructured meshes.

\section{Test cases}\label{sec:validation}

In this section, we present a series of test cases to illustrate the
performance of the present numerical formulation for solving the Leucothea
model for dispersive waves. A first case is the wave propagation in
an approximately linear regime considering a non-trivial bathymetry.
A second is the simulation of a periodic travelling wave in a nonlinear
regime using the numerically computed solution of the Leucothea model
\eqref{eq:model}. A third one reproduces the well-known experiments
of Vincent \& Briggs \citep{vincent1989}, a classical benchmark for
Boussinesq-like numerical solvers. Finally, the potential of the solver
is demonstrated in a realistic coastal geometry: the bay of Saint-Malo
on the macro-tidal coast of Brittany, built on real SHOM bathymetric
and coastline data, with a directional TMA (Texel--Marsen--Arsloe)
wave spectrum imposed at the open boundary. For this last case we
focus on the convergence of the wave-height statistics with the mesh
resolution. Although limited in number, these four cases were chosen
to span the regimes targeted by the solver: quasi-linear dispersive
propagation over variable bathymetry, fully nonlinear permanent waves,
experimentally measured refraction--focusing, and a realistic operational
configuration combining real bathymetry, a complex coastline, spectral
forcing and extensive wet/dry areas.

Unless otherwise stated, all computations use the calibrated CPR scheme
for the gravity step, with $\theta_{e}$ from \eqref{eq:theta_mesh}
and $\gamma_{\ell}=\alpha_{\ell}=\Gamma(Fr_{e})$ from \eqref{eq:gamma_law}.
Its viscosity constants are fixed by the calibration, the remaining
parameters ($C_{g}$, the edge averages, the relaxation layer) being
algorithmic choices. For the acoustic step, the computations use the
weights $\theta_{u}=\theta_{w}=1$ and $\theta_{p}=1/2$, the stabilization
$\beta=1/2$ ($\beta_{\ell}=2$), and the Courant parameter $C_{a}=0.8$.
The Mach number $M$ and the gravity Courant parameter $C_{g}$ are
specified for each case. These settings satisfy the linear stability
conditions of Section~\ref{subsec:linear_acoustic}: $\beta_{\ell}\geq1$
and $\beta_{\ell}\sigma^{2}_{a}\leq0.32<1/2$, with $\sigma_{a}\leq C_{a}/2$
(see Remark \ref{rem:courant_conventions}). The value $\beta_{\ell}=2$
keeps a margin over the minimum $\beta_{\ell}=1$, at a small cost
in dissipation. By design, the gravity calibration is not covered
by the sufficient nonlinear estimate (see Section~\ref{subsec:nonlinear-sw}).
The acoustic settings satisfy the sufficient conditions of Proposition
\ref{prop:acoustic_estim2} on a uniform square mesh at constant depth.
The geometric checks reported in Remark \ref{rem:acoustic_meshes}
concern interior cells at constant depth; certification of the nonlinear
test cases additionally requires checking both local conditions with
the computed depths and accounting for the boundary treatment. No
growth of the total discrete energy was observed in the periodic wave-propagation
tests.

For comparison purposes, four numerical schemes were used for the
first and second test cases. The first three all use the proposed
splitting strategy, which separates the gravity part from the acoustic
part, the acoustic step being treated identically in the three; they
differ only in the discretization of the gravity part. The first two
use the CPR scheme \eqref{eq:scheme_sw}: in its fully upwind variant
($\theta_{e}=1$, $\gamma_{\ell}=1$) and in its calibrated variant
($\theta_{e}$ from \eqref{eq:theta_mesh}, $\gamma_{\ell}=\Gamma(Fr_{e})$).
The third one uses a classical Rusanov flux with MUSCL reconstructions
without any limitation (since the solutions considered are largely
smooth) and an RK2-Heun time stepping, yielding a second-order discretization
of the gravity part. The fourth scheme does not use the proposed splitting:
it is an implicit--explicit (IMEX) strategy, in the spirit of \citep{richard_extension_2021},
based on a strong-stability-preserving IMEX-SSP(3,3,2) L-stable time
stepping with a Rusanov flux and MUSCL reconstructions for the gravity
part and an implicit treatment of the auxiliary variables $w$ and
$p$. As in \citep{richard_extension_2021}, the velocity divergence
is evaluated explicitly within the implicit stage, which breaks the
formal second-order accuracy in time. The IMEX strategy advances the
whole system with a time-step estimate based on $\sqrt{gh+a^{2}}\simeq a$,
that is at the acoustic Courant number, rather than $\sqrt{gh}$ for
the splitting schemes. Comparing the first three schemes isolates
the influence of the gravity discretization within the same splitting
strategy, while the comparison with the fourth assesses the overall
performance of the splitting and IMEX strategies. The accuracy must
be weighed against the computational cost, which differs markedly
between the schemes.

The cost ratios below were measured on the approximately linear wave-propagation
test case, using the same mesh, final physical time and prescribed
Courant number ($C_{g}=0.8$ for the splitting schemes). The splitting
schemes take approximately the same number of global time steps. The
RK2-Heun/Rusanov/MUSCL scheme is about $2.7$ times more expensive
than the CPR scheme; its intrinsic cost is in fact about four times
higher, but this is partly masked by the acoustic sub-iterations shared
by all three splitting schemes. The IMEX strategy is about $17$ times
more expensive than the CPR scheme. Its time-step estimate based on
$\sqrt{gh+a^{2}}$ yields, at the same Courant number, about $\sqrt{1+1/M^{2}}\approx5$
times more time steps in this test ($M=0.2$), so that the measured
factor combines this larger number of steps with a higher cost per
step.

\subsection{Approximately linear regime}\label{subsec:linear_regime}

Here we present the ability of the solver to simulate wave propagation
in a nearly linear regime over a non-trivial bathymetry. We consider
a simple square box computational domain with periodic boundary conditions
in both directions. Since we are mainly interested in the simulation
of coastal flow, a side length $L=1\ \mathrm{km}$ and a reference
depth $H=20\ \mathrm{m}$ are chosen, which is a typical configuration
of interest. To initialize the water surface, we consider it as a
sum of linear waves for which the amplitudes, lengths and phase shifts
are chosen randomly from a certain range of values.

To construct a non-trivial bathymetry, a trigonometric basis function
with compact support is considered such that

\begin{equation}
\begin{gathered}b(\boldsymbol{x}_{K})=-20+b_{\max}\sum^{n_{b}}_{i=1}\mathscr{B}_{i}(\boldsymbol{x}_{K}),\\
\mathscr{B}_{i}(\boldsymbol{x})=r_{i}\sum_{\boldsymbol{p}\in\{-L,0,L\}^{2}}\psi\!\left(\frac{\|\boldsymbol{x}-\boldsymbol{p}-\boldsymbol{x}_{0,i}\|}{\varepsilon_{b}L}\right),\qquad\psi(s)=\begin{cases}
\dfrac{1+\cos(\pi s)}{2}, & 0\le s\le1,\\
0, & s>1.
\end{cases}
\end{gathered}
\end{equation}

\begin{figure}[htbp]
\centering{}\includegraphics[width=0.5\textwidth]{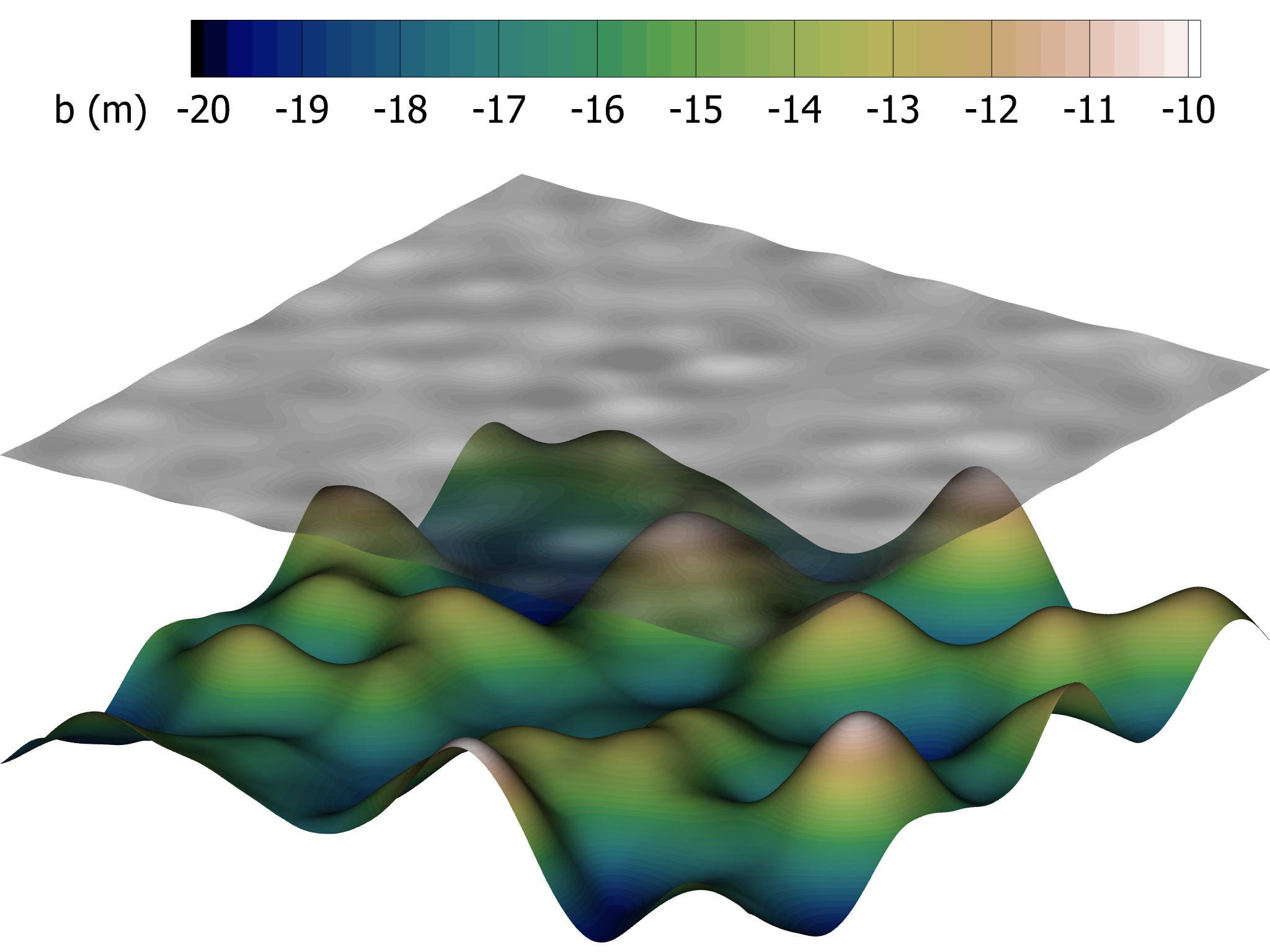}\caption{Test case of wave propagation in an approximately linear regime: built
non-trivial bathymetry and initial water surface.}\label{fig:test-linear-bathy}
\end{figure}

The nine translation vectors ensure the periodicity of the basis function
as long as $\varepsilon_{b}<0.5$. The resulting bottom is $C^{1}$,
which suffices here since the model involves only first derivatives
of $b$. The centres $\boldsymbol{x}_{0,i}$ are random coordinate
vectors in the square box, and $r_{i}\in\left[0,1\right]$ are random
amplitudes. We choose $\varepsilon_{b}=0.15$ and $n_{b}=100$. Finally,
given a set of random values $\boldsymbol{x}_{0,i}$ and $r_{i}$,
$b_{\max}$ is defined such that $\max_{\boldsymbol{x}_{K}}\left(b\left(\boldsymbol{x}_{K}\right)\right)=-10$,
so that the bottom remains submerged for small wave amplitudes, while
its variations are large enough to influence the flow significantly,
in particular through dispersion. The non-trivial bathymetry thus
constructed is shown in Fig. \ref{fig:test-linear-bathy}.

\begin{figure}[htbp]
\centering{}\includegraphics[width=0.5\textwidth]{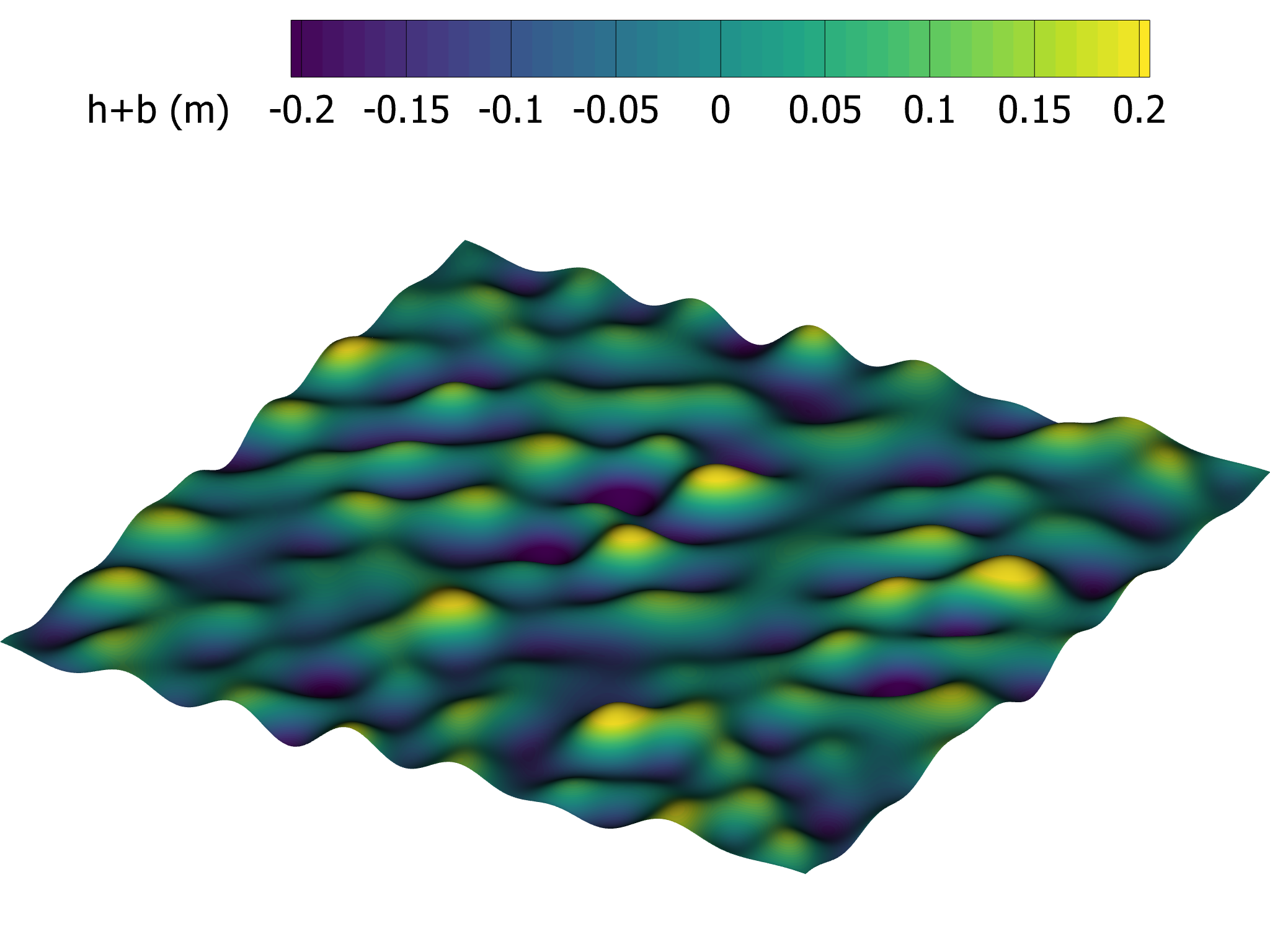}\includegraphics[width=0.5\textwidth]{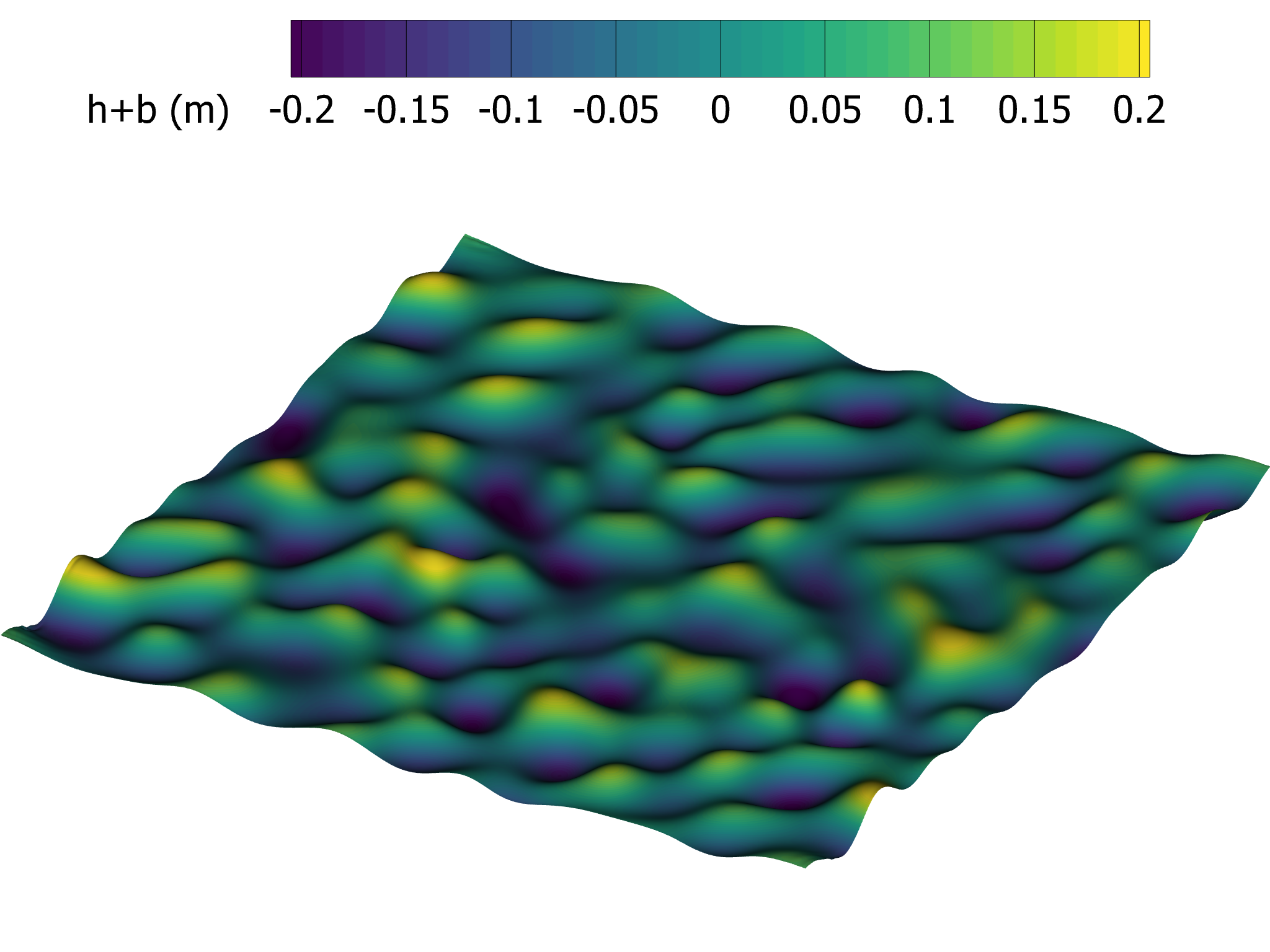}\caption{Test case of wave propagation in an approximately linear regime: initial
(left) and final (right) water surface $h+b$ deformations.}\label{fig:test-linear-init-final}
\end{figure}

Variables of the model $\left(h,u,v,w,p\right)$ are initialized according
to the linear wave solution derived in \ref{sec:linear_waves}:

\begin{subequations}
\begin{equation}
{\displaystyle h\left(\boldsymbol{x}_{K}\right)+b\left(\boldsymbol{x}_{K}\right)=\sum^{n_{w}}_{i=1}\:\varepsilon_{h,i}\:H\:\sin\left(\boldsymbol{k}_{i}\cdot\boldsymbol{x}_{K}-\varphi_{i}\right)}\,,
\end{equation}

\begin{equation}
u\left(\boldsymbol{x}_{K}\right)=\sum^{n_{w}}_{i=1}\:\frac{\omega_{i}k_{i,x}}{\left\Vert \boldsymbol{k}_{i}\right\Vert ^{2}}\:\varepsilon_{h,i}\:\sin\left(\boldsymbol{k}_{i}\cdot\boldsymbol{x}_{K}-\varphi_{i}\right)\,,
\end{equation}

\begin{equation}
v\left(\boldsymbol{x}_{K}\right)=\sum^{n_{w}}_{i=1}\:\frac{\omega_{i}k_{i,y}}{\left\Vert \boldsymbol{k}_{i}\right\Vert ^{2}}\:\varepsilon_{h,i}\:\sin\left(\boldsymbol{k}_{i}\cdot\boldsymbol{x}_{K}-\varphi_{i}\right)\,,
\end{equation}

\begin{equation}
w\left(\boldsymbol{x}_{K}\right)=\sum^{n_{w}}_{i=1}\:\frac{3}{2}\left(\frac{\omega_{i}}{\left\Vert \boldsymbol{k}_{i}\right\Vert ^{2}H}-\frac{g}{\omega_{i}}\right)\:\varepsilon_{h,i}\:\cos\left(\boldsymbol{k}_{i}\cdot\boldsymbol{x}_{K}-\varphi_{i}\right)\,,
\end{equation}

\begin{equation}
p\left(\boldsymbol{x}_{K}\right)=\sum^{n_{w}}_{i=1}\:\left(\frac{\omega^{2}_{i}}{\left\Vert \boldsymbol{k}_{i}\right\Vert ^{2}}-gH\right)\:\varepsilon_{h,i}\:\sin\left(\boldsymbol{k}_{i}\cdot\boldsymbol{x}_{K}-\varphi_{i}\right)\,,
\end{equation}
\end{subequations}

where the angular frequency $\omega_{i}$ is the solution of the dispersion
relation derived in \ref{sec:linear_waves} assuming the mean water
depth $H$ and a Mach number $M=0.2$, which leads to an acoustic
speed $a=5\sqrt{gH}\approx70\;\mathrm{m/s}$ for a standard gravitational
acceleration $g=9.81\;\mathrm{m/s^{2}}$, and a number of random waves
$n_{w}=100$. Note that this sum of linear waves is not an exact solution
over the variable bathymetry, but provides a suitable initial condition.
The sets of random values are an amplitude $\varepsilon_{h,i}\in\left[0,10^{-3}\right]$,
a wavenumber $\boldsymbol{k}_{i}$ such that $\boldsymbol{k}_{i}\in(2\pi/L)\left(\left\{ -8,\ldots,8\right\} ^{2}\setminus\left\{ (0,0)\right\} \right)$
and a phase shift $\varphi_{i}\in\left[0,2\pi\right]$. It turns out
that $\max_{i}\left(\left\Vert \boldsymbol{k}_{i}\right\Vert H\right)\approx1.4$,
which corresponds to a deviation of the phase velocity of about 20\%
compared to the shallow-water model. The simulation time is set to
ten periods of the smallest generated angular frequency $T=20\pi/\min_{i\in\left\llbracket 1,n_{w}\right\rrbracket }\left(\omega_{i}\right)\approx715\thinspace\mathrm{s}$.
The resulting initial and final water surfaces are shown in Fig. \ref{fig:test-linear-init-final}.

\begin{figure}[htbp]
\begin{centering}
\includegraphics[width=0.5\textwidth]{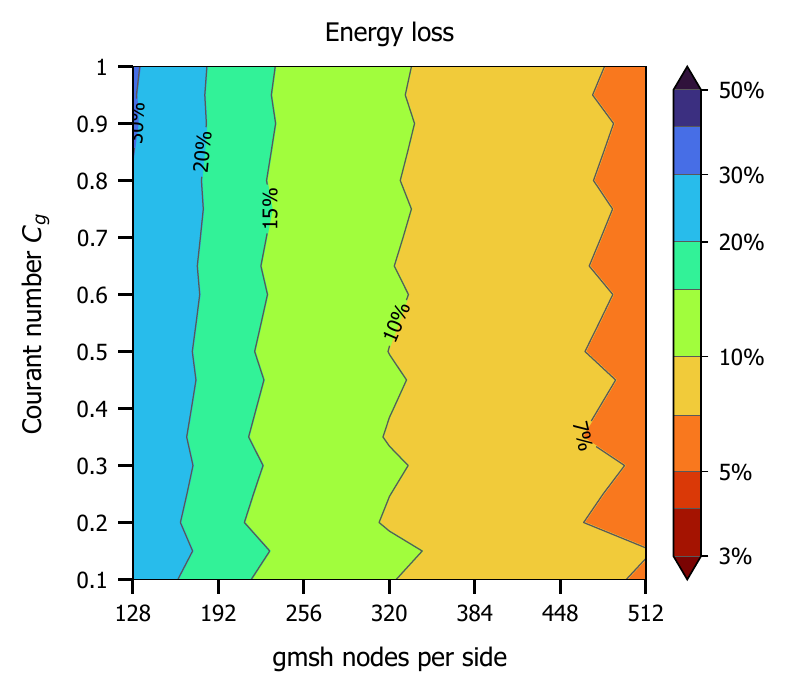}\includegraphics[width=0.5\textwidth]{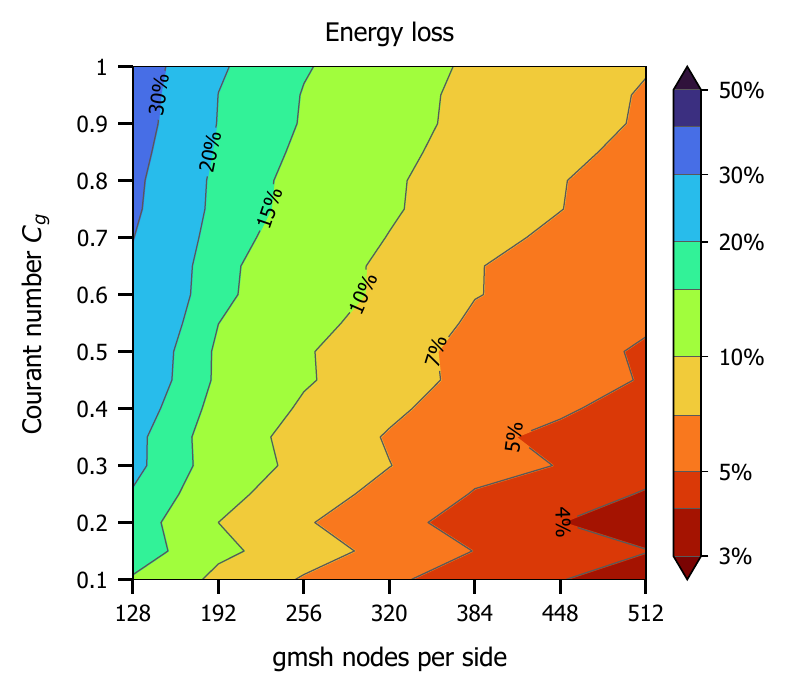}
\par\end{centering}
\centering{}\includegraphics[width=0.5\textwidth]{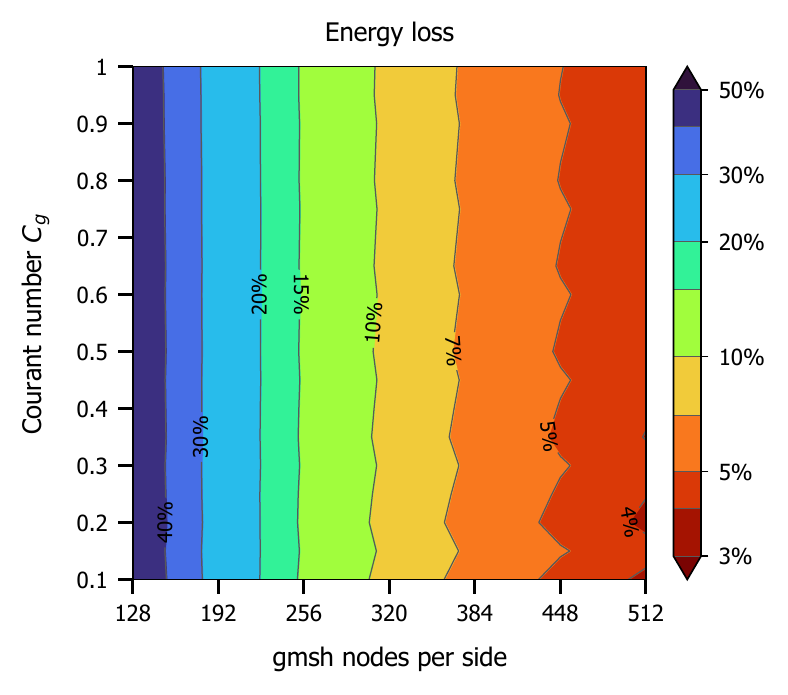}\includegraphics[width=0.5\textwidth]{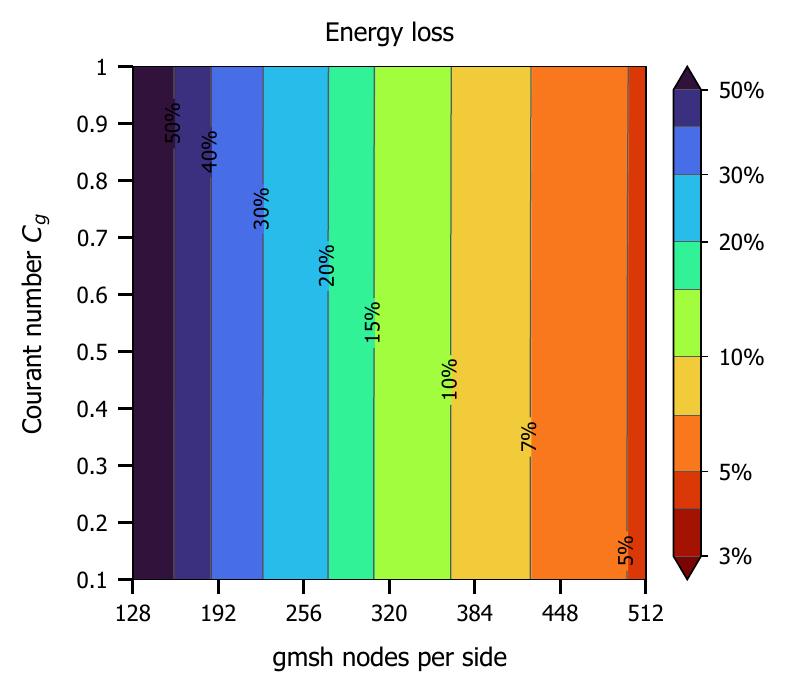}\caption{Contour plots of energy dissipation in percent as a function of mesh
size (see Table \ref{tab:gmsh-meshes}) and the Courant number; (top-left)
proposed splitting with the fully upwind CPR scheme ($\theta_{e}=1$,
$\gamma_{\ell}=1$); (top-right) proposed splitting with the calibrated
CPR scheme ($\theta_{e}$ from \eqref{eq:theta_mesh}, $\gamma_{\ell}=\Gamma(Fr_{e})$);
(bottom-left) proposed splitting with a Rusanov flux and MUSCL reconstructions;
(bottom-right) IMEX-SSP(3,3,2) strategy with a Rusanov flux and MUSCL
reconstructions for the gravity part, an implicit treatment of the
auxiliary variables and an explicit velocity divergence.}\label{fig:test-linear-energy}
\end{figure}

The software gmsh \citep{geuzaine2009} was used to generate seven
meshes with a number of nodes initially evenly distributed on each
side of the square box before applying the so-called Frontal-Delaunay
algorithm to perform a triangulation inside. The number of generated
triangles for each mesh is given in Table \ref{tab:gmsh-meshes}.
Since the dissipation of our numerical schemes depends on the Courant
number, we calculate the amount of energy loss for each generated
mesh and a range of the Courant number. The reported energy loss is
$100\bigl(1-\mathcal{E}_{\mathrm{wave}}(T)/\mathcal{E}_{\mathrm{wave}}(0)\bigr)$,
where
\[
\mathcal{E}^{n}_{\mathrm{wave}}=\mathcal{E}^{n}+\frac{g}{2}\sum_{K}m_{K}b^{2}_{K}
\]
is the energy relative to the lake-at-rest state $h+b=0$. This constant
shift leaves the discrete energy balance unchanged and normalizes
the loss by the initial wave energy. The results are summarized in
Fig. \ref{fig:test-linear-energy}, where four schemes are compared.
\begin{center}
\begin{table}[htbp]
\begin{centering}
\begin{tabular}{>{\raggedleft}p{0.3\textwidth}>{\raggedleft}p{0.3\textwidth}}
\toprule 
Number of gmsh nodes per side & Number of triangles\tabularnewline
\midrule 
$128$ & $37972$\tabularnewline
$192$ & $85356$\tabularnewline
$256$ & $151720$\tabularnewline
$320$ & $237014$\tabularnewline
$384$ & $341172$\tabularnewline
$448$ & $464360$\tabularnewline
$512$ & $606476$\tabularnewline
\bottomrule
\end{tabular}
\par\end{centering}
\caption{Test case of wave propagation in an approximately linear regime: number
of triangles generated by gmsh for each mesh used.}\label{tab:gmsh-meshes}
\end{table}
\par\end{center}

As expected, the dissipated energy decreases monotonically as the
mesh is refined, for all four schemes. In this nearly linear regime,
$Fr_{e}\lesssim10^{-2}$, so that $\theta_{e}\lesssim10^{-2}$ and
the calibrated scheme reduces in practice to a centred transport with
$\gamma_{\ell}=\Gamma(Fr_{e})\simeq1$. On coarse and moderate meshes,
it is by far the least dissipative of the four schemes: for $C_{g}\leq0.5$,
it loses two to four times less energy than the RK2-Heun/Rusanov/MUSCL
and IMEX schemes at the lowest resolution, and it remains the least
dissipative up to about $384$ nodes per side. Together with the lower
computational cost measured at a common Courant number (see the introduction
of this section), these results suggest a favourable compromise between
energy dissipation and cost for the CPR scheme on coarse and moderate
meshes.

The calibrated scheme also dissipates less as the Courant number decreases.
This is a valuable property in practice: on realistic meshes, the
time step is set by the most constraining cell, so that the local
Courant number is lower, often much lower, almost everywhere else,
precisely where the scheme is least dissipative. The fully upwind
variant is less sensitive, and almost insensitive on the finest meshes,
while the RK2-Heun/Rusanov/MUSCL and IMEX schemes are insensitive
to $C_{g}$. The sensitivity of the calibrated scheme does not come
from the transport weight, negligible here, but is consistent with
the time-step scaling of the gravity stabilization. Although $\gamma_{\ell}=1$
is the rest-state long-wave threshold, it does not cancel the leading-order
dissipation in every direction, even on a uniform Cartesian grid.
The small excess $\Gamma(Fr_{e})-1$ adds stabilization required for
uniform long-wave stability in the presence of a current.

Conversely, the spatial convergence of the CPR schemes at fixed Courant
number is slow, their stabilization being proportional to the time
step, hence to the mesh size: on the finest meshes, the second-order
RK2-Heun/Rusanov/MUSCL discretization becomes slightly less dissipative
than the calibrated scheme as soon as $C_{g}\gtrsim0.3$. Superimposed
on this overall trend, small oscillations of the dissipation with
the Courant number are observed, especially for the CPR variants.
They originate from the way the acoustic time step is adapted to the
gravity one: since the acoustic speed is larger than the characteristic
velocity of the gravity part, the acoustic step is sub-iterated and
its number of sub-steps is obtained by rounding up the ratio of the
two time steps, as explained in Section~\ref{subsec:time_stepping}.
Varying the  gravity time step with the Courant number therefore produces
jumps in the number of acoustic sub-steps, and hence small jumps in
the dissipated energy.

\subsection{Periodic travelling waves}\label{subsec:cnoidal}

We now want to investigate the behaviour of the numerical solver for
nonlinear wave regimes. To this end, we study the error accumulated
over a finite integration time when simulating periodic travelling
waves, also called cnoidal waves. Over a flat bottom, such a wave
is a stationary profile, for the water depth $h$ as for all other
model variables, translated at a constant speed $c$: the reference
solution at any time $t$ is therefore the initial profile shifted
by $ct$. Cnoidal waves form a family of solutions ranging from linear
to solitary waves, obtained by varying the two dimensionless parameters
$\varepsilon=A/H$, the nonlinearity parameter, and $\mu=k^{2}H^{2}$,
the dispersion or shallowness parameter. We do not have analytical
solutions of the Leucothea model \eqref{eq:model}, even in the soliton
case; \ref{sec:cnoidal} describes how to compute them numerically.
They are obtained by continuation in the Mach number $M=\sqrt{gH}/a$,
starting from the analytical cnoidal solution of the Serre--Green--Naghdi
model, to which the Leucothea model reduces as $a\to\infty$, that
is $M\to0$.

\begin{figure}[htbp]
\centering{}\includegraphics[width=0.8\textwidth]{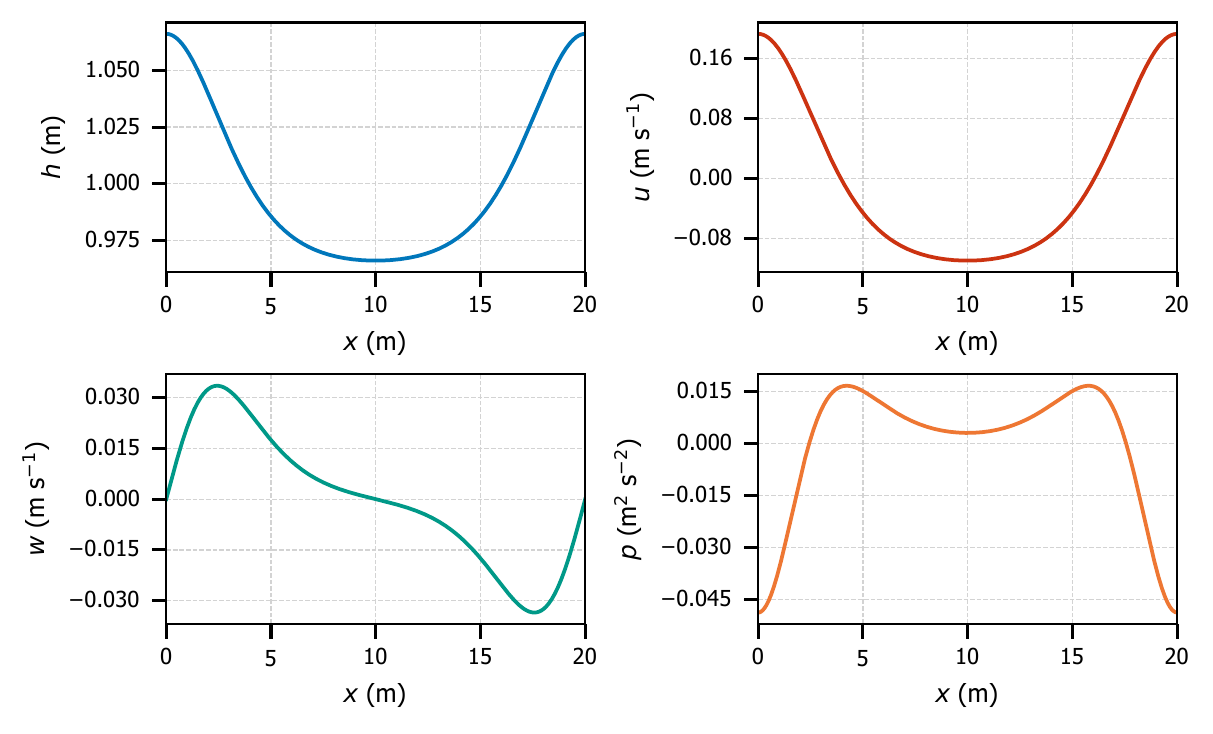}\caption{Cnoidal test case: solution computed numerically from~\eqref{eq:model}
with the algorithm of \ref{sec:cnoidal}.}\label{fig:test-cnoidal-sol-1d}
\end{figure}

\begin{figure}[htbp]
\centering{}\includegraphics[width=0.6\textwidth]{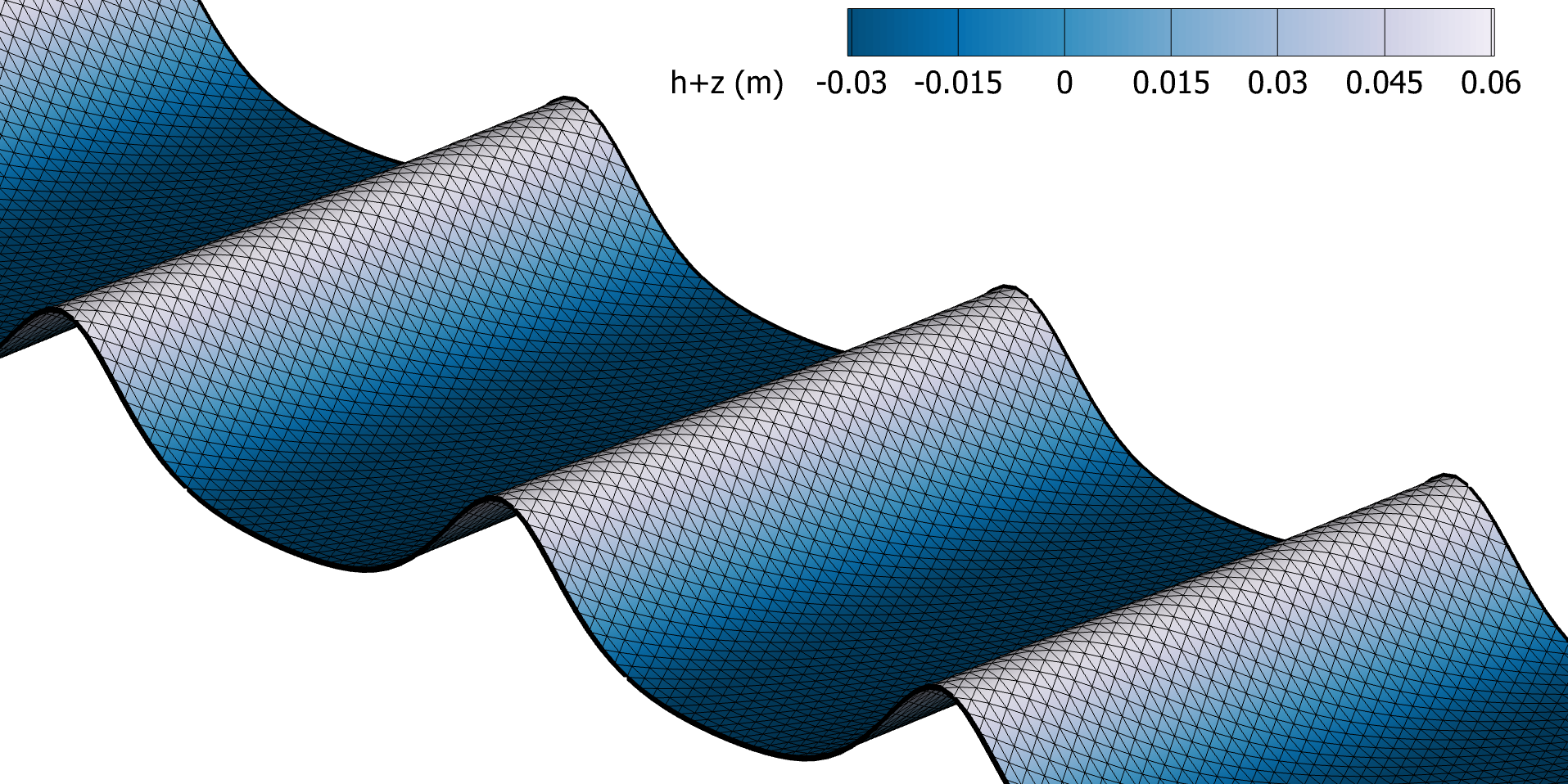}\caption{Cnoidal initial condition on the mesh with $35$ nodes per wavelength.}
\label{fig:test-cnoidal-sol-2d}
\end{figure}

For our test case, we arbitrarily consider a water depth $H=1\ \mathrm{m}$,
an amplitude $A=0.1\ \mathrm{m}$ and a wavelength $L=20\ \mathrm{m}$,
so that $\varepsilon=0.1$ and $\mu\approx0.098$. Such a parametrized
cnoidal solution is shown in Fig. \ref{fig:test-cnoidal-sol-1d} and
shows a periodic wave that lies between the linear regime --- the
profile is clearly not trigonometric --- and the solitary-wave limit,
with an elliptic parameter $m\approx0.94$ for the SGN profile ($m\to0$
and $m\to1$ in these two limits, respectively). We consider a rectangular
computational domain with periodic boundary conditions in both directions
such that the length is equal to $10$ wavelengths. The software gmsh
\citep{geuzaine2009} was used to generate meshes with a given number
of nodes distributed along the horizontal direction before applying
the so-called Frontal-Delaunay algorithm to perform triangulation
in the domain. The computations are pseudo-one-dimensional: the wave
propagates along $x$, and the transverse extent of the domain, in
the $y$-direction, is adjusted to the mesh size so as to keep well-shaped
triangles at a moderate computational cost. All the mesh data are
summarized in Table \ref{tab:gmsh-meshes-cnoidal}, and the initial
condition on the mesh with 35 nodes per wavelength is shown in Fig.
\ref{fig:test-cnoidal-sol-2d}. Finally, the Mach number $M=0.2$
is used throughout. For the CPR schemes, we compare $C_{g}=0.2,\ 0.5,\ 0.8$,
while $C_{g}=0.8$ is used for the other two schemes. The time of
simulation $T\approx64.21\ \mathrm{s}$ is set so that the wave travels
across the whole computational domain, that is $10$ wavelengths.
\begin{center}
\begin{table}[htbp]
\begin{centering}
\begin{tabular}{>{\centering}p{0.2\textwidth}>{\centering}p{0.2\textwidth}>{\centering}p{0.2\textwidth}>{\centering}p{0.2\textwidth}}
\toprule 
Nodes per wavelength & Triangles & Nodes per wavelength & Triangles\tabularnewline
\midrule 
$10$ & $9222$ & $45$ & $41498$\tabularnewline
$13$ & $12028$ & $58$ & $53388$\tabularnewline
$17$ & $15672$ & $75$ & $69028$\tabularnewline
$21$ & $19380$ & $97$ & $89264$\tabularnewline
$27$ & $24894$ & $124$ & $114144$\tabularnewline
$35$ & $32242$ & $160$ & $147294$\tabularnewline
\bottomrule
\end{tabular}
\par\end{centering}
\caption{Cnoidal test case: the unstructured meshes used for the convergence
study, with their resolution in nodes per wavelength and total triangle
count (the domain y-length is set to $40\,\delta$, where $\delta=L/n$
is the prescribed gmsh mesh size for $n$ nodes per wavelength).}
\label{tab:gmsh-meshes-cnoidal}
\end{table}
\par\end{center}

\begin{figure}[htbp]
\centering{}\includegraphics[width=0.6\textwidth]{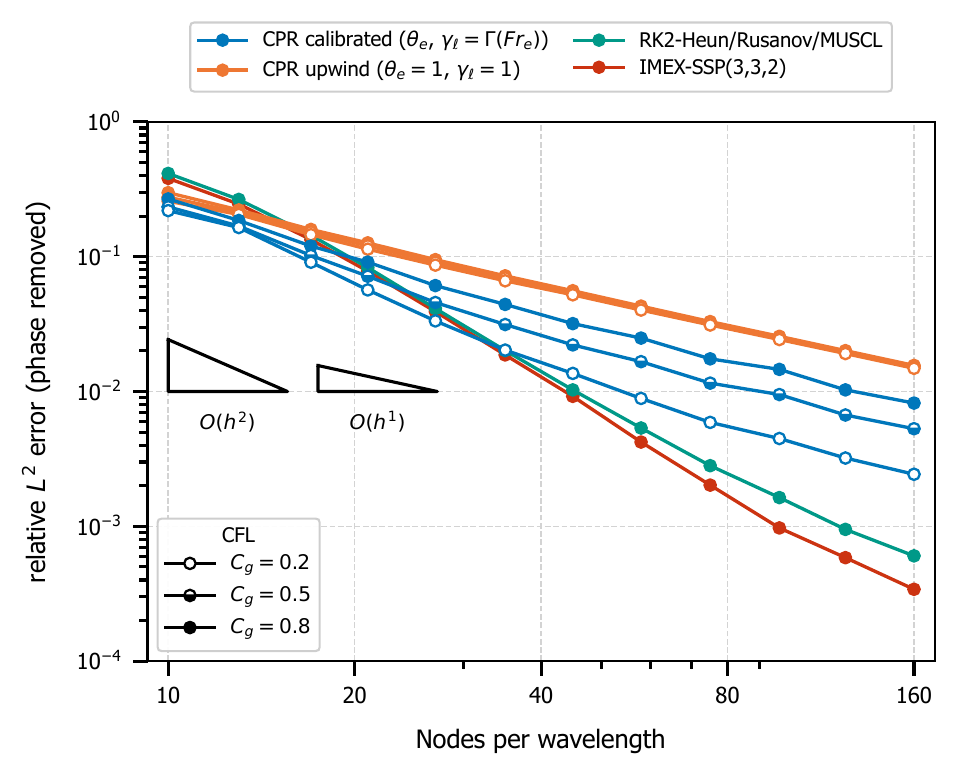}\caption{Convergence of the phase-aligned relative~$L^{2}$ error $e_{h}$
on the water depth after~$10$ wave periods, for the four schemes
and the meshes of Table~\ref{tab:gmsh-meshes-cnoidal}, at $C_{g}=0.2$,
$0.5$ and $0.8$ for the CPR schemes (open, half-filled and filled
markers) and $C_{g}=0.8$ for the other two.}
\label{fig:test-cnoidal-conv}
\end{figure}

The $L^{2}$ error norm on the water depth is computed at the end
of the simulation, the exact solution being the (periodic) initial
state since the wave travels an integer number of wavelengths. The
phase shift possibly introduced by the scheme is removed before computing
the error: the reference profile is translated along $x$ by the shift
that minimizes the $L^{2}$ error, found by a golden-section search
over $[-L/2,L/2]$. The reported error therefore measures the deformation
of the profile, the phase error being excluded. The phase-aligned
error is measured using the relative, cell-area-weighted discrete
$L^{2}$ norm
\[
e_{h}(\phi)=\left(\frac{\sum_{K}m_{K}\left|h_{K}-h_{\mathrm{ref}}(x_{K}+\phi)\right|^{2}}{\sum_{K}m_{K}\left|h_{\mathrm{ref}}(x_{K}+\phi)-H\right|^{2}}\right)^{1/2},
\]
where the periodic reference profile is evaluated by piecewise linear
interpolation. This normalization removes the dependence on the transverse
domain width for the one-dimensional reference wave. The shape error
is shown in Fig.~\ref{fig:test-cnoidal-conv}. What matters in practice
is the accuracy reached at the resolutions used in operational coastal
simulations, typically $20$ to $40$ nodes per wavelength. In this
range, the calibrated CPR scheme is competitive: it remains the most
accurate up to about $30$ nodes per wavelength at small Courant numbers,
its error decreasing with the Courant number as in the previous test
case, and it is the most accurate and the least dissipative of the
four schemes on coarser meshes. Its advantage is however much smaller
than in the nearly linear regime, where $Fr_{e}\lesssim10^{-2}$:
the CPR scheme loses ground as the Froude number becomes moderate,
and the second-order RK2-Heun/Rusanov/MUSCL and IMEX schemes take
over beyond a few tens of nodes per wavelength, the sooner the larger
the Courant number. At practical resolutions, the calibrated CPR scheme
therefore offers a favourable balance between computational cost and
profile accuracy, although its phase error remains larger than that
of the second-order schemes.

Beyond this range, the convergence of the CPR schemes is clearly too
slow: the calibrated scheme starts with a nearly second-order slope
but settles to first order, and the fully upwind variant, insensitive
to the Courant number, is first order throughout and the least accurate
on fine meshes. As for the nearly linear waves, improving the accuracy
of the gravity step is a priority for future work.

The phase error, removed from the shape error, is mostly a matter
of resolution. All schemes produce a slightly slow wave; on coarse
meshes, the CPR schemes accumulate over the ten periods a lag of a
sizable fraction of a wavelength, four to ten times larger than that
of the second-order schemes, which dominates their raw error, and
this lag becomes negligible on fine meshes.

The second-order convergence observed with RK2-Heun/Rusanov/MUSCL,
which uses the same acoustic step and splitting, indicates that the
acoustic step and the splitting do not dominate the error over the
tested resolution range.

\subsection{Vincent \& Briggs M1 experiment}

We now assess the solver against the experiments of Vincent \& Briggs
\citep{vincent1989} on wave propagation over a submerged elliptical
mound in a directional wave basin. Refraction over the mound focuses
the incident wave behind it, which produces a strong amplification
at the focal point, surrounded by lateral lobes. This benchmark tests
both the accuracy of the solver on a two-dimensional pattern of interfering
waves and the wave generation at the boundary.

\begin{figure}[htbp]
\centering{}\includegraphics[width=0.5\textwidth]{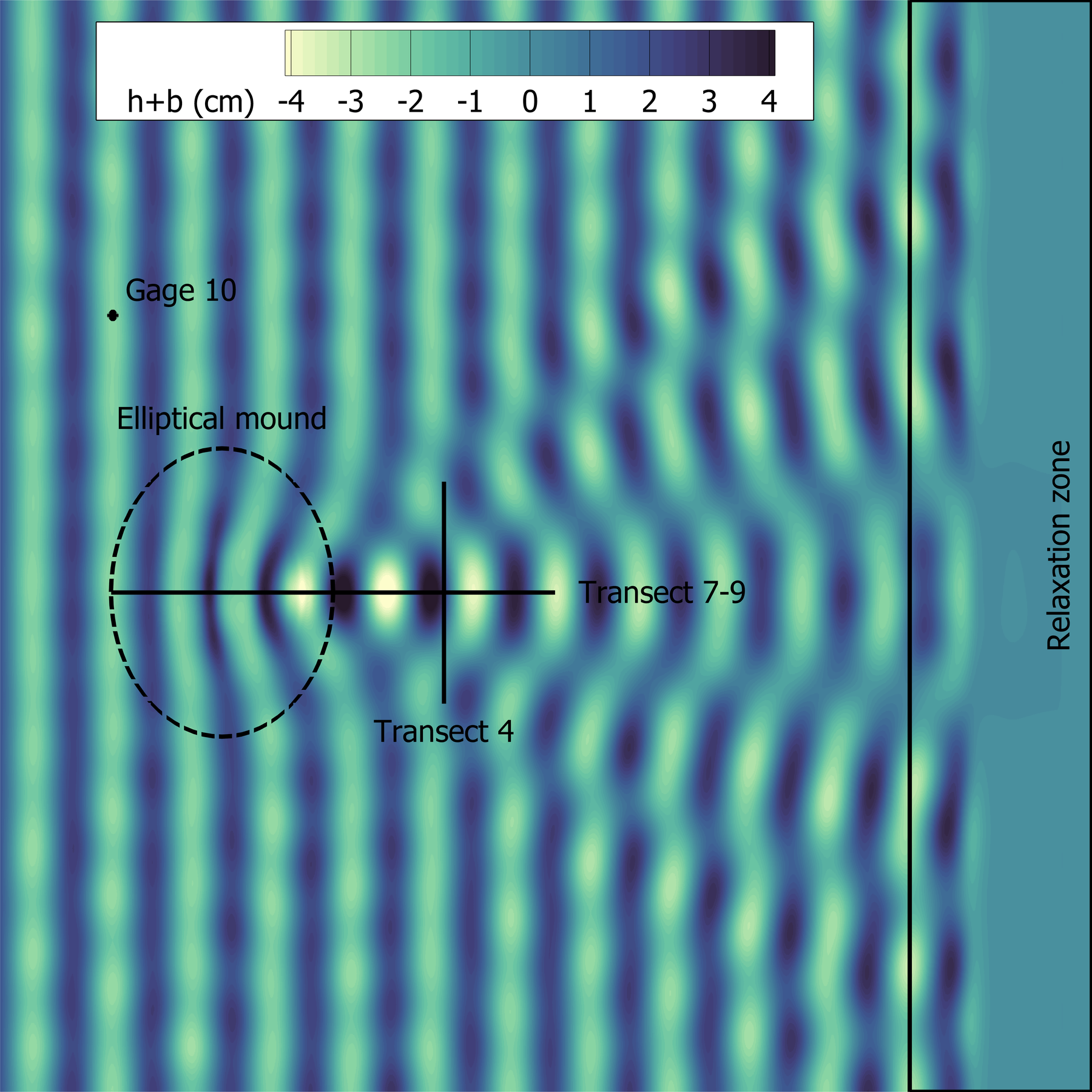}\caption{Water surface $h+b$ in the quasi-steady regime of the M1 Vincent
\& Briggs experiment: monochromatic wave propagation over a submerged
elliptical mound. The wave is generated at the left boundary and absorbed
by a relaxation layer at the right boundary.}\label{fig:Vincent-Briggs-field}
\end{figure}

The experimental directional wave channel is $25$ m long and $30$
m wide. The submerged elliptical mound is centred at $\left(x_{0},y_{0}\right)=\left(6.1\ \mathrm{m},13.72\ \mathrm{m}\right)$,
taking the origin at the bottom-left corner, its perimeter is defined
such that

\begin{equation}
\left(\frac{x-x_{0}}{3.05}\right)^{2}+\left(\frac{y-y_{0}}{3.96}\right)^{2}=1,
\end{equation}

and the bottom elevation inside this perimeter is

\begin{equation}
b\left(x,y\right)=-2h_{0}+0.762\sqrt{1-\left(\frac{x-x_{0}}{3.81}\right)^{2}-\left(\frac{y-y_{0}}{4.95}\right)^{2}}\,,
\end{equation}

where $h_{0}=45.72\ \mathrm{cm}$ is the water depth at rest around
the mound, where $b=-h_{0}$. Note that the water depth at the maximum
amplitude of the elliptical mound is $15.24\ \mathrm{cm}$, resulting
in a relatively large deformation of the bottom elevation. Among the
experiments, which range from monochromatic waves to directional spectra,
we consider the M1 case: a monochromatic wave of period $T=1.3\ \mathrm{s}$
and height $4.4\ \mathrm{cm}$. This gives $kh_{0}\approx1.27$ in
the constant-depth region and $kh\approx0.64$ over the crest of the
mound, with wavelengths of $2.26\ \mathrm{m}$ and $1.49\ \mathrm{m}$.
The wave is generated at the left boundary with the characteristic
condition of Section~\ref{subsec:bc}, its amplitude of $2.2\ \mathrm{cm}$
being raised over ten periods by a sine ramp. This amplitude corresponds
to the height measured at gauge 10, $\mathbf{x}_{10}=\left(3.05\ \mathrm{m},21.34\ \mathrm{m}\right)$,
in the experiment, and it is imposed directly without calibration.
At the right boundary, a relaxation layer of $5\ \mathrm{m}$, slightly
more than two wavelengths, absorbs the outgoing waves. The computations
use the common settings of Section~\ref{sec:validation}, with $C_{g}=0.2$
and $M=0.2$. After a transition period of approximately $40\ \mathrm{s}$,
a quasi-steady regime is reached, in the sense that only very low
amplitude and frequency modulations remain in the time series, most
likely due to wave interactions and reflections on the lateral walls.
Figure \ref{fig:Vincent-Briggs-field} shows the water surface field
at an arbitrary time in the quasi-steady regime, illustrating the
strong deformation of the incident monochromatic wave by refraction
over the submerged elliptical mound, since the phase speed depends
on the water depth.

The computational domain is discretized with a series of unstructured
triangular meshes of decreasing cell size, generated with gmsh \citep{geuzaine2009},
whose characteristics are summarized in Table~\ref{tab:vb-meshes}.
The prescribed gmsh sizes $\Delta$, from $0.4$ to $0.025\ \mathrm{m}$,
correspond to roughly $6$ to $90$ cells per incident wavelength,
and $4$ to $60$ over the crest, as orders of magnitude.
\begin{center}
\begin{table}[htbp]
\begin{centering}
\begin{tabular}{>{\centering}p{0.3\textwidth}>{\centering}p{0.3\textwidth}}
\toprule 
gmsh mesh size (m) & Number of triangles\tabularnewline
\midrule 
$0.4$ & $13106$\tabularnewline
$0.2$ & $52304$\tabularnewline
$0.1$ & $208152$\tabularnewline
$0.05$ & $831996$\tabularnewline
$0.025$ & $3327054$\tabularnewline
\bottomrule
\end{tabular}
\par\end{centering}
\caption{Vincent \& Briggs M1 test case: the unstructured meshes used for the
convergence study, with the prescribed gmsh mesh size (before triangulation)
and the resulting triangle count.}
\label{tab:vb-meshes}
\end{table}
\par\end{center}

\begin{figure}[htbp]
\centering{}\includegraphics[width=0.75\textwidth]{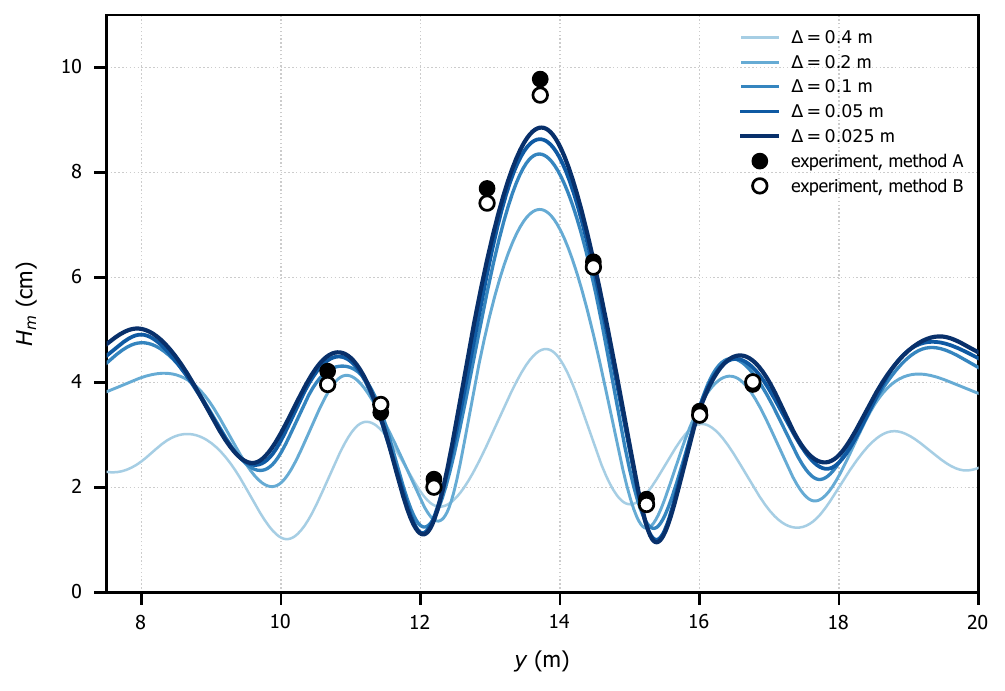}\caption{Mean wave height $H_{m}$ along transect 4 (see Fig.~\ref{fig:Vincent-Briggs-field})
for the M1 Vincent \& Briggs experiment, on the meshes of Table~\ref{tab:vb-meshes},
compared with the measurements analysed by the two methods of \citet[App. A and B]{vincent1989}.}\label{fig:vincent-briggs-mean-wave-height}
\end{figure}

A zero-up-crossing analysis gives the mean wave height $H_{m}$ over
the quasi-steady window, at regularly spaced points along transect
4, as in the experiment (Fig.~\ref{fig:vincent-briggs-mean-wave-height}).
The agreement is good overall. The lateral lobes, formed by the waves
diffracted on both sides of the focal zone, are captured in position
and amplitude from $\Delta=0.1\ \mathrm{m}$ on. The main discrepancy
is at the focal point, where the peak is underestimated by about 8\%
on the finest mesh. This peak still increases slowly with refinement,
by $0.3\ \mathrm{cm}$ and then $0.2\ \mathrm{cm}$ over the last
two halvings of $\Delta$, as observed for the cnoidal wave. This
points to a significant contribution of the numerical dissipation
of the gravity step. The dispersion relation of the model may also
play a part. Over the crest, the wave steepens (Ursell number  $H\lambda^{2}/h^{3}\approx28$,
with the values given above) and higher harmonics are released behind
the mound. The fundamental remains in a range where the model is accurate
to within 2\%, but the second harmonic reaches $kh\approx4.4$ in
the constant-depth region, where the model underestimates its phase
speed by more than 20\%. A fully nonlinear and fully dispersive potential-flow
model reproduces this peak more closely \citep{raoult_2019}. Separating
the two contributions requires finer meshes, or a version of Leucothea
with enhanced dispersive properties \citep{duran2024}, which will
be considered in future work.

\subsection{Realistic case: bay of Saint-Malo}\label{subsec:saint-malo}

\begin{figure}[htbp]
\begin{centering}
\includegraphics[width=0.48\textwidth]{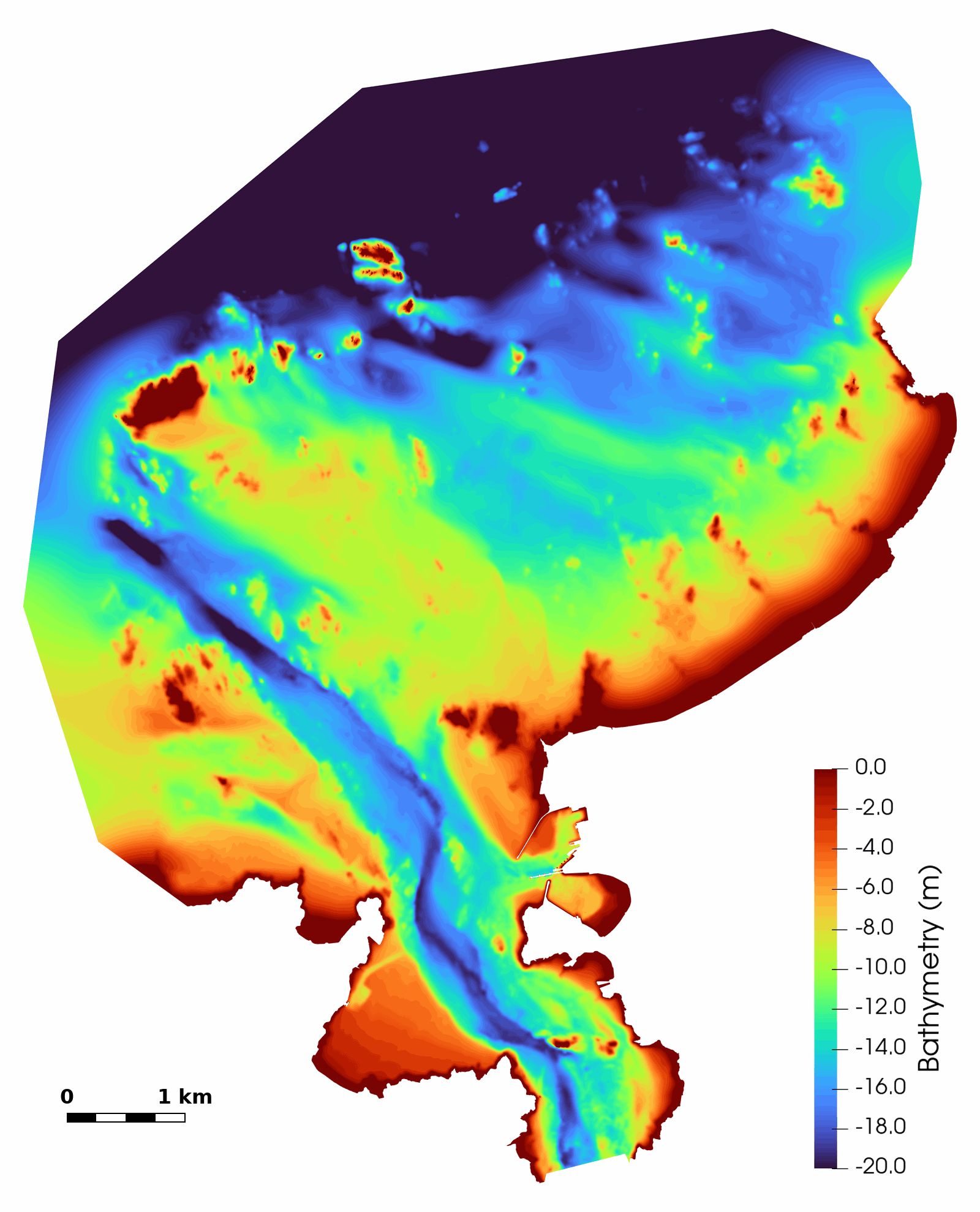}\includegraphics[width=0.48\textwidth]{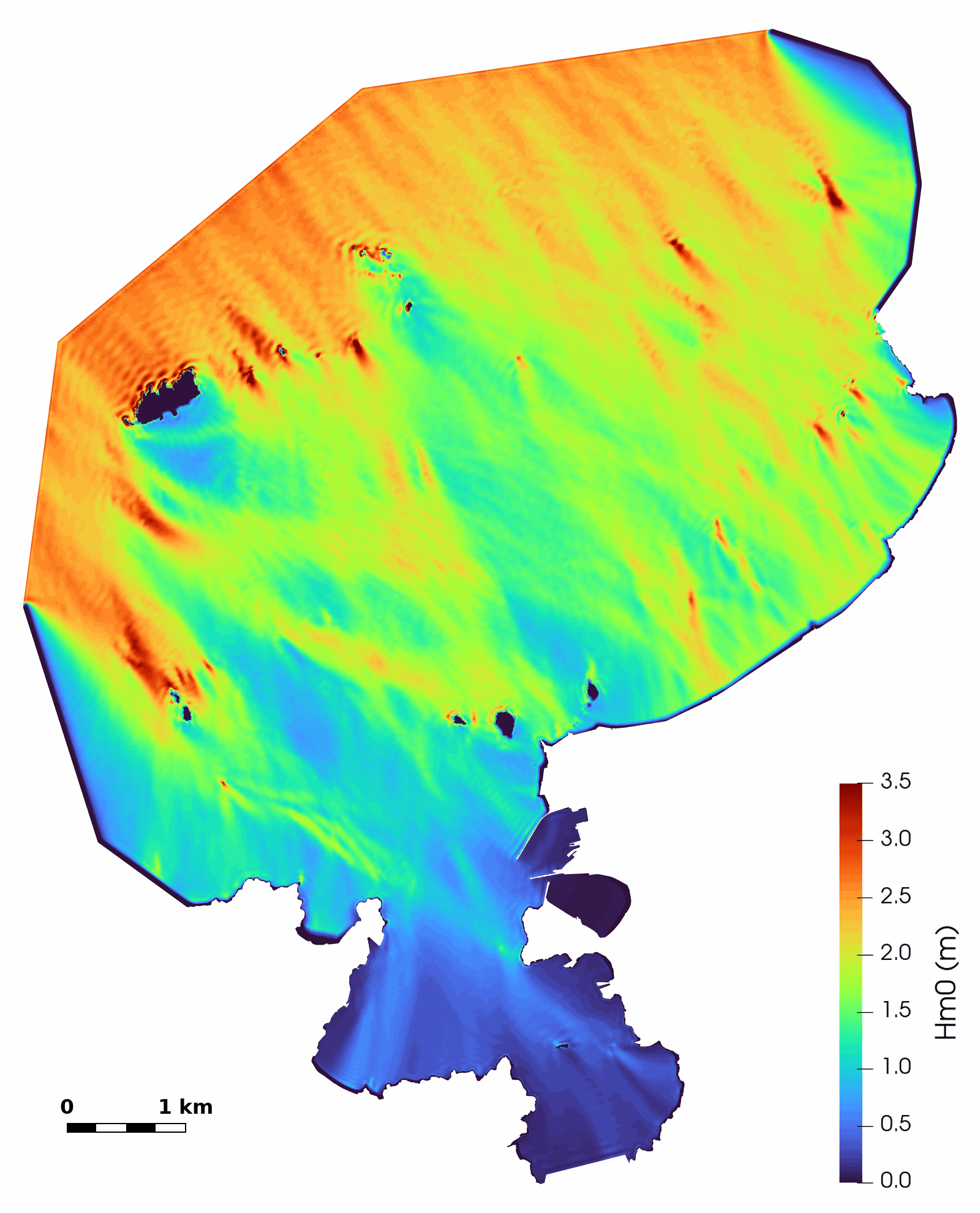}
\par\end{centering}
\begin{centering}
\includegraphics[width=0.48\textwidth]{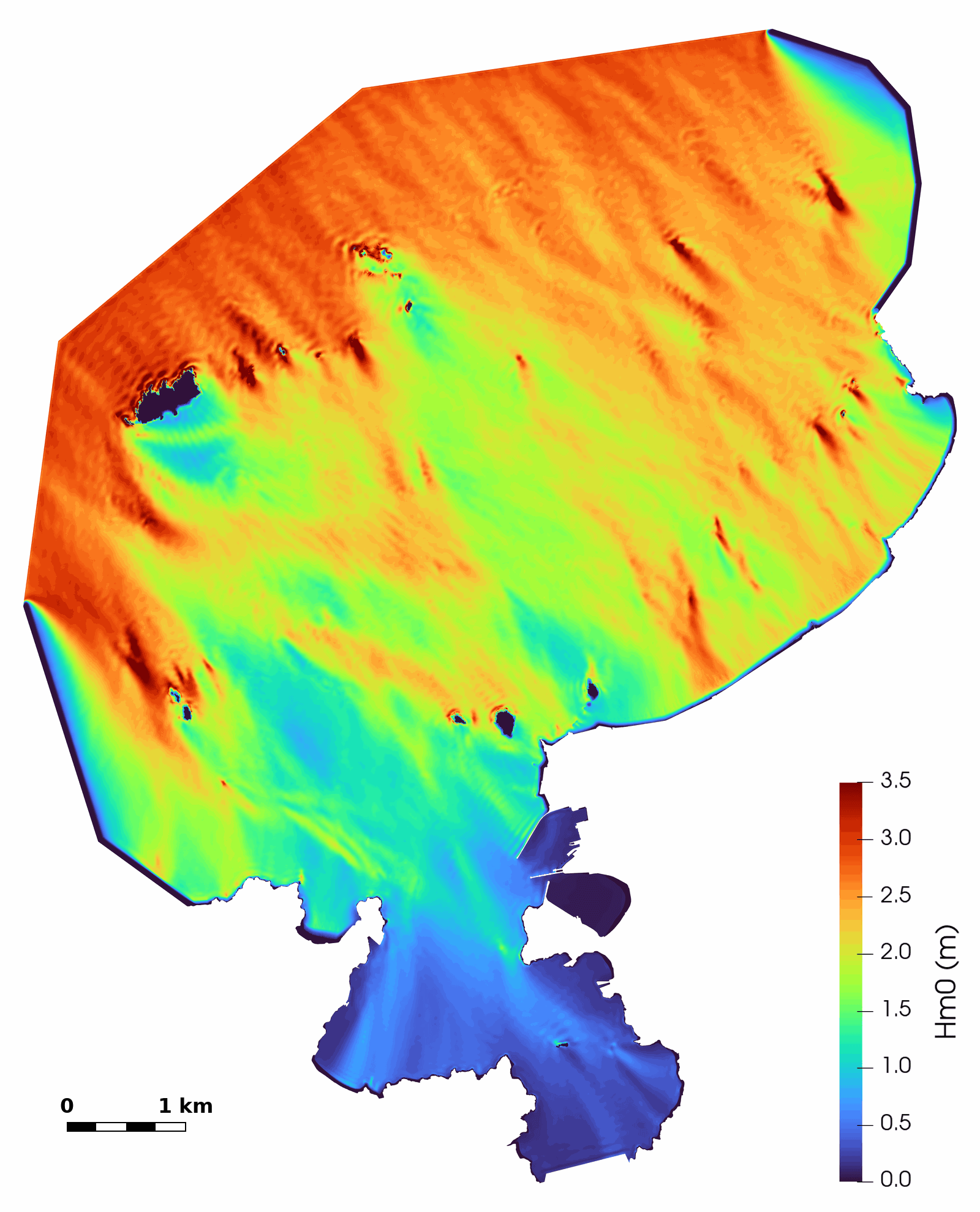}\includegraphics[width=0.48\textwidth]{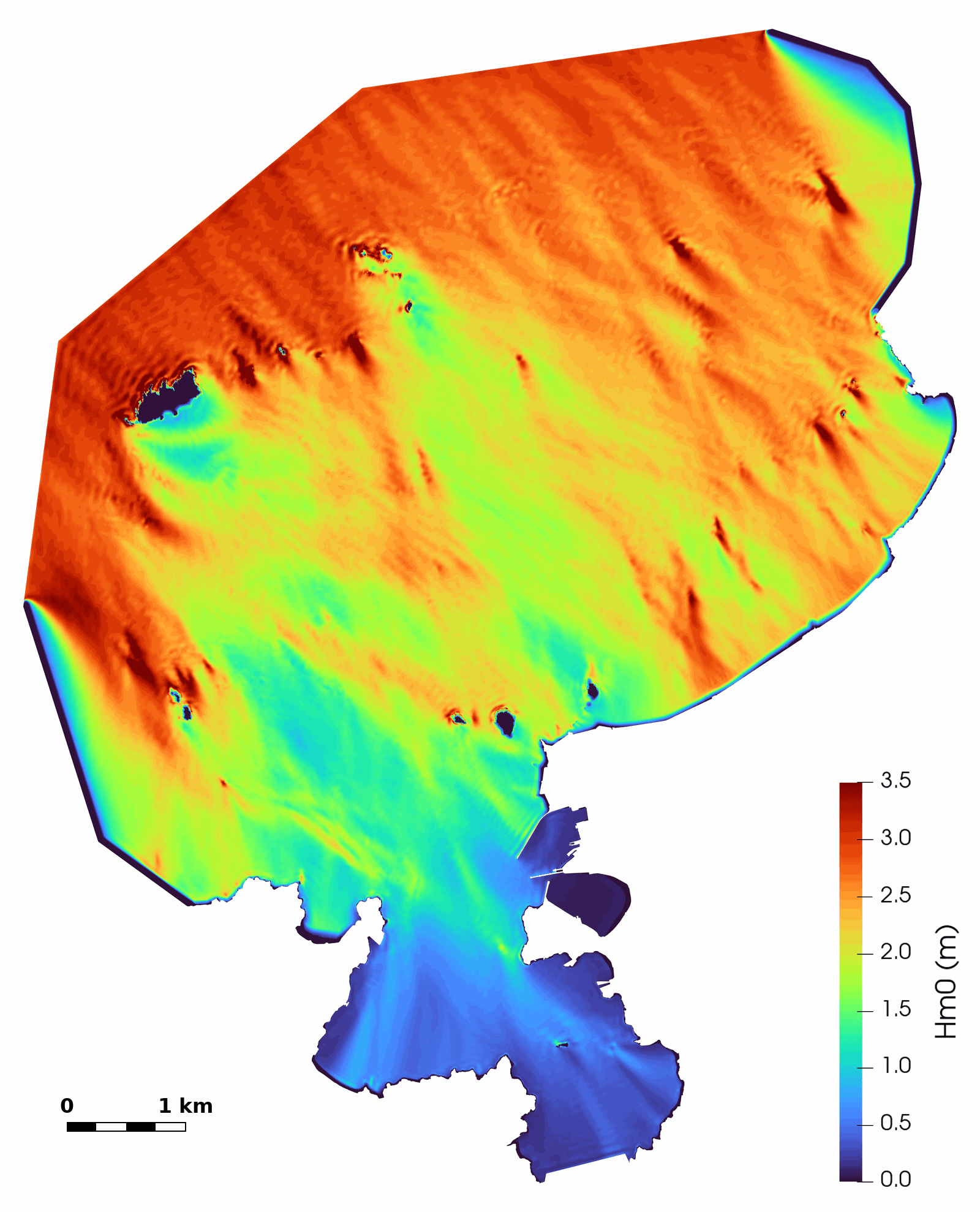}
\par\end{centering}
\caption{Saint-Malo test case: bathymetry and $H_{m0}=4\sigma_{\eta}$ wave-height
statistics for three mesh resolutions, under a directional TMA spectrum
imposed at the inlet ($f_{p}=0.08\ \mathrm{Hz}$, $H_{s}=3\ \mathrm{m}$)
and a two-hour simulation. (top-left) bathymetry ($0$ to $-20\ \mathrm{m}$);
(top-right) low resolution ($1.3$ M elements, $2$--$13\ \mathrm{m}$);
(bottom-left) intermediate resolution ($5.1$ M elements, $1$--$7\ \mathrm{m}$);
(bottom-right) high resolution ($20$ M elements, $0.5$--$3\ \mathrm{m}$).}\label{fig:saint-malo-hm0}
\end{figure}

As a final application, we consider the realistic coastal geometry
of the bay of Saint-Malo (Brittany, France), and more specifically
the Sillon, a long sandy beach protected along its whole length by
a sea wall (the so-called digue du Sillon). This beach is a recurrent
flooding hot spot during storm-tide events, which makes it a natural
target for the present solver. The computational domain stops at the
official SHOM coastline \citep{shom_limtm_2021}. The bathymetry is
obtained by merging the digital terrain model of the PAPI (French
flood prevention action programme) ($5\ \mathrm{m}$ resolution) \citep{shom_papi_2020}
and the Litto3D data set ($1\ \mathrm{m}$ resolution) \citep{shom_litto3d_2024},
the merged model being then uniformly smoothed to filter out the highest
spatial frequencies and remove sub-grid noise. It features the strong,
irregular variations characteristic of this macro-tidal coast: rocky
shoals, tidal channels and a complex coastline (see the bathymetry
panel of Fig. \ref{fig:saint-malo-hm0}, ranging from $0$ to $-20\ \mathrm{m}$).
This case is well beyond the academic configurations considered so
far: it combines real bathymetry, a complex domain boundary and a
directional wave forcing, and it is precisely the kind of coastal-engineering
situation the present solver is intended for.
\begin{center}
\begin{table}[htbp]
\begin{centering}
\begin{tabular}{>{\centering}p{0.175\textwidth}>{\centering}p{0.15\textwidth}>{\centering}p{0.17\textwidth}>{\centering}p{0.2\textwidth}>{\centering}p{0.18\textwidth}}
\toprule 
Resolution & Elements (M) & Cell size (m) & Nodes per wavelength & Time step (ms)\tabularnewline
\midrule 
Low & $\sim1.3$ & $\sim2\text{--}13$ & $\sim10$ & $\sim18$\tabularnewline
Intermediate & $\sim5.1$ & $\sim1\text{--}7$ & $\sim20$ & $\sim8$\tabularnewline
High & $\sim20$ & $\sim0.5\text{--}3$ & $\sim40$ & $\sim3.5$\tabularnewline
\bottomrule
\end{tabular}
\par\end{centering}
\caption{Saint-Malo test case: the three unstructured meshes used for the convergence
study, with their number of elements, cell-size range, resolution
in nodes per wavelength and resulting time step at $C_{g}=0.2$.}\label{tab:saint-malo-meshes}
\end{table}
\par\end{center}

A directional TMA spectrum is imposed along the offshore boundary
with the characteristic condition of Section~\ref{subsec:bc}. Its
parameters are a peak frequency $f_{p}=0.08\ \mathrm{Hz}$, a significant
wave height $H_{s}=3\ \mathrm{m}$, a peak enhancement factor $\gamma=3.3$
($\sigma_{a}=0.07$, $\sigma_{b}=0.09$), a cut-off frequency of $0.2\ \mathrm{Hz}$,
and a depth of 20 m for the TMA transformation. The spectrum is discretized
with $1000$ frequencies and $360$ directions. The directional distribution
is $\cos^{2s}\bigl((\theta-\theta_{m})/2\bigr)$, with $s=25$, typical
of swell, and a mean propagation direction $\theta_{m}=-50^{\circ}$
from the eastward axis, that is, waves coming from the north-west.
Relaxation layers of width 100 m absorb the outgoing waves along the
two lateral open boundaries, and the rest of the boundary, along the
coastline, is a solid wall. The simulation covers two hours of physical
time with the common settings of Section~\ref{sec:validation}, with
$M=0.2$ and $C_{g}=0.2$. This value of $C_{g}$ reduces the upwind
part $\theta_{e}$ of the gravity scheme, since $\theta_{e}$ scales
with the time step. The quantity of interest is the significant wave
height $H_{m0}=4\sigma_{\eta}$, where $\sigma_{\eta}$ is the standard
deviation of the free-surface elevation in each cell. It is accumulated
with Welford's online algorithm over the last $90$ minutes, after
a $30$-minute spin-up. $H_{m0}$ is the standard output of coastal
wave models, and the quantity measured by in situ and satellite altimetry
instruments.

The physical setting remains academic: neither tidal forcing nor a
wave-breaking model is included. Bottom friction is applied to the
velocity after the gravity step, $h\,\partial_{t}\boldsymbol{u}=-\bigl(C_{l}+C_{b}(h)\lVert\boldsymbol{u}\rVert\bigr)\boldsymbol{u}$,
with a small linear coefficient $C_{l}=10^{-5}\ \mathrm{m/s}$ and
a quadratic coefficient deduced from a depth-averaged logarithmic
velocity profile of roughness length $z_{0}=10^{-3}\ \mathrm{m}$,
$C_{b}(h)=\bigl[\kappa(h-z_{0})/\bigl(h\ln(h/z_{0})-h+z_{0}\bigr)\bigr]^{2}$,
$\kappa=0.41$ (about $2.5\times10^{-3}$ at $10\ \mathrm{m}$ depth),
increased smoothly below $h=0.1\ \mathrm{m}$; it is integrated implicitly
and only decreases the kinetic energy. In the shallowest areas, where
the incident waves would break, the numerical dissipation of the scheme
is sufficient to keep the computation stable. Wave heights are not
physically controlled there, however, so the mesh comparison below
assesses the numerical convergence of the statistics, not the physical
realism of the surf zone. Further validations including wave-breaking
models will be considered in future work.

To assess the convergence of these statistics with the mesh resolution,
three unstructured meshes of increasing density are compared, all
generated with gmsh \citep{geuzaine2009} using the Frontal-Delaunay
algorithm to guarantee high-quality triangular elements. The local
mesh size is prescribed as a function of both the distance to the
coastline and the target wavelength of linear Airy waves. The coarsest
mesh resolves about $10$ nodes per wavelength for Airy waves at $0.1\ \mathrm{Hz}$,
in the upper part of the spectrum, and the two finer meshes are obtained
by halving the characteristic mesh size; their characteristics are
summarized in Table~\ref{tab:saint-malo-meshes}. The resulting $H_{m0}$
fields are shown in Fig.~\ref{fig:saint-malo-hm0}.

At low resolution, the $H_{m0}$ map already reveals the diagonal
banded structure caused by refraction over rocky shoals and wave sheltering
by the islands. However, the incoming wave field is already weaker
near the open boundary than prescribed by the TMA forcing, and wave
heights progressively decrease towards the coast. At intermediate
resolution, the refraction patterns and wave shadows become more sharply
defined, while the incident wave field is stronger and remains energetic
over a much larger area. At high resolution, the contrast between
exposed areas ($H_{m0}>3\ \mathrm{m}$ offshore) and sheltered areas
($H_{m0}<0.5\ \mathrm{m}$ inside the bay) is sharp, and the focusing
on the rocky tips and channels is clearly resolved. The solver thus
captures, on a real and geometrically complex domain, the refraction
and sheltering processes that govern the local wave climate.

\begin{figure}[htbp]
\begin{centering}
\includegraphics[width=0.48\textwidth]{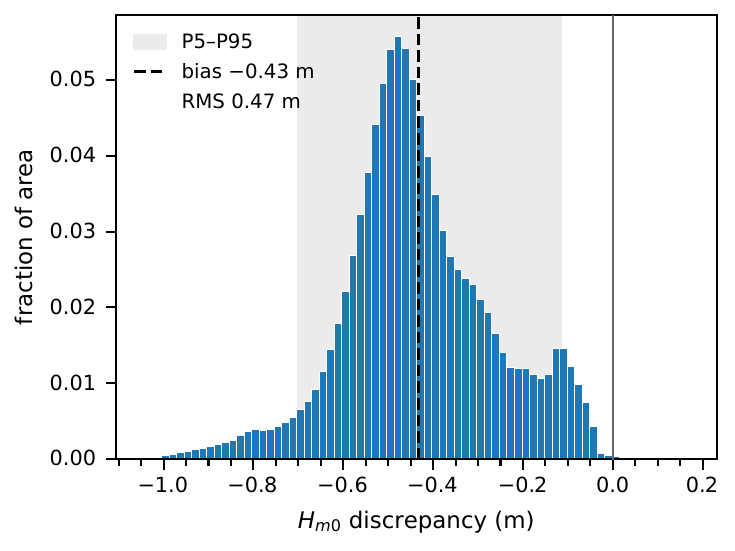}\includegraphics[width=0.48\textwidth]{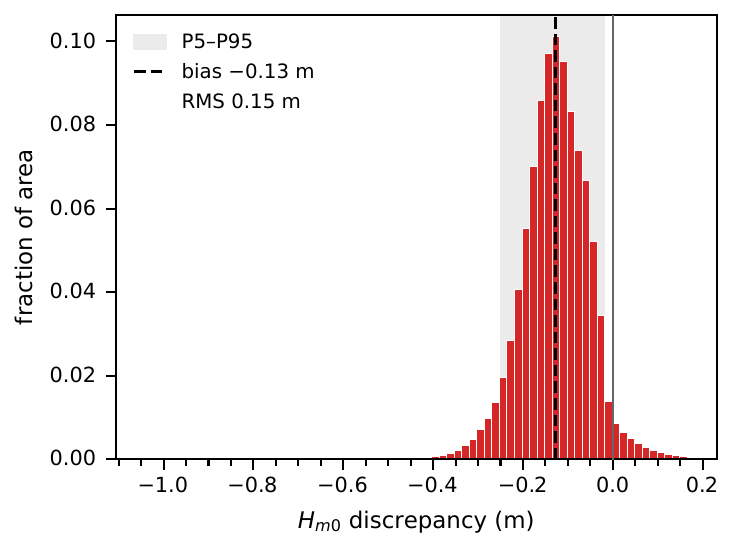}
\par\end{centering}
\caption{Saint-Malo test case: area-weighted distribution of the pointwise
$H_{m0}$ discrepancy with respect to the high-resolution mesh, over
the cells where $H_{m0}\protect\geq0.1\ \mathrm{m}$ on the fine mesh.
(left) low-resolution mesh; (right) intermediate mesh. The dashed
line marks the area-weighted bias and the shaded band the $[\mathrm{P5},\mathrm{P95}]$
interval.}\label{fig:saint-malo-error-hist}
\end{figure}

This convergence is quantified by the pointwise discrepancy of $H_{m0}$
between each coarse mesh and the high-resolution one (Fig.~\ref{fig:saint-malo-error-hist}).
To compare the fields on the same mesh, each coarse-mesh field is
transferred to the high-resolution mesh by assigning to each fine
cell the value of the coarse triangle containing its centroid. This
uses the piecewise-constant coarse field, without interpolation between
coarse-cell values. The discrepancy is then computed as the transferred
coarse value minus the high-resolution value in each fine cell. The
distributions are weighted by the fine-cell areas and restricted to
cells where $H_{m0}\geq0.1\ \mathrm{m}$ on the fine mesh, about 98\%
of the domain area. On the low-resolution mesh, the discrepancy is
markedly negative: the bias is $-0.43\ \mathrm{m}$ (−22\%) and the
root-mean-square $0.47\ \mathrm{m}$. On the intermediate mesh, these
figures drop to $-0.13\ \mathrm{m}$ (−7\%) and $0.15\ \mathrm{m}$,
and the distribution sharpens around its mode. Both biases being measured
against the finest mesh, an error of the form $C\delta^{p}$ gives
a ratio $(4^{p}-1)/(2^{p}-1)=2^{p}+1$ between them; the observed
ratio, about 3.3, corresponds to $p\approx1.2$ ($p\approx1.1$ for
the root-mean-square). Assuming that this rate persists under further
refinement, extrapolation suggests that the high-resolution mesh may
still underestimate the converged $H_{m0}$ by about $0.13/(2^{p}-1)\approx0.1\ \mathrm{m}$,
roughly 5\%. As for the Vincent \& Briggs and cnoidal cases, the error
is a systematic underestimation that decreases slowly with refinement,
which points to the numerical dissipation of the gravity step. The
intermediate mesh already offers a good compromise between accuracy
and cost for this configuration. With a throughput of about $0.1\ \mu\mathrm{s}$
per cell, per time step and per core, the two-hour simulation takes
about 12 min on 72 cores, 43 min on 180 cores and 3.2 h on 360 cores
for the low, intermediate and high resolutions, that is, 10 and 3
times faster than real time for the first two meshes.

\section{Conclusion}

In this work we have introduced an entropy-stable splitting strategy
for the numerical solution of the Leucothea model, a hyperbolic relaxation
of the SGN equations derived under a weak compressibility assumption,
in two dimensions and on unstructured meshes. The splitting separates
the model into a shallow-water subsystem, supplemented with two passively
transported scalars, and an acoustic subsystem carrying the dispersive
effects, each of the two inheriting its own exact energy balance.
The gravity part relies on the CPR scheme, in its calibrated or fully
upwind variants, while the acoustic part is advanced by a scheme whose
implicit $w$--$p$ coupling reduces, after elementary algebra, to
a fully explicit update. Local energy estimates have been established
for both steps in the nonlinear, fully discrete framework, yielding
sufficient conditions on the time step and stabilization parameters.
For the acoustic step on a uniform square mesh at constant depth and
for $\theta_{p}=1/2$, these conditions coincide with the necessary
and sufficient conditions for energy contraction of the linearized
scheme over all Fourier modes. The nonlinear estimate is therefore
sharp in this configuration. The complementary linear analyses in
the energy norm guide the calibration of these constants and quantify
the numerical dissipation as a function of the Courant number, the
spatial resolution and the Froude number.

The resulting solver is fully explicit, local, and free of any elliptic
step or global linear system. The fast acoustic waves are handled
by sub-iterating the explicit acoustic update after each gravity time
step. The overall cost remains three to six times that of a classical
shallow-water solver, in the measured configurations. On academic
benchmarks of wave propagation in a nearly linear regime over a non-trivial
bathymetry and of nonlinear periodic travelling waves computed from
the model itself, the proposed splitting combined with the calibrated
CPR scheme offers a favourable balance between computational cost
and profile accuracy at the resolutions used in practice, about 20
to 40 nodes per wavelength. Its advantage is largest at low Froude
and Courant numbers, but phase errors remain significant on coarse
meshes, and the second-order schemes take over on finer meshes. The
Vincent \& Briggs experiment confirms the ability of the solver to
reproduce refraction--focusing patterns with good accuracy. Finally,
the Saint-Malo configuration, built on real bathymetric and coastline
data with a directional spectrum forcing, demonstrates the feasibility
of wave-resolving simulations on a large coastal domain with complex
geometry, open boundaries and wet/dry areas.

Across the test cases, the mesh convergence of the calibrated CPR
scheme is slow. The analyses and numerical comparisons indicate that
the gravity discretization governs most of the numerical dissipation
in the regimes investigated, while the acoustic contribution remains
small. This interpretation is supported by the second-order convergence
obtained with another gravity discretization using the same splitting
and acoustic step. The proposed acoustic treatment can therefore be
retained while choosing a shallow-water solver suited to the accuracy,
robustness and cost required by the application; CPR is one such choice.
An energy guarantee for the complete method requires compatible energy
control from both steps. A higher-order reconstruction in the gravity
step is a natural next step to improve the convergence. Particular
attention should be paid to the pressure-gradient discretization and
its associated stabilization on unstructured meshes, while preserving
the lake-at-rest equilibrium and discrete energy control. Further
work will also examine how mesh geometry affects the dissipation of
the gravity step, with the aim of reducing its directional dependence
while preserving the energy estimates.

A complementary direction concerns the model itself: including the
full topography terms and the enhanced dispersive properties while
preserving an exact energy equation remains non-trivial in two dimensions,
although a first step has recently been taken in one dimension \citep{duran2024}.
The underestimation observed at the focal point of the Vincent \&
Briggs experiment provides a concrete motivation for this extension.
Improving the model's dispersive accuracy would complement the reduction
of numerical dissipation, particularly for the shorter-wave components
generated by nonlinear propagation.

Further perspectives include the treatment of wave breaking, required
for a complete description of nearshore dynamics \citep{kazakova_new_2019,richard_new_2019,duran_modelling_2020,hung2025}.
As numerical dissipation is reduced, an explicit representation of
the physical energy loss associated with breaking becomes increasingly
important. The dissipation of the numerical scheme cannot substitute
for a breaking model in the surf zone, and coupling a physical breaking
model to the solver should account consistently for the additional
energy transfers and dissipation.

Boundary conditions should also be improved to better accommodate
dispersive wave generation and absorption. Reducing spurious reflections
for oblique waves and broad directional spectra, and accounting consistently
for energy exchanges through open boundaries, are important objectives
for coastal applications. Finally, the explicit and local structure
of the method supports the exploitation of massively parallel architectures,
notably GPUs, with the aim of reaching real-time operational forecasting
of coastal hazards. 

\section*{CRediT authorship contribution statement}

\textbf{Frédéric Couderc:} Conceptualization, Formal analysis, Investigation,
Software, Validation, Visualization, Writing -- original draft, Writing
-- review and editing. \textbf{Arnaud Duran:} Conceptualization,
Formal analysis, Methodology, Writing -- original draft, Writing
-- review and editing. \textbf{Maria Kazakova:} Formal analysis,
Methodology. \textbf{Bastien Lombard:} Data curation, Investigation,
Resources, Visualization. \textbf{Pascal Noble:} Formal analysis,
Project administration, Supervision. \textbf{Rémy Baraille:} Conceptualization,
Formal analysis, Funding acquisition, Project administration, Resources,
Supervision. \textbf{Jean-Paul Vila:} Conceptualization, Formal analysis,
Methodology.

\section*{Declaration of competing interest}

The authors declare that they have no known competing financial interests
or personal relationships that could have appeared to influence the
work reported in this paper.

\section*{Data and code availability}

The numerical results were obtained with the Tolosa-lct solver developed
at the Institut de Mathématiques de Toulouse. Its source code is hosted
in a private repository; the source code, the post-processing scripts
and the processed data supporting the findings of this study are available
from the corresponding author upon reasonable request. The bathymetric
and coastline data are SHOM products \citep{shom_papi_2020,shom_litto3d_2024,shom_limtm_2021}.
The figures were generated directly from the numerical output of Tolosa-lct
by scripts written with the assistance of OpenAI GPT-6 Astra and Anthropic
Claude Opus 5.5. These scripts read the result files without modifying
the data, and no figure was created or altered by a generative image
tool.

\section*{Declaration of generative AI and AI-assisted technologies in the
manuscript preparation process}

The numerical model, the schemes, the Tolosa-lct solver and the test
cases result from several years of work by the authors. During the
final preparation of this manuscript, the authors used OpenAI GPT-6
Astra and Anthropic Claude Opus 5.5 to harmonize the notation and
presentation of the linear and nonlinear stability analyses, to check
the proofs and, in a few places, to complete them or sharpen their
bounds, to improve the language, structure and consistency of the
text, and to write the Python scripts used to post-process the numerical
results and generate the figures (see Data and code availability).
After using these tools, the authors reviewed and edited the content
as needed. In particular, every step of the proofs was checked by
the authors, and the figures were verified against the underlying
data. The authors take full responsibility for the content of the
published article.

\section*{Acknowledgements}

Funding: This work was supported by SHOM {[}research contract No 21CP06,
supporting F. Couderc and P. Noble, and research contract No 21CP05,
supporting A. Duran{]} and by the French National Research Agency
(ANR) {[}project MOTIONS, grant number ANR-23-CE56-0006-02, awarded
to A. Duran; project BOURGEONS, grant number ANR-23-CE40-0014-01,
supporting M. Kazakova{]}. Computations were performed using HPC resources
from CALMIP {[}grant number 2026-P1234{]}.

\appendix

\section{Global energy balance for the gravity step}

\label{app:gravity-global-balance}

We derive the balance \eqref{eq:gravity-global-balance} for the discrete
energy \eqref{eq:gravity-total-energy}, under the periodic boundary
conditions assumed above. Each interface is counted once, with $\boldsymbol{n}_{e}$
oriented from $K$ to $K_{e}$. With the half-jump convention, pairing
the contributions of the two neighbouring cells gives, for any conservative
oriented flux $F^{n}_{e,K}$, 
\begin{equation}
\sum_{K}z^{n}_{K}\sum_{e\in\partial K}m_{e}F^{n}_{e,K}=-2\sum_{e}m_{e}\llbracket z^{n}\rrbracket_{e}F^{n}_{e}.\label{eq:gravity-edge-pairing}
\end{equation}
The same identity holds with a scalar product for vector-valued quantities.

\paragraph{Temporal energy increment}

For any scalar or vector state $z_{K}$, the exact identity 
\begin{equation}
\frac{h^{n+1}_{K}}{2}\|z^{n+1}_{K}\|^{2}-\frac{h^{n}_{K}}{2}\|z^{n}_{K}\|^{2}=z^{n}_{K}\cdot\Delta(hz)_{K}-\frac{\|z^{n}_{K}\|^{2}}{2}\Delta h_{K}+\frac{h^{n+1}_{K}}{2}\|\Delta z_{K}\|^{2}\label{eq:gravity-kinetic-identity}
\end{equation}
separates the first-order variation from the non-negative temporal
remainder. Likewise, since $\phi^{n}_{K}=g(h^{n}_{K}+b_{K})$, 
\[
\frac{g}{2}\bigl((h^{n+1}_{K})^{2}-(h^{n}_{K})^{2}\bigr)+gb_{K}\Delta h_{K}=\phi^{n}_{K}\Delta h_{K}+\frac{g}{2}(\Delta h_{K})^{2}.
\]
With $\boldsymbol{Z}^{n}_{K}=(\boldsymbol{u}^{n}_{K},2w^{n}_{K}/\sqrt{3},p^{n}_{K}/a)$,
the energy density reads 
\[
h^{n}_{K}E^{n}_{K}=\frac{1}{2}h^{n}_{K}\|\boldsymbol{Z}^{n}_{K}\|^{2}+\frac{g}{2}(h^{n}_{K})^{2}+gb_{K}h^{n}_{K}.
\]
Applying \eqref{eq:gravity-kinetic-identity} to $\boldsymbol{Z}_{K}$
and adding the potential-energy increment therefore yields 
\begin{equation}
\mathcal{E}^{n+1}-\mathcal{E}^{n}=\sum_{K}m_{K}\left[\boldsymbol{Z}^{n}_{K}\cdot\Delta(h\boldsymbol{Z})_{K}+\left(\phi^{n}_{K}-\frac{\|\boldsymbol{Z}^{n}_{K}\|^{2}}{2}\right)\Delta h_{K}\right]+\mathcal{R}_{\Delta t}.\label{eq:gravity-first-variation}
\end{equation}
Here $\mathcal{R}_{\Delta t}$ is exactly the remainder in \eqref{eq:gravity-time-remainder}.

\paragraph{Transport contribution}

The transport flux of $h\boldsymbol{Z}$ is 
\[
\mathcal{G}^{n}_{e}=\overline{\boldsymbol{Z}^{n}}_{e}\mathcal{F}^{*,n}_{e}-\theta^{n}_{e}|\mathcal{F}^{*,n}_{e}|\llbracket\boldsymbol{Z}^{n}\rrbracket_{e}.
\]
Substituting the mass and transport updates into \eqref{eq:gravity-first-variation}
and applying \eqref{eq:gravity-edge-pairing} gives the transport
contribution 
\[
2\Delta t\sum_{e}m_{e}\left(\llbracket\boldsymbol{Z}^{n}\rrbracket_{e}\cdot\mathcal{G}^{n}_{e}-\frac{1}{2}\llbracket\|\boldsymbol{Z}^{n}\|^{2}\rrbracket_{e}\mathcal{F}^{*,n}_{e}\right).
\]
Since $\llbracket\|\boldsymbol{Z}^{n}\|^{2}\rrbracket_{e}=2\overline{\boldsymbol{Z}^{n}}_{e}\cdot\llbracket\boldsymbol{Z}^{n}\rrbracket_{e}$,
the centred transport terms cancel interface by interface, leaving
\begin{equation}
-2\Delta t\sum_{e}m_{e}\theta^{n}_{e}|\mathcal{F}^{*,n}_{e}|\,\|\llbracket\boldsymbol{Z}^{n}\rrbracket_{e}\|^{2}.\label{eq:gravity-transport-contribution}
\end{equation}

\paragraph{Potential energy and pressure force}

The mass update gives the potential-energy contribution 
\[
\sum_{K}m_{K}\phi^{n}_{K}\Delta h_{K}=2\Delta t\sum_{e}m_{e}\llbracket\phi^{n}\rrbracket_{e}\mathcal{F}^{*,n}_{e}.
\]
The pressure force acts only on the horizontal momentum. Its contribution
is 
\begin{align*}
 & -\Delta t\sum_{K}h^{n}_{K}\boldsymbol{u}^{n}_{K}\cdot\sum_{e\in\partial K}m_{e}(\overline{\phi^{n}}_{e}-\Lambda^{n}_{e})\boldsymbol{n}_{e,K}\\
 & \qquad=2\Delta t\sum_{e}m_{e}\bigl(\llbracket h^{n}\boldsymbol{u}^{n}\rrbracket_{e}\cdot\boldsymbol{n}_{e}\bigr)(\overline{\phi^{n}}_{e}-\Lambda^{n}_{e}),
\end{align*}
where \eqref{eq:gravity-edge-pairing} has again been used. Inserting
$\mathcal{F}^{*,n}_{e}=\overline{h^{n}\boldsymbol{u}^{n}}_{e}\cdot\boldsymbol{n}_{e}-\boldsymbol{\Pi}^{n}_{e}\cdot\boldsymbol{n}_{e}$,
the sum of these two contributions separates into a centred part and
the stabilization terms. The centred part is 
\begin{align*}
 & 2\Delta t\sum_{e}m_{e}\left(\llbracket\phi^{n}\rrbracket_{e}\overline{h^{n}\boldsymbol{u}^{n}}_{e}+\overline{\phi^{n}}_{e}\llbracket h^{n}\boldsymbol{u}^{n}\rrbracket_{e}\right)\cdot\boldsymbol{n}_{e}\\
 & \hspace{25mm}=2\Delta t\sum_{e}m_{e}\llbracket\phi^{n}h^{n}\boldsymbol{u}^{n}\rrbracket_{e}\cdot\boldsymbol{n}_{e}.
\end{align*}
It vanishes by geometric closure of each cell: 
\[
2\sum_{e}m_{e}\llbracket\phi^{n}h^{n}\boldsymbol{u}^{n}\rrbracket_{e}\cdot\boldsymbol{n}_{e}=-\sum_{K}\phi^{n}_{K}h^{n}_{K}\boldsymbol{u}^{n}_{K}\cdot\underbrace{\sum_{e\in\partial K}m_{e}\boldsymbol{n}_{e,K}}_{=\,0}.
\]
The remaining terms are therefore 
\begin{equation}
-2\Delta t\sum_{e}m_{e}\left[(\boldsymbol{\Pi}^{n}_{e}\cdot\boldsymbol{n}_{e})\llbracket\phi^{n}\rrbracket_{e}+\Lambda^{n}_{e}\llbracket h^{n}\boldsymbol{u}^{n}\rrbracket_{e}\cdot\boldsymbol{n}_{e}\right].\label{eq:gravity-stabilization-contribution}
\end{equation}
Combining \eqref{eq:gravity-transport-contribution} and \eqref{eq:gravity-stabilization-contribution}
with the temporal remainder gives \eqref{eq:gravity-global-balance}.
For nonperiodic boundaries, the associated boundary contributions
must be retained in the pairing of cell contributions. 

\section{Energy dissipation for the acoustic subsystem}

\label{sec:dissipation}

\subsection{Preliminaries}

Before establishing the balance \eqref{eq:acoustic_estim1} and Proposition~\ref{prop:acoustic_estim2},
it is necessary to introduce some useful notations, as well as duality
formulas. First, for any sequence of interface scalar or vectorial
values $(\omega_{e})_{e\in\mathbb{E}}$, $(\boldsymbol{\omega}_{e})_{e\in\mathbb{E}}$,
we introduce the discrete operators

\[
\begin{array}{l}
{\displaystyle \boldsymbol{\nabla}_{K}\left(\omega\right)=\frac{1}{m_{K}}\sum_{e\in\partial K}\omega_{e}\,\boldsymbol{n}_{e,K}\,m_{e}},\\
{\displaystyle \boldsymbol{\nabla}_{K}\cdot\left(\boldsymbol{\omega}\right)=\frac{1}{m_{K}}\sum_{e\in\partial K}\boldsymbol{\omega}_{e}\cdot\boldsymbol{n}_{e,K}\,m_{e}},
\end{array}
\]
and make use of the following specific notation for centred gradient
and divergence: 
\begin{equation}
\boldsymbol{\nabla}_{K}\left(\overline{\omega}\right)={\displaystyle \frac{1}{m_{K}}\sum_{e\in\partial K}\left(\frac{\omega_{K}+\omega_{K_{e}}}{2}\right)\,\boldsymbol{n}_{e,K}\,m_{e}}={\displaystyle \frac{1}{m_{K}}\sum_{e\in\partial K}\left(\frac{\omega_{K_{e}}-\omega_{K}}{2}\right)\,\boldsymbol{n}_{e,K}\,m_{e}}\,,%
\label{eq:centered_gradient}
\end{equation}

\begin{equation}
\boldsymbol{\nabla}_{K}\cdot\left(\overline{\boldsymbol{\omega}}\right)={\displaystyle \frac{1}{m_{K}}\sum_{e\in\partial K}\left(\frac{\boldsymbol{\omega}_{K}+\boldsymbol{\omega}_{K_{e}}}{2}\right)\cdot\boldsymbol{n}_{e,K}\,m_{e}}={\displaystyle \frac{1}{m_{K}}\sum_{e\in\partial K}\left(\frac{\boldsymbol{\omega}_{K_{e}}-\boldsymbol{\omega}_{K}}{2}\right)\cdot\boldsymbol{n}_{e,K}\,m_{e}}\,.\label{eq:centered_divergence}
\end{equation}

Also, in the forthcoming developments, we will regularly employ the
notation below: 
\begin{equation}
\mathcal{S}_{K}\left(\omega\right)=\frac{1}{m_{K}}\sum_{e\in\partial K}\omega_{e}m_{e}\,.\label{eq:SK}
\end{equation}
The notations \eqref{eq:centered_gradient}, \eqref{eq:centered_divergence}
and \eqref{eq:SK} admit natural generalizations to interface products.
We will for instance denote 
\[
\boldsymbol{\nabla}_{K}\left(\overline{a},b\right)={\displaystyle \frac{1}{m_{K}}\sum_{e\in\partial K}\overline{a}_{e}b_{e}\,\boldsymbol{n}_{e,K}\,m_{e}}\,,
\]
where we recall that $\overline{a}_{e}=\left(a_{K_{e}}+a_{K}\right)/2$.
Similar extensions will stand for 
\[
\boldsymbol{\nabla}_{K}\cdot\left(\overline{\boldsymbol{a}},b\right)={\displaystyle \frac{1}{m_{K}}\sum_{e\in\partial K}b_{e}\overline{\boldsymbol{a}}_{e}\cdot\boldsymbol{n}_{e,K}\,m_{e}}\,\,\text{ and }\,\,\mathcal{S}_{K}\left(\overline{a},b\right)={\displaystyle \frac{1}{m_{K}}\sum_{e\in\partial K}\overline{a}_{e}b_{e}\,m_{e}}\,.
\]
For vector-valued edge quantities, $\mathcal{S}_{K}\left(a,(\boldsymbol{\delta}b),(\boldsymbol{\delta}c)\right):=(1/m_{K})\sum_{e\in\partial K}a_{e}\,(\boldsymbol{\delta}_{e}b)\cdot(\boldsymbol{\delta}_{e}c)\,m_{e}$,
with the scalar product of $\mathbb{R}^{2}$. The following lemma
provides some useful duality formulas satisfied by the operators introduced
previously:
\begin{lem}
\label{lem:duality} Let $(a_{e})_{e\in\mathbb{E}}$ be a sequence
of interface scalar values, $(b_{K})_{K\in\mathbb{T}}$, $(c_{K})_{K\in\mathbb{T}}$
two sequences of collocated scalar values, and $(\boldsymbol{a}_{K})_{K\in\mathbb{T}}$
a collocated sequence of two-component vectors. The following duality
identities hold: 
\begin{equation}
{\displaystyle b_{K}\boldsymbol{\nabla}_{K}\cdot\left(a,\boldsymbol{\delta}c\right)}=-\mathcal{S}_{K}\left(a,\left(\boldsymbol{\delta}b\right),\left(\boldsymbol{\delta}c\right)\right)+\mathcal{R}^{2}_{K}\left(a,b,c\right)\,,\label{eq:duality2}
\end{equation}
with $\mathcal{R}^{2}_{K}\left(a,b,c\right)=(1/m_{K})\sum_{e\in\partial K}\mathcal{F}^{2}_{K,K_{e}}\left(a,b,c\right)m_{e}\quad,\quad\mathcal{F}^{2}_{K,K_{e}}\left(a,b,c\right)=a_{e}\overline{b}_{e}(\boldsymbol{\delta}_{e}c)\cdot\boldsymbol{n}_{e,K}$.

\begin{equation}
{\displaystyle c_{K}\boldsymbol{\nabla}_{K}\cdot\left(a,\boldsymbol{\delta}b\right)=b_{K}\boldsymbol{\nabla}_{K}\cdot\left(a,\boldsymbol{\delta}c\right)+\mathcal{R}^{3}_{K}\left(a,b,c\right)}\,,\label{eq:duality3}
\end{equation}
with $\mathcal{R}^{3}_{K}\left(a,b,c\right)=(1/m_{K})\sum_{e\in\partial K}\mathcal{F}^{3}_{K,K_{e}}\left(a,b,c\right)m_{e}\quad,\quad\mathcal{F}^{3}_{K,K_{e}}\left(a,b,c\right)=a_{e}\left(\overline{c}_{e}(\boldsymbol{\delta}_{e}b)-\overline{b}_{e}(\boldsymbol{\delta}_{e}c)\right)\cdot\boldsymbol{n}_{e,K}$.

\begin{equation}
{\displaystyle \boldsymbol{a}_{K}\cdot\boldsymbol{\nabla}_{K}(\overline{b})+b_{K}\boldsymbol{\nabla}_{K}\cdot(\overline{\boldsymbol{a}})=\boldsymbol{\nabla}_{K}\cdot(\overline{\boldsymbol{a}},\overline{b})+\mathcal{R}^{4}_{K}(\boldsymbol{a},b)}\,,\label{eq:duality4}
\end{equation}
with $\mathcal{R}^{4}_{K}(\boldsymbol{a},b)=(1/m_{K})\sum_{e\in\partial K}\mathcal{F}^{4}_{K,K_{e}}\left(\boldsymbol{a},b\right)m_{e}\quad,\quad\mathcal{F}^{4}_{K,K_{e}}\left(\boldsymbol{a},b\right)=-(\boldsymbol{\delta}_{e}b)(\delta_{e}\boldsymbol{a})\cdot\boldsymbol{n}_{e,K}$.

The fluxes $(\mathcal{F}^{i})_{i=2,3,4}$ are consistent with zero
in the sense $\mathcal{F}^{i}_{K,K}=0$ for all $K\in\mathbb{T}$,
and satisfy the conservativity property $\mathcal{F}^{i}_{K,K_{e}}=-\mathcal{F}^{i}_{K_{e},K}$
for $i=2,3,4$.
\end{lem}
\begin{proof}
Using the jump notations \eqref{eq:notations_interface}, \eqref{eq:notations_operators}
and definitions \eqref{eq:centered_gradient}, \eqref{eq:centered_divergence}
the proofs of these results are elementary. For \eqref{eq:duality2}
we just write 
\[
\begin{array}{lll}
b_{K}\boldsymbol{\nabla}_{K}\cdot\left(a,\boldsymbol{\delta}c\right) & = & {\displaystyle \frac{1}{m_{K}}\left(\sum_{e\in\partial K}a_{e}\left(\frac{b_{K}-b_{K_{e}}}{2}\right)\left(\boldsymbol{\delta}_{e}c\right)\cdot\boldsymbol{n}_{e,K}m_{e}+\sum_{e\in\partial K}a_{e}\left(\frac{b_{K}+b_{K_{e}}}{2}\right)\left(\boldsymbol{\delta}_{e}c\right)\cdot\boldsymbol{n}_{e,K}m_{e}\right)}\,.\end{array}
\]
From this, \eqref{eq:duality3} follows directly by symmetry: 
\[
\begin{array}{lll}
c_{K}\boldsymbol{\nabla}_{K}\cdot\left(a,\boldsymbol{\delta}b\right) & = & {\displaystyle -\mathcal{S}_{K}\left(a,\left(\boldsymbol{\delta}b\right),\left(\boldsymbol{\delta}c\right)\right)+\frac{1}{m_{K}}\sum_{e\in\partial K}a_{e}\overline{c}_{e}(\boldsymbol{\delta}_{e}b)\cdot\boldsymbol{n}_{e,K}m_{e}}\\
 & = & {\displaystyle b_{K}\boldsymbol{\nabla}_{K}\cdot\left(a,\boldsymbol{\delta}c\right)+\frac{1}{m_{K}}\sum_{e\in\partial K}a_{e}\left(\overline{c}_{e}(\boldsymbol{\delta}_{e}b)-\overline{b}_{e}(\boldsymbol{\delta}_{e}c)\right)\cdot\boldsymbol{n}_{e,K}m_{e}}.
\end{array}
\]
Regarding the last identity \eqref{eq:duality4}, it results from
the sum of 
\begin{align*}
\boldsymbol{a}_{K}\cdot\boldsymbol{\nabla}_{K}\left(\overline{b}\right) & =\frac{1}{m_{K}}\sum_{e\in\partial K}\left(\frac{b_{K_{e}}-b_{K}}{2}\right)\left(\frac{\boldsymbol{a}_{K_{e}}+\boldsymbol{a}_{K}}{2}\right)\cdot\boldsymbol{n}_{e,K}m_{e}\\
 & -\frac{1}{m_{K}}\sum_{e\in\partial K}\left(\frac{b_{K_{e}}-b_{K}}{2}\right)\left(\frac{\boldsymbol{a}_{K_{e}}-\boldsymbol{a}_{K}}{2}\right)\cdot\boldsymbol{n}_{e,K}m_{e}\,,
\end{align*}
and $b_{K}\,\boldsymbol{\nabla}_{K}\cdot\left(\overline{\boldsymbol{a}}\right)=(1/m_{K})\sum_{e\in\partial K}b_{K}\bigl((\boldsymbol{a}_{K_{e}}+\boldsymbol{a}_{K})/2\bigr)\cdot\boldsymbol{n}_{e,K}m_{e}$. 
\end{proof}

\subsection{Proof of the local energy balance \eqref{eq:acoustic_estim1}}

\label{subsec:proof_acoustic_estim1} Since the water depth remains
unchanged during the acoustic step, we write $h_{K}=h^{n}_{K}=h^{n+1}_{K}$.
With the discrete operators introduced above, the scheme reads

\begin{align}
\boldsymbol{u}^{n+1}_{K}-\boldsymbol{u}^{n}_{K} & =-\frac{\Delta\tau}{h_{K}}\boldsymbol{\nabla}_{K}\left(\overline{hp^{n+1}}\right),\label{eq:acoustic_u}\\
w^{n+1}_{K}-w^{n}_{K} & =\frac{3\Delta\tau}{2h_{K}}p^{n+1}_{K},\label{eq:acoustic_w}\\
p^{n+1}_{K}-p^{n}_{K} & =-a^{2}\Delta\tau\left[\frac{2}{h_{K}}\bigl(\theta_{p}w^{n+1}_{K}+(1-\theta_{p})w^{n}_{K}\bigr)+\boldsymbol{\nabla}_{K}\cdot(\boldsymbol{u}^{*})\right].\label{eq:acoustic_p}
\end{align}

Using the polarization identities, with the initial horizontal velocity
and the final vertical velocity and pressure, gives

\begin{align}
\frac{h_{K}}{2}\left(\|\boldsymbol{u}^{n+1}_{K}\|^{2}-\|\boldsymbol{u}^{n}_{K}\|^{2}\right)={} & \frac{h_{K}}{2}\|\boldsymbol{u}^{n+1}_{K}-\boldsymbol{u}^{n}_{K}\|^{2}\nonumber \\
 & -\Delta\tau\boldsymbol{u}^{n}_{K}\cdot\boldsymbol{\nabla}_{K}\left(\overline{hp^{n+1}}\right),\label{eq:estim_energy_K}\\
\frac{2h_{K}}{3}\left((w^{n+1}_{K})^{2}-(w^{n}_{K})^{2}\right)={} & -\frac{2h_{K}}{3}(w^{n+1}_{K}-w^{n}_{K})^{2}+2\Delta\tau p^{n+1}_{K}w^{n+1}_{K},\label{eq:estim_energy_W}\\
\frac{h_{K}}{2a^{2}}\left((p^{n+1}_{K})^{2}-(p^{n}_{K})^{2}\right)={} & -\frac{h_{K}}{2a^{2}}(p^{n+1}_{K}-p^{n}_{K})^{2}-2\Delta\tau p^{n+1}_{K}w^{n+1}_{K}\nonumber \\
 & -\Delta\tau h_{K}p^{n+1}_{K}\boldsymbol{\nabla}_{K}\cdot(\boldsymbol{u}^{*})\nonumber \\
 & +2\Delta\tau(1-\theta_{p})p^{n+1}_{K}(w^{n+1}_{K}-w^{n}_{K}).\label{eq:estim_energy_P}
\end{align}

Adding these identities cancels the terms $2\Delta\tau p^{n+1}_{K}w^{n+1}_{K}$.
The remaining vertical terms give $-\mathfrak{T}_{K}(\theta_{p})$,
with $\mathfrak{T}_{K}$ defined in \eqref{eq:W-upwind-residual}.

Moreover, the local duality identity \eqref{eq:duality4} and $\boldsymbol{u}^{*}_{e}=\overline{\boldsymbol{u}^{n}}_{e}-\boldsymbol{\Gamma}_{e}$
yield

\begin{align*}
 & \boldsymbol{u}^{n}_{K}\cdot\boldsymbol{\nabla}_{K}\left(\overline{hp^{n+1}}\right)+h_{K}p^{n+1}_{K}\boldsymbol{\nabla}_{K}\cdot(\boldsymbol{u}^{*})\\
 & \qquad=\frac{1}{m_{K}}\sum_{e\in\partial K}m_{e}\widetilde{(hp\boldsymbol{u})}_{e}\cdot\boldsymbol{n}_{e,K}-h_{K}p^{n+1}_{K}\boldsymbol{\nabla}_{K}\cdot(\boldsymbol{\Gamma}).
\end{align*}

Here $\widetilde{(hp\boldsymbol{u})}_{e}$ is the numerical energy
flux defined in the proposition. Since the height-dependent terms
in the energy remain constant, we obtain

\begin{align*}
 & (hE)^{n+1}_{K}-(hE)^{n}_{K}+\frac{\Delta\tau}{m_{K}}\sum_{e\in\partial K}m_{e}\widetilde{(hp\boldsymbol{u})}_{e}\cdot\boldsymbol{n}_{e,K}\\
 & \qquad=\frac{h_{K}}{2}\|\boldsymbol{u}^{n+1}_{K}-\boldsymbol{u}^{n}_{K}\|^{2}-\mathfrak{T}_{K}(\theta_{p})-\frac{h_{K}}{2a^{2}}(p^{n+1}_{K}-p^{n}_{K})^{2}\\
 & \qquad\quad+\Delta\tau h_{K}p^{n+1}_{K}\boldsymbol{\nabla}_{K}\cdot(\boldsymbol{\Gamma}),
\end{align*}

which proves the claimed identity.

\subsection{Proof of Proposition \ref{prop:acoustic_estim2}}

\label{subsec:proof_acoustic_estim2}In this proof, $h_{K}=h^{n}_{K}=h^{n+1}_{K}$,
$\widehat{p}_{K}=h_{K}p^{n+1}_{K}$ is the depth-weighted pressure
at the new time level and $\Delta p_{K}=p^{n+1}_{K}-p^{n}_{K}$, so
that $h_{K}p^{n}_{K}=\widehat{p}_{K}-h_{K}\Delta p_{K}$. We also
write $D_{K}=\boldsymbol{\nabla}_{K}\cdot(\overline{\mu},\boldsymbol{\delta}\widehat{p})$.
We start from the identity \eqref{eq:acoustic_estim1}, whose right-hand
side reads

\begin{equation}
\mathfrak{R}_{K}=\frac{h_{K}}{2}\|\boldsymbol{u}^{n+1}_{K}-\boldsymbol{u}^{n}_{K}\|^{2}-\mathfrak{T}_{K}(\theta_{p})-\frac{h_{K}}{2a^{2}}(\Delta p_{K})^{2}+\Delta\tau\,\widehat{p}_{K}\,\boldsymbol{\nabla}_{K}\cdot(\boldsymbol{\Gamma}).\label{eq:RK_residual}
\end{equation}

\paragraph{Step one: exchange and stabilization terms}

Using \eqref{eq:acoustic_w}, $\mathfrak{T}_{K}(\theta_{p})=(3\Delta\tau^{2}/(2h_{K}))(p^{n+1}_{K})^{2}(2\theta_{p}-1)\geq0$
for $\theta_{p}\geq1/2$.

Since $\boldsymbol{\Gamma}_{e}=\beta\Delta\tau\,\overline{\mu}_{e}\boldsymbol{\delta}_{e}(hp^{n})$
and $hp^{n}=\widehat{p}-h\Delta p$, the identities \eqref{eq:duality2},
with $a=\overline{\mu}$ and $b=c=\widehat{p}$, and \eqref{eq:duality3},
with $a=\overline{\mu}$, $b=h\Delta p$ and $c=\widehat{p}$, give

\[
\Delta\tau\,\widehat{p}_{K}\boldsymbol{\nabla}_{K}\cdot(\boldsymbol{\Gamma})=-\beta\Delta\tau^{2}\mathcal{S}_{K}\left(\overline{\mu},\|\boldsymbol{\delta}\widehat{p}\|^{2}\right)-\beta\Delta\tau^{2}h_{K}\Delta p_{K}\,D_{K}+\frac{\beta\Delta\tau^{2}}{m_{K}}\sum_{e\in\partial K}m_{e}\left(F^{(1)}_{K,e}-F^{(2)}_{K,e}\right),
\]

with the interface fluxes

\begin{equation}
F^{(1)}_{K,e}=\overline{\mu}_{e}\,\frac{\widehat{p}_{K}+\widehat{p}_{K_{e}}}{2}\,\frac{\widehat{p}_{K_{e}}-\widehat{p}_{K}}{2},\qquad F^{(2)}_{K,e}=\frac{\overline{\mu}_{e}}{2}\left(\widehat{p}_{K}\,h_{K_{e}}\Delta p_{K_{e}}-h_{K}\Delta p_{K}\,\widehat{p}_{K_{e}}\right).\label{eq:acoustic_flux_corrections}
\end{equation}

Moving the last sum to the left-hand side of \eqref{eq:acoustic_estim1}
defines the energy flux of Proposition \ref{prop:acoustic_estim2}:

\begin{equation}
\mathcal{G}_{a}(W_{K},W_{K_{e}},\boldsymbol{n}_{e,K})=\widetilde{(hp\boldsymbol{u})}_{e}\cdot\boldsymbol{n}_{e,K}-\beta\Delta\tau\left(F^{(1)}_{K,e}-F^{(2)}_{K,e}\right).\label{eq:Ga}
\end{equation}

\paragraph{Step two: exact local balance}

The pressure increment is kept and combined with the cross term:

\[
-\frac{h_{K}}{2a^{2}}(\Delta p_{K})^{2}-\beta\Delta\tau^{2}h_{K}\Delta p_{K}\,D_{K}=-\frac{h_{K}}{2a^{2}}\left(\Delta p_{K}+\beta a^{2}\Delta\tau^{2}D_{K}\right)^{2}+\frac{\beta^{2}a^{2}\Delta\tau^{4}h_{K}}{2}D^{2}_{K}.
\]

With $\boldsymbol{u}^{n+1}_{K}-\boldsymbol{u}^{n}_{K}=-\Delta\tau\boldsymbol{\nabla}_{K}(\overline{\widehat{p}})/h_{K}$,
gathering these terms gives

\begin{equation}
(hE)^{n+1}_{K}-(hE)^{n}_{K}+\frac{\Delta\tau}{m_{K}}\sum_{e\in\partial K}m_{e}\,\mathcal{G}_{a}=-\mathfrak{T}_{K}(\theta_{p})-\frac{h_{K}}{2a^{2}}\left(\Delta p_{K}+\beta a^{2}\Delta\tau^{2}D_{K}\right)^{2}+\frac{\Delta\tau^{2}}{2m_{K}}\mathcal{Q}_{K},\label{eq:acoustic_exact_local}
\end{equation}

with the quadratic form

\begin{equation}
\mathcal{Q}_{K}=\frac{m_{K}}{h_{K}}\left\Vert \boldsymbol{\nabla}_{K}(\overline{\widehat{p}})\right\Vert ^{2}+\beta^{2}a^{2}\Delta\tau^{2}m_{K}h_{K}\,D^{2}_{K}-2\beta\,m_{K}\,\mathcal{S}_{K}\left(\overline{\mu},\|\boldsymbol{\delta}\widehat{p}\|^{2}\right).\label{eq:acoustic_local_quadratic}
\end{equation}

It remains to show that $\mathcal{Q}_{K}\leq0$.

\paragraph{Step three: matrix form}

For a fixed cell, introduce the vector of weighted jumps $\boldsymbol{d}_{K}$,
the gradient matrix $\mathsf{B}_{K}$, with one column per edge, and
the summation vector $\boldsymbol{s}_{K}$:

\[
(\boldsymbol{d}_{K})_{e}=\sqrt{m_{e}\overline{\mu}_{e}}\,\llbracket\widehat{p}\rrbracket_{e},\qquad(\mathsf{B}_{K})_{e}=\sqrt{\frac{m_{e}}{m_{K}h_{K}\overline{\mu}_{e}}}\,\boldsymbol{n}_{e,K},\qquad(\boldsymbol{s}_{K})_{e}=\sqrt{\frac{h_{K}m_{e}\overline{\mu}_{e}}{m_{K}}}.
\]

By \eqref{eq:centered_gradient}, $\mathsf{B}_{K}\boldsymbol{d}_{K}=\sqrt{m_{K}/h_{K}}\,\boldsymbol{\nabla}_{K}(\overline{\widehat{p}})$;
moreover $\boldsymbol{s}^{T}_{K}\boldsymbol{d}_{K}=\sqrt{m_{K}h_{K}}\,D_{K}$
and $\|\boldsymbol{d}_{K}\|^{2}=m_{K}\mathcal{S}_{K}(\overline{\mu},\|\boldsymbol{\delta}\widehat{p}\|^{2})$.
Hence

\[
\mathcal{Q}_{K}=\|\mathsf{B}_{K}\boldsymbol{d}_{K}\|^{2}+\beta^{2}a^{2}\Delta\tau^{2}\left(\boldsymbol{s}^{T}_{K}\boldsymbol{d}_{K}\right)^{2}-2\beta\|\boldsymbol{d}_{K}\|^{2}.
\]

The definitions \eqref{eq:acoustic_geometry} give $\mathsf{B}_{K}\mathsf{B}^{T}_{K}=\mathsf{M}_{K}$
and $\|\boldsymbol{s}_{K}\|^{2}=\mathsf{T}_{K}$, and the closure
of the cell, $\sum_{e\in\partial K}m_{e}\boldsymbol{n}_{e,K}=0$,
gives

\[
\mathsf{B}_{K}\boldsymbol{s}_{K}=\frac{1}{m_{K}}\sum_{e\in\partial K}m_{e}\,\boldsymbol{n}_{e,K}=0.
\]

\paragraph{Step four: bound}

We decompose $\boldsymbol{d}_{K}$ into its orthogonal projections
onto $\operatorname{span}(\boldsymbol{s}_{K})$ and its orthogonal
complement:

\[
\boldsymbol{d}_{\parallel}:=\frac{\boldsymbol{s}^{T}_{K}\boldsymbol{d}_{K}}{\|\boldsymbol{s}_{K}\|^{2}}\boldsymbol{s}_{K},\qquad\boldsymbol{d}_{\perp}:=\boldsymbol{d}_{K}-\boldsymbol{d}_{\parallel}.
\]

Thus $\boldsymbol{s}^{T}_{K}\boldsymbol{d}_{\perp}=0$. Since $\mathsf{B}_{K}\boldsymbol{s}_{K}=0$,
the gradient only sees the orthogonal part, $\mathsf{B}_{K}\boldsymbol{d}_{K}=\mathsf{B}_{K}\boldsymbol{d}_{\perp}$,
while the divergence only sees the collinear one: $\boldsymbol{s}^{T}_{K}\boldsymbol{d}_{K}=\boldsymbol{s}^{T}_{K}\boldsymbol{d}_{\parallel}$,
hence $(\boldsymbol{s}^{T}_{K}\boldsymbol{d}_{K})^{2}=\|\boldsymbol{s}_{K}\|^{2}\|\boldsymbol{d}_{\parallel}\|^{2}=\mathsf{T}_{K}\|\boldsymbol{d}_{\parallel}\|^{2}$.
Since $\mathsf{B}_{K}\mathsf{B}^{T}_{K}=\mathsf{M}_{K}$, the operator
norm of $\mathsf{B}_{K}$ satisfies $\|\mathsf{B}_{K}\|^{2}_{2}=\lambda_{\max}(\mathsf{M}_{K})$.
Hence,

\[
\|\mathsf{B}_{K}\boldsymbol{d}_{\perp}\|^{2}\leq\lambda_{\max}(\mathsf{M}_{K})\|\boldsymbol{d}_{\perp}\|^{2}.
\]

With $\|\boldsymbol{d}_{K}\|^{2}=\|\boldsymbol{d}_{\parallel}\|^{2}+\|\boldsymbol{d}_{\perp}\|^{2}$,
we obtain

\[
\mathcal{Q}_{K}\leq\left(\lambda_{\max}(\mathsf{M}_{K})-2\beta\right)\|\boldsymbol{d}_{\perp}\|^{2}+\beta\left(\beta a^{2}\Delta\tau^{2}\mathsf{T}_{K}-2\right)\|\boldsymbol{d}_{\parallel}\|^{2}\leq0
\]

under \eqref{eq:CFL_acoustic}. The first condition thus controls
the jumps that carry a gradient, and the second one the mode in which
the depth-weighted pressure in $K$ differs by the same amount from
that in every neighbouring cell, a local analogue of the checkerboard
pattern on a Cartesian mesh. Both bounds are attained: by $\boldsymbol{d}_{K}=\boldsymbol{s}_{K}$
for the second one, and for the first one by an eigenvector of $\mathsf{B}^{T}_{K}\mathsf{B}_{K}$
associated with $\lambda_{\max}(\mathsf{M}_{K})$, which lies in the
range of $\mathsf{B}^{T}_{K}$ and is therefore orthogonal to $\boldsymbol{s}_{K}$.
Since the jumps on the edges of $K$ are independent, \eqref{eq:CFL_acoustic}
are exactly the conditions ensuring $\mathcal{Q}_{K}\leq0$ for arbitrary
local values of $\widehat{p}$, with the depths fixed. The remaining
terms of \eqref{eq:acoustic_exact_local} are non-positive, which
proves \eqref{eq:acoustic_estim2}. Finally, $F^{(1)}$ and $F^{(2)}$
are fluxes of the form $\mathcal{F}^{2}$ and $\mathcal{F}^{3}$ of
Lemma \ref{lem:duality}: they are conservative and vanish for equal
states, so that $\mathcal{G}_{a}$ is conservative and consistent
with $hp\boldsymbol{u}\cdot\boldsymbol{n}$.

\section{Linear waves and dispersion relation}\label{sec:linear_waves}

The two-dimensional system of equations \eqref{eq:model} is linearized
around the equilibrium state $\left(H,0,0,0,0\right)$:

\begin{subequations}
\begin{equation}
{\displaystyle \frac{\partial h}{\partial t}+H\frac{\partial u}{\partial x}+H\frac{\partial v}{\partial y}=0}\,,
\end{equation}

\begin{equation}
{\displaystyle \frac{\partial u}{\partial t}+g\frac{\partial h}{\partial x}+\frac{\partial p}{\partial x}=0}\,,
\end{equation}

\begin{equation}
{\displaystyle \frac{\partial v}{\partial t}+g\frac{\partial h}{\partial y}+\frac{\partial p}{\partial y}=0}\,,
\end{equation}

\begin{equation}
H{\displaystyle \frac{\partial w}{\partial t}=\frac{3}{2}p}\,,
\end{equation}

\begin{equation}
{\displaystyle H\frac{\partial p}{\partial t}=-a^{2}\left(2w+H\frac{\partial u}{\partial x}+H\frac{\partial v}{\partial y}\right)}\,,
\end{equation}
\end{subequations}
 where $\left(h,u,v,w,p\right)$ is now the set of small perturbation
variables around the equilibrium state for simplifying notations.
Now, if we seek solutions of the form $\left(h,u,v,w,p\right)=\left(A_{h},A_{u},A_{v},A_{w},A_{p}\right)e^{i\left(\boldsymbol{k}\cdot\boldsymbol{x}-\omega t\right)}$,
the dispersion relation is given by the roots of the following polynomial
in $\omega^{2}$, which give the two propagating branches:

\begin{equation}
{\displaystyle \frac{H^{2}}{3a^{2}}\omega^{4}-\left(1+\frac{\Vert\boldsymbol{k}\Vert^{2}H^{2}}{3}\left(1+\frac{gH}{a^{2}}\right)\right)\omega^{2}+gH\Vert\boldsymbol{k}\Vert^{2}=0}\,.
\end{equation}

The stationary transverse mode does not appear in this polynomial.
Note that in the limit $a\to+\infty$ we recover the dispersion relation
emanating from the original SGN equations, that is

\begin{equation}
\left(1+\frac{k^{2}H^{2}}{3}\right)\omega^{2}=gHk^{2}\,.
\end{equation}

\begin{figure}[htbp]
\centering{}\includegraphics[width=0.95\textwidth]{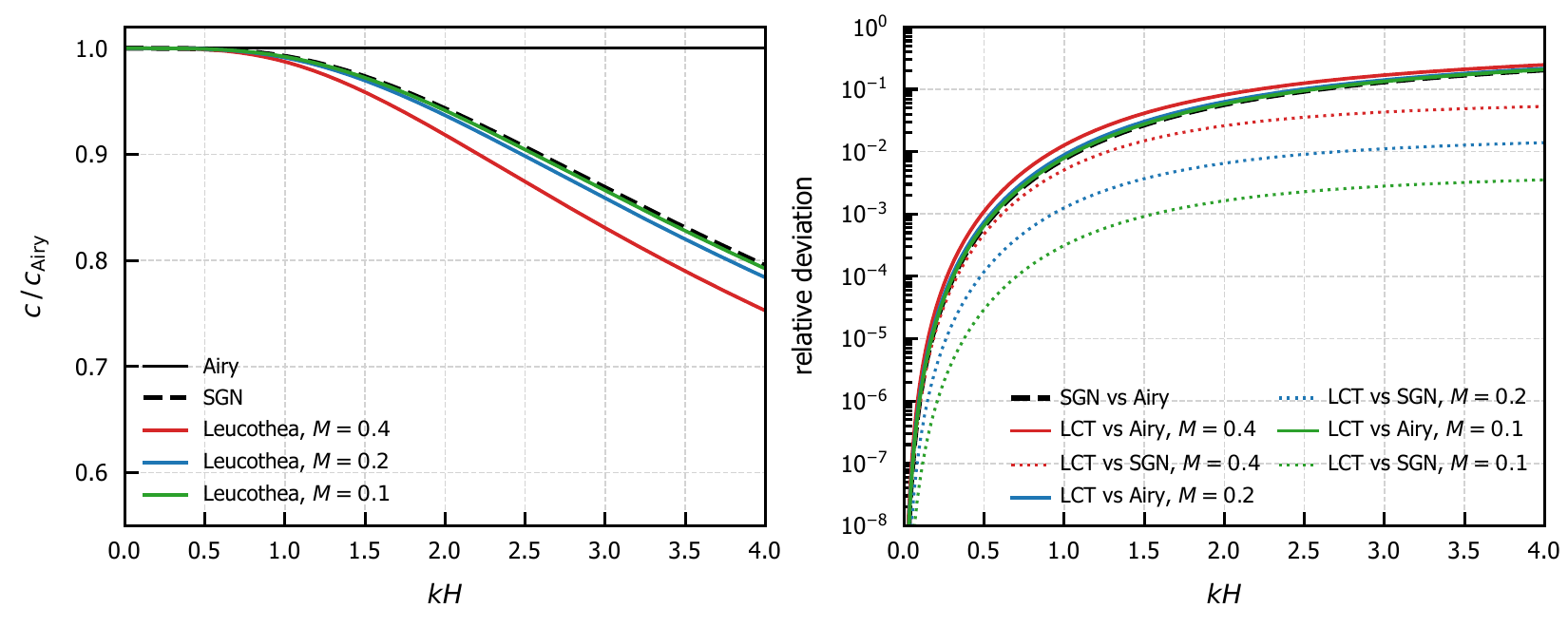}\caption{Dispersion properties of the Leucothea model: (left) normalized phase
speed $c/c_{\mathrm{Airy}}$ of the gravity branch as a function of
$kH$ for $M=0.4,\,0.2,\,0.1$, compared with the SGN dispersion relation;
(right) relative deviations with respect to the Airy reference (solid
lines, SGN dashed) and with respect to the SGN phase speed (dotted
lines). The dispersion error with respect to Airy is essentially that
of the SGN family, while the relaxation error remains uniformly below
$M^{2}/2$. }\label{fig:dispersion-M}
\end{figure}

Denoting $\kappa=kH$ the shallowness parameter and $M=\sqrt{gH}/a$
the Mach number, Fig.~\ref{fig:dispersion-M} compares the phase
speed of the gravity branch, normalized by the Airy reference $c_{\mathrm{Airy}}=\sqrt{g\tanh(kH)/k}$,
with the SGN one, for the values of the Mach number used in practice.
Two distinct errors must be read on this figure. The deviation from
the Airy reference measures the intrinsic dispersion error of the
SGN family, a first-order approximation in the shallowness parameter.
The deviation of the Leucothea model from its SGN target, that is
the relaxation error induced by the finite acoustic speed, is bounded
by $1-(1+M^{2})^{-1/2}<M^{2}/2$ over the whole spectrum, the supremum
being reached in the short-wave limit (1.94\% for $M=0.2$, 7.15\%
for $M=0.4$). Indeed, with $x=(kH)^{2}/3$, $m=M^{2}$ and $z=c^{2}_{\mathrm{lct}}/c^{2}_{\mathrm{SGN}}$,
the gravity branch corresponds to the smallest root of 
\[
P(z)=\frac{mx}{(1+x)^{2}}\,z^{2}-\left(1+\frac{mx}{1+x}\right)z+1\,.
\]
Since $P(1)=-mx^{2}/(1+x)^{2}\leq0$ and $P(1/(1+m))=m/((1+m)(1+x))+mx/((1+x)^{2}(1+m)^{2})\geq0$,
this root lies in $[1/(1+m),1]$, and it tends to $1/(1+m)$ as $kH\to\infty$,
which gives the bound and its optimality. For $M=0.2$, this relaxation
error is more than one order of magnitude smaller than the SGN dispersion
error at the largest $kH$ shown. For $M=0.4$, it remains smaller,
although the separation is less pronounced.

\section{Periodic and solitary travelling waves}\label{sec:cnoidal}

We look for the stationary solutions in the reference frame $\xi=x-ct$,
where $c$ is a constant velocity, to find the so-called cnoidal waves,
the periodic travelling wave solutions for the one-dimensional version
of the model \eqref{eq:model} over a flat bottom, written here in
primitive-variable quasilinear form:

\begin{subequations}
\begin{equation}
{\displaystyle \partial_{t}h+\partial_{x}\left(hu\right)=0}\,,
\end{equation}

\begin{equation}
{\displaystyle h\partial_{t}u+hu\partial_{x}u+\left(gh+p\right)\partial_{x}h+h\partial_{x}p=0}\,,
\end{equation}

\begin{equation}
{\displaystyle h\partial_{t}w+hu\partial_{x}w=\frac{3p}{2}}\,,
\end{equation}

\begin{equation}
{\displaystyle h\partial_{t}p+hu\partial_{x}p=-a^{2}\left(2w+h\partial_{x}u\right)}\,.
\end{equation}
\end{subequations}

Given that $\partial/\partial t=-c\,\mathrm{d}/\mathrm{d}\xi=-c\left(.\right)'$
and $\partial/\partial_{x}=d/d\xi=\left(.\right)'$, the following
expressions are derived from the above system:

\begin{subequations}
\begin{equation}
{\displaystyle \left[h\left(u-c\right)\right]'=0}\,,
\end{equation}

\begin{equation}
{\displaystyle \left[h\left(u-c\right)u+g\frac{h^{2}}{2}+hp\right]'=0}\,,
\end{equation}

\begin{equation}
{\displaystyle \left[h\left(u-c\right)w\right]'=\frac{3p}{2}}\,,
\end{equation}

\begin{equation}
{\displaystyle \left[h\left(u-c\right)p\right]'=-a^{2}\left(2w+hu'\right)}\,.
\end{equation}
\end{subequations}

The integration of the first two equations of the above system by
introducing two integration constants $\mathcal{D}$ and $\mathcal{P}$
yields

\begin{subequations}
\begin{equation}
{\displaystyle h\left(u-c\right)=\mathcal{D}}\,,
\end{equation}

\begin{equation}
{\displaystyle \mathcal{D}u+g\frac{h^{2}}{2}+hp=\mathcal{P}+c\mathcal{D}\Rightarrow p=-\frac{gh}{2}+\frac{\mathcal{P}}{h}-\frac{\mathcal{D}^{2}}{h^{2}}:=P(h)}\,,
\end{equation}

\begin{equation}
{\displaystyle \mathcal{D}w'=\frac{3p}{2}}\,,
\end{equation}

\begin{equation}
{\displaystyle \mathcal{D}p'=-a^{2}\left(2w+hu'\right)\Rightarrow\mathcal{D}P'(h)h'=-a^{2}\left(2w-\mathcal{D}\frac{h'}{h}\right)}\,,
\end{equation}
\end{subequations}
 resulting in the final system of nonlinear ordinary differential
equations to solve with $h$ and $w$ as unknowns:

\begin{subequations}
\label{eq:lct-ode-cnoidal} 
\begin{equation}
2\mathcal{D}\:w'=3P(h)\,,
\end{equation}

\begin{equation}
{\displaystyle \frac{\mathcal{D}}{2}\left(\frac{1}{h}-\frac{P'(h)}{a^{2}}\right)\:h'=w}\,,
\end{equation}
\end{subequations}

with $P(h):=-gh/2+\mathcal{P}/h-\mathcal{D}^{2}/h^{2}$. This system
admits an additional integral of motion:

\begin{equation}
{\displaystyle \frac{3P(h)}{2}\left(\frac{1}{h}-\frac{P'(h)}{a^{2}}\right)h'=2ww'}\,.
\end{equation}

Then we deduce that there exists another constant of integration $\mathcal{R}$
such that

\begin{equation}
w^{2}+\frac{3}{4a^{2}}P^{2}(h)+\frac{3gh}{4}+\frac{3\mathcal{P}}{2h}-\frac{3\mathcal{D}^{2}}{4h^{2}}-\frac{3\mathcal{R}}{2}=0\,.
\end{equation}

Replacing $w$ by its expression as a function of $h'$, we finally
find the nonlinear ordinary differential equation satisfied by the
water depth $h$ 

\begin{equation}
\left(h'\right)^{2}+\frac{3g}{\mathcal{D}^{2}}\frac{{\displaystyle h^{3}-\frac{2\mathcal{R}}{g}h^{2}+\frac{2\mathcal{P}}{g}h-\frac{\mathcal{D}^{2}}{g}+\frac{h^{2}}{a^{2}g}P^{2}(h)}}{{\displaystyle \left(1-\frac{h}{a^{2}}P'(h)\right)^{2}}}=0\,.
\end{equation}

In the limit $a\to+\infty$, this equation is the same as that found
for the Green--Naghdi equations, for which the exact analytical solution
is well known \citep{carter_kinematics_2011}. Considering three water
depths $0<h_{1}\leq h_{2}\leq h_{3}$, the equation can be rewritten
as

\begin{equation}
{\displaystyle (h')^{2}+\frac{3g}{\mathcal{D}^{2}}(h-h_{1})(h-h_{2})(h-h_{3})=0,}
\end{equation}

with expressions

\begin{equation}
{\displaystyle h_{1}h_{2}h_{3}=\frac{\mathcal{D}^{2}}{g},\quad h_{1}+h_{2}+h_{3}=\frac{2\mathcal{R}}{g},\quad h_{1}h_{2}+h_{2}h_{3}+h_{3}h_{1}=\frac{2\mathcal{P}}{g}.}\label{eq:cnoidal_cst_int}
\end{equation}

The solution $h$ is then parametrized by $(h_{i})_{i=1,2,3}$ and
is given by the explicit formula 
\begin{equation}
{\displaystyle h(\xi)=h_{2}+(h_{3}-h_{2})\:{\rm cn}^{2}\left(\sqrt{\frac{3(h_{3}-h_{1})}{4h_{1}h_{2}h_{3}}}\xi,m\right),}\label{eq:h-cn}
\end{equation}
where ${\rm cn}$ is a Jacobi elliptic function with parameter $m=(h_{3}-h_{2})/(h_{3}-h_{1})\in[0,1]$.
In order to relate the parameters $h_{1},h_{2},h_{3}$ with classical
parameters of water wave theory, we introduce the wavenumber $k$,
the amplitude $A$ and the mean water depth $H$, given by the properties
of the Jacobi elliptic functions:

\begin{subequations}
\begin{equation}
k=\frac{\pi}{K(m)}\sqrt{\frac{3(h_{3}-h_{1})}{4h_{1}h_{2}h_{3}}}\,,
\end{equation}

\begin{equation}
A=h_{3}-h_{2}=m(h_{3}-h_{1})\,,
\end{equation}

\begin{equation}
H=h_{1}+(h_{3}-h_{1})\:\frac{E(m)}{K(m)}\,,
\end{equation}
\end{subequations}

where $K(m)$ and $E(m)$ are respectively the complete elliptic integrals
of the first and second kind. We are now in a position to introduce
the two dimensionless parameters usually found in the theory of water
waves. We first set $\varepsilon=A/H$ the nonlinearity parameter
and $\mu=k^{2}H^{2}$ the dispersion parameter (or shallowness parameter).
After a few calculations, we can find that 
\begin{subequations}
\begin{equation}
{\displaystyle \frac{h_{3}}{h_{1}}-1=\frac{\varepsilon K(m)}{mK(m)-\varepsilon E(m)}:=\frac{1}{f(m,\varepsilon)},}\label{eq:u3}
\end{equation}
\begin{equation}
\mu=\frac{3\pi^{2}}{4K^{4}(m)}\frac{\left(K(m)f(m,\varepsilon)+E(m)\right)^{2}}{f(m,\varepsilon)^{3}+(2-m)f(m,\varepsilon)^{2}+(1-m)f(m,\varepsilon)}\,.
\end{equation}
\end{subequations}

Then, $m$ must satisfy the equation 
\begin{equation}
{\displaystyle G(m,\varepsilon,\mu)=f(m,\varepsilon)^{3}+(2-m)f(m,\varepsilon)^{2}+(1-m)f(m,\varepsilon)-\frac{3\pi^{2}}{4\mu\:K(m)^{4}}\left(K(m)\:f(m,\varepsilon)+E(m)\right)^{2}=0.}\label{eq:gm}
\end{equation}
Once the solution $m$ of $G(m,\varepsilon,\mu)=0$ has been found
numerically, we can easily deduce $(h_{i})_{i=1,2,3}$ from previous
relations: 
\begin{equation}
{\displaystyle h_{1}=\frac{H\:K(m)\:f(m,\varepsilon)}{K(m)\:f(m,\varepsilon)+E(m)},\quad h_{3}=h_{1}\left(1+\frac{1}{f(m,\varepsilon)}\right),\quad h_{2}=mh_{1}+(1-m)h_{3}.}\label{eq:eta_i}
\end{equation}

Coming back to the numerical resolution of the system \eqref{eq:lct-ode-cnoidal},
we introduce a new variable $\tilde{w}:=(2\mathcal{D}/g)\,w$ and
define $Q(h):=P(h)/g$ to solve a new system:

\begin{subequations}
\label{eq:lct-ode-cnoidal-reform} 
\begin{equation}
\tilde{w}'=3Q(h)\,,
\end{equation}

\begin{equation}
{\displaystyle \left(\frac{1}{h}-\frac{M^{2}}{H}Q'(h)\right)\:h'=\frac{g}{\mathcal{D}^{2}}\tilde{w}}\,,
\end{equation}
\end{subequations}

where the Mach number $M=\sqrt{gH}/a$ is defined using the mean water
depth $H$. Such a formulation allows us to express the constants
directly from \eqref{eq:cnoidal_cst_int}. It is easily seen that
we can restrict our attention to periodic solutions with $h$ an even
function and $\tilde{w}$ an odd function. Selecting a mean water
depth $H$, a wavelength $L$, an amplitude $A$ and a Mach number
$M$, we numerically seek the root of

\begin{equation}
\mathcal{F}\left(h_{0},\mathcal{D},\mathcal{P}\right)=\left(\begin{array}{l}
\tilde{w}(L/2)\\
h_{0}-h\left(L/2\right)-A\\
{\displaystyle \frac{2}{L}\int^{L/2}_{0}h\left(x\right)dx}-H
\end{array}\right)\,,
\end{equation}

starting from the values $h_{0}=h_{3}$, $\mathcal{D}$ and $\mathcal{P}$
found in the Serre--Green--Naghdi limit when $a\to+\infty$ and
where $\left(h,\tilde{w}\right)$ is the solution of the numerical
integration of \eqref{eq:lct-ode-cnoidal-reform} over the interval
$[0,L/2]$ with the initial condition $(h(0),\tilde{w}(0))=(h_{0},0)$.
The root of $\mathcal{F}$ yields the parametrized periodic travelling
wave, its three components ensuring respectively the prescribed wavelength
$L$, amplitude $A$ and mean water depth $H$ for a selected Mach
number $M$. In practice, we use continuation in the Mach number to
facilitate convergence of the Newton iterations: the analytical Serre--Green--Naghdi
solution provides the exact starting point at $M=0$, and the target
value is reached through successive resolutions with an adaptive step,
halved whenever the Newton iterations fail and enlarged after fast
convergence. With this strategy, the procedure converged systematically
over the whole range of parameters considered in this work. The ODEs
are integrated using DOP853 with relative and absolute tolerances
of $10^{-13}$ and $10^{-14}$, respectively; the Newton iterations
stop when the Euclidean norm of the residual falls below $10^{-8}$.
The reference frame is fixed by imposing zero mean discharge over
one wavelength. Since $h(u-c)=\mathcal{D}$ and the mean depth is
$H$, this gives $c=-\mathcal{D}/H$, and the horizontal velocity
is reconstructed as $u=c+\mathcal{D}/h$. Finally, let us mention
that for solitary waves the shooting strategy becomes intrinsically
ill-conditioned, the sensitivities behaving exponentially near the
saddle point $h=H$, and is replaced by a global collocation method
solving the boundary-value problem on $[0,L/2]$ with $(\mathcal{P},\mathcal{D})$
as unknown parameters, combined with the same continuation in $M$.

\begin{figure}[htbp]
\begin{centering}
\includegraphics[width=0.6\textwidth]{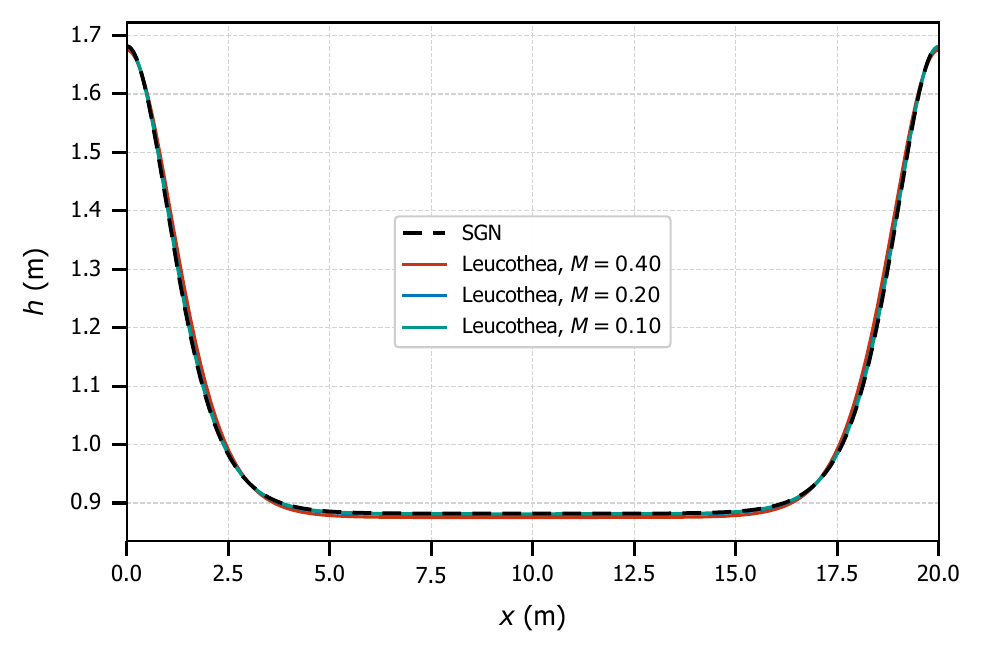}
\par\end{centering}
\caption{Comparison of the parametrized periodic travelling wave solutions
of the Serre--Green--Naghdi model (dashed) and of the Leucothea
model at $M=0.4,\,0.2,\,0.1$, with $H=1\ \mathrm{m}$, $L=20\ \mathrm{m}$
and $A=0.8\ \mathrm{m}$.}\label{fig:cnoidal-vs-sgn}
\end{figure}

To illustrate the discrepancy between the Serre--Green--Naghdi model
and the Leucothea model, Fig.~\ref{fig:cnoidal-vs-sgn} shows the
solution for each model. Very little deviation is observed, even at
fairly large Mach numbers, confirming the robustness of the Leucothea
model with respect to the Mach number, consistently with the dispersion
analysis of \ref{sec:linear_waves}.

\bibliographystyle{plainnat}
\phantomsection\addcontentsline{toc}{section}{\refname}\bibliography{main}

\end{document}